\documentstyle[amsfonts,amssymb]{article}    

\newtheorem{theorem}{Theorem}[section]
\newtheorem{lemma}[theorem]{Lemma} 
 
\newtheorem{corollary}[theorem]{Corollary} 
\newtheorem{remark}[theorem]{Remark} 
\newtheorem{proposition}[theorem]{Proposition}

\title{\bf Fundamental Solutions for the Klein-Gordon Equation \\ in  de~Sitter Spacetime
}

\author{{\bf Karen Yagdjian, 
Anahit Galstian}
 }

\begin{document}

\date{}
\maketitle
\thispagestyle{empty} 
\vspace{-0.3cm}

\begin{center}
{\it Department of Mathematics, 
University of Texas-Pan American,  
1201 W.~University Drive, \\
Edinburg, TX 78541-2999,  
USA }
\end{center}
\medskip

\addtocounter{section}{-1}
\renewcommand{\theequation}{\thesection.\arabic{equation}}
\setcounter{equation}{0}
\pagenumbering{arabic}
\setcounter{page}{1}
\thispagestyle{empty}

\hspace{2cm}\begin{abstract}
\begin{small}
 In this article we construct the fundamental solutions  for the Klein-Gordon  equation  in 
de~Sitter spacetime.  We use these fundamental solutions  to  represent     solutions of the Cauchy problem and to prove   $L^p-L^q$ estimates 
for the solutions of the equation with and without a source term. 
\medskip



\end{small}
\end{abstract}

\section{Introduction and Statement of Results}
\label{S1}

\setcounter{equation}{0}
\renewcommand{\theequation}{\thesection.\arabic{equation}}

In this paper we construct the fundamental solutions  for the Klein-Gordon equation in the de~Sitter 
 spacetime and use these fundamental solutions  to find representations of the solutions 
to the Cauchy problem as well as   $L^p-L^q$ estimates for them.
\smallskip

After averaging on a suitable scale, our universe is homogeneous 
and isotropic; therefore, the properties of the universe can be properly described by treating the matter as a perfect homogeneous fluid. 
In the models of the universe proposed by Einstein \cite{Einstein} and de~Sitter \cite{{Sitter}} the line
element is  connected with the proper 
mass density  and the proper pressure in the universe by the field equations
for a perfect fluid. There are two alternatives  
which lead to the solutions of Einstein    and de~Sitter, respectively \cite[Sec.132]{Moller}. 
\smallskip

In the models proposed by Einstein \cite{Einstein} and de Sitter \cite{{Sitter}} the universe is  assumed to be a {\sl static} 
system, which means that we can introduce a system of coordinates $x^i=(r,\theta ,\phi ,ct)$ in which the line
element has the static and spherically symmetric form
\[
ds^2=-b(r)c^2dt^2 + a(r)dr^2 + r^2(d \theta ^2+ \sin^2 d\phi ^2),
\]
where $a$ and $b$ are functions of $r$ only. On account of the assumed homogeneity of the universe any point in the space may be taken 
as the origin $r=0$ of the spatial system of coordinates.
The functions $a(r)$ and $b(r)$ are  connected with the proper 
mass density $ \mu $ and the proper pressure $ p$ in the universe by the field equations
for a perfect fluid
\begin{eqnarray}
\label{0.1} 
&  &
( \mu  c^2 + p )b'=0, \\
\label{0.2} 
&  &
\frac{b' }{abr} -\frac{1}{r^2} \left(1-\frac{1}{a} \right)+\Lambda  =\kappa  p, \\
\label{0.3} 
&  &
\frac{a' }{a^2r} +\frac{1}{r^2} \left(1-\frac{1}{a} \right)-\Lambda  =\kappa \mu  c^2,
\end{eqnarray}
where $\Lambda $ is {\it cosmological constant}, while $p$ and $\mu $ are constants. The general solution of the equation (\ref{0.3}) is 
\begin{equation}
a(r)= \left( 1-\frac{2M_{bh}}{r} - \frac{\Lambda r^2}{3}\right)^{-1}.
\end{equation}
The constant of integration $M_{bh}$ may have a meaning of the ``mass of the black holes''. 
 There are two alternatives  $b'=0 $ or $ \mu c^2+ p=0$,
which lead to the solutions of Einstein   and de Sitter, respectively. 
In the case of de Sitter universe $ab=constant$. 
By a trivial change of scale of the time variable, this constant can, of course, always be 
made equal to 1, which means 
\[
b(r)= 1-\frac{2M_{bh}}{r} - \frac{\Lambda r^2}{3}\,.
\]
The corresponding metric with the line element
\[
ds^2= - \left( 1-\frac{2M_{bh}}{r} - \frac{\Lambda r^2}{3}\right) c^2\, dt^2
+ \left( 1-\frac{2M_{bh}}{r} - \frac{\Lambda r^2}{3}\right)^{-1}dr^2 + r^2(d\theta ^2 + \sin^2 \theta \, d\phi ^2)
\]
is called the Schwarzschild - de~Sitter metric. The Cauchy problem for the linear wave equation without source term   on  
the maximally extended 
Schwarzschild - de Sitter spacetime in the case of non-extremal black-hole corresponding to parameter values
\[ 
0<M_{bh}< \frac{1}{3\sqrt{\Lambda }}, 
\]
is considered  by Dafermos and Rodnianski \cite{Dafermos-Rodnianski_1}.  
They   proved that in the region bounded by a set of black/white hole horizons and cosmological horizons, 
solutions converge pointwise to a constant faster than any given polynomial rate, 
where the decay is measured with respect to natural future-directed advanced and retarded time coordinates.
\smallskip

There is an important question of local energy decay for the solution of the wave equation and Klein-Gordon equation
in black hole spacetime.
Results on the decay of local energy can provide a proof of the global nonlinear stability of the spacetime.
The global nonlinear stability we believe is only known for  Minkowski spacetime. 
Bony and Hafner \cite{Bony-Hafner}  describe an expansion of the solution of the 
wave equation in the Schwarzschild - de~Sitter metric in terms of resonances. 
The resonances correspond to the frequencies and rates of dumping of signals emitted by the black hole in the presence of 
perturbations (see \cite[Chapter~4.35]{Chandrasekhar}).
The main term in the expansion obtained in \cite{Bony-Hafner} is due to a zero resonance. 
The error term decays polynomially if one permits a logarithmic derivative loss in the angular directions and exponentially if one permits a small derivative loss in the angular directions. 
\smallskip

In the present paper we set $M_{bh}=0$ to exclude black holes. The case of the Klein-Gordon equation in the presence  of the black holes
will be discussed in the forthcoming paper.  
Thus,  the de Sitter line element has  the form
\[
ds^2= - \left( 1- \frac{r^2}{R^2}\right) c^2\, dt^2+ \left( 1- \frac{r^2}{R^2}\right)^{-1}dr^2 + r^2(d\theta ^2 + \sin^2 \theta \, d\phi ^2)\,.
\]
The transformation 
\begin{eqnarray*}
r'=\frac{r}{\sqrt{1-r^2/R^2}} e^{-ct/R}\,, \quad t'=t+\frac{R}{2c} \ln \left( 1- \frac{r^2}{R^2}\right) \,, \quad \theta '=\theta \,, 
\quad \phi '=\phi 
\end{eqnarray*}
leads to the following form for the line element:
\begin{equation}
\label{M101}
ds^2= -   c^2\, d{t'}^2+e^{2ct'/R}( d{r'}^2   + r'^2\,d{\theta'} ^2 + r'^2\sin^2 \theta' \, d{\phi '}^2)\,.
\end{equation}
Finally, defining new space coordinates $x'$, $y'$, $z'$ connected with $r'$, $\theta '$, $\phi '$
by the usual equations connecting Cartesian  coordinates and polar coordinates in a Euclidean space, (\ref{M101}) may be written 
\cite[Sec.134]{Moller}
\begin{equation}
\label{M102}
ds^2= -   c^2\, d{t'}^2+ e^{2ct'/R}( d{x'}^2   + d{y'} ^2 +  d{z '}^2)\,.
\end{equation}
The new coordinates $x'$, $y '$, $z '$, $t'$ can take all values from $-\infty$ to $\infty$.
Here $R$ is the ``radius'' of the universe. 
The de~Sitter model allows us to get an explanation of the actual red shift 
of spectral lines observed by Hubble and Humanson \cite{Moller}.  The de~Sitter model  also enjoys the advantage of being 
the only time-dependent cosmological model for which both particle 
creation and vacuum stress have been explicitly evaluated by all known techniques \cite{Birrell}. 
In a certain sense all  solutions  look like the de~Sitter solution at late times \cite{Heinzle-Rendall}.  
\medskip

The homogeneous and isotropic cosmological models possess highest symmetry that makes them  more amenable to rigorous study.
Among them we mention FLRW (Friedmann-Lematre-Robertson-Walker) models. The simplest class of cosmological models 
can be obtained if we assume additionally that the metric of the slices of constant time is flat and that 
the spacetime metric can be written in the form
\[
ds^2= -dt^2+ a^2(t)( d{x}^2   + d{y} ^2 +  d{z}^2 )
\]
with an appropriate scale factor $a(t)$. Although on the made  assumptions,   the  spatially flat FLRW models appear to give a good 
explanation of our universe. The assumption that the universe is expanding leads to the positivity of  
  the time derivative $\frac{d }{dt}a (t)$. A further assumption that the universe obeys the accelerated expansion suggests 
that the second derivative  $\frac{d^2 }{dt^2}a (t)$  is positive.  
A substantial amount of the observational material can be satisfactorily interpreted in terms of the models 
which take into account existing acceleration of   the recession of distant galaxies.
\smallskip

The time dependence of the function $a(t)$ is  determined by  the Einstein field equations for gravity. 
The Einstein equations with the {\it cosmological constant} $\Lambda $ have form
\[
R_{\mu \nu }-\frac{1}{2} g_{\mu \nu }R = -8\pi GT_{\mu \nu }-\Lambda g_{\mu \nu },
\] 
where term $ \Lambda g_{\mu \nu }$ can be interpreted as an energy-momentum of the vacuum. 
 Even a small value of 
$\Lambda $ could have drastic effects on the evolution of the universe. 
Under the assumption of FLRW symmetry the equation of motion in the case of positive
cosmological constant \, $\Lambda $ \, leads to the solution
\[
a(t)=a(0)e^{\sqrt{\frac{\Lambda }{3}}t},
\]
which produces models with exponentially accelerated expansion. The model described by the last equation is usually called the {\it de Sitter
model}.
\smallskip

The unknown of principal importance in the Einstein equations is 
a metric $g$. It comprises the basic geometrical feature of the gravitational field, 
and consequently explains the phenomenon of the mutual gravitational attraction of substance. 
In the presence of matter these equations contain a non-vanishing right hand side $-8\pi GT_{\mu \nu }$. 
In general the matter fields described by the function  $\phi $  must satisfy  equations of motion and 
in the case of the massive scalar field the equation of motion is that $\phi $ 
should satisfy the  Klein-Gordon equation generated by the metric $g$. 
In the  de~Sitter universe  the equation for the scalar field with mass \,  $m$\,   and potential function \, $V$   \, 
written out explicitly in coordinates is (See, e.g. \cite{Friedrich-Rendall,Rendall_2004}.)
\begin{equation}
\label{1.3}
  \phi_{tt} +   n H \phi_t - e^{-2Ht} \bigtriangleup \phi + m^2\phi=   - V'(\phi )\,.
\end{equation}
Here $x \in {\mathbb R}^n$, $t \in {\mathbb R}$, and $\bigtriangleup $ is 
 the Laplace operator on the flat metric, $\bigtriangleup := \sum_{j=1}^n \frac{\partial^2 }{\partial x_j ^2} $, 
while $H = \sqrt{\Lambda /3}$ is Hubble constant.  If we introduce new unknown function $u = e^{-\frac{n}{2}Ht}\phi$, then the semilinear 
Klein-Gordon  equation for $u$ on de~Sitter spacetime takes the form 
\begin{equation}
\label{K_G_semilinear}
u_{tt} - e^{-2Ht} \bigtriangleup u  + M^2 u=  - e^{\frac{n}{2}t}V'(e^{-\frac{n}{2}t}u ),
\end{equation}
where
$M^2:= m^2 - \frac{n^2}{4}H^2$.  
Henceforth the quantity $M$, with nonnegative real part $\Re M \geq 0$, defined by 
\[
M^2:= m^2 - \frac{n^2}{4}H^2 \,,
\] 
will be called the ``curved mass'' of particle. 
We extract a linear part of the   equation (\ref{K_G_semilinear}) as an initial model that must be treated
first:
\begin{equation}
\label{K_G_linear}
u_{tt} - e^{-2Ht} \bigtriangleup u  + M^2 u=  0\,.
\end{equation}
The fundamental solutions and the Cauchy problem for the equation   with $M=0 $ in the backward direction of time are considered in \cite{Yag_Galst_Potsdam}.
\smallskip

The de~Sitter line element in the   higher dimensional analogue of de~Sitter space is
\[
ds^2 =-dt^2+e^{2Ht}\big( (dx^1)^2+ \ldots + (dx^n)^2\big)\,.
\]
It is a simplified version of the multidimensional cosmological models with the metric tensor given by 
\[
g =-e^{2\gamma (t)}dt^2+  e^{2\phi ^1(t)}g_1+ \ldots +  e^{2\phi ^n(t)}g_n \,,
\]
and can be chosen as a starting point for the study.  
The multidimensional cosmological models have attracted a lot of attention during recent years in constructing 
mathematical models of an anisotropic universe (see, e.g. \cite{Brozos-Vázquez,Heinzle-Rendall} and references therein). 
\smallskip

The  equation (\ref{K_G_linear}) is strictly hyperbolic. That implies the well-posedness of the Cauchy problem for (\ref{K_G_linear})  
in the different functional spaces. The coefficient of the equation is an analytic function and 
Holmgren's theorem  implies a local uniqueness in the space of distributions. 
Moreover, the speed of propagation is finite, namely, 
it is equal to $e^{-Ht} $ for every $ t \in {\mathbb R}$. 
The second-order strictly hyperbolic  equation (\ref{K_G_linear}) possesses two fundamental solutions 
resolving the Cauchy problem. They can be written microlocally in terms of the Fourier integral operators \cite{Horm}, which
give a complete description of the wave front sets of the solutions.
The distance between two characteristic roots $\lambda_1 (t,\xi ) $ and $\lambda_2 (t,\xi ) $ of the equation (\ref{K_G_linear}) is  
$|\lambda_1 (t,\xi ) - \lambda_2 (t,\xi )| = e^{-Ht}|\xi |$, $ t \in {\mathbb R}$, $\xi \in {\mathbb R}^n$. 
It tends to zero as $t $ approaches $\infty$. Thus, the operator is not uniformly (that is for all $t \in {\mathbb R} $) strictly hyperbolic. 
Moreover, the finite integrability
of the characteristic roots, $\int_0^{\infty}|\lambda_i (t,\xi )| dt <  \infty $,  leads to the 
existence of so-called ``horizon'' for that equation. More precisely,  any  signal emitted from the spatial point $x_0 \in {\mathbb R}^n$
at time $t_0 \in {\mathbb R} $ remains inside the ball $ |x -x_0 | <\frac{1}{H}e^{-Ht_0} $ 
for all time $t \in ( t_0,\infty) $. The  equation (\ref{K_G_linear})
is neither Lorentz  invariant nor  invariant with respect to usual scaling and that brings additional difficulties. 
In particular, it can cause a nonexistence of the $L^p-L^q$ decay 
for the solutions in the direction of time. In \cite{yagdjian_birk} it is mentioned 
the model equation with permanently bounded domain of influence,   power decay of characteristic roots,  and without
$L^p-L^q$ decay 
for the solutions that illustrates that phenomenon. The above mentioned $L^p-L^q$ decay 
estimates are one of the important tools for studying nonlinear equations
(see, e.g. \cite{Racke,Shatah}).
\smallskip

The time inversion 
 transformation $t \to -t$ reduces the  equation (\ref{K_G_linear}) to the mathematically equivalent  equation
\begin{equation}
\label{Gal_nl}
u_{tt} - e^{2Ht} \bigtriangleup u  + M^2 u=  0\,.
\end{equation}
The  wave equation,  that is equation  (\ref{Gal_nl}) with  $M=0$, 
 was investigated in \cite{Galstian} by the second author. More precisely,  in \cite{Galstian} the resolving operator for the Cauchy problem
\begin{equation}
\label{0.5new_intr}
u_{tt}   - e^{2t}\bigtriangleup u =0, \qquad u(x,0)= \varphi_0 (x), \quad u_t(x,0)= \varphi_1 (x)\,,
\end{equation}
is written as a sum of the Fourier integral operators with the amplitudes given in terms of the 
Bessel functions and in terms of confluent hypergeometric functions. In particular,  it is proved  in \cite{Galstian} that for $t>0$ 
the solution of the Cauchy problem  (\ref{0.5new_intr})   
is given by 
\begin{eqnarray*}
\hspace*{-0.5cm}  u(x,t) 
& = &
- i\frac{2}{(2\pi)^{n}} \int_{{\mathbb R}^n}
\Big\{ 
e^{i[x\cdot \xi +( e^t-1) |\xi|] } H_{+}\big(\frac{1}{2} ;1;2ie^t |\xi| \big) 
H_{-}\big(\frac{3}{2} ; 3;2i|\xi| \big) \nonumber \\
&  &
\hspace{1.8cm}-\,
e^{i[x\cdot \xi  -(e^t-1) |\xi| ]} H_{-}\big(\frac{1}{2} ;1;2ie^t |\xi| \big)
H_{+}\big(\frac{3}{2} ;3;2i|\xi|  \big)\Big\}
|\xi|^2   \hat  \varphi  _0 (\xi) d\xi  \nonumber \\
&   &
- i\frac{1}{(2\pi)^{n}} \int_{{\mathbb R}^n}
 \Big\{ 
e^{i[x\cdot \xi +( e^t-1) |\xi|] } H_{+}\big(\frac{1}{2} ;1;2ie^t |\xi| \big) 
H_{-}\big(\frac{1}{2} ; 1;2i|\xi| \big) \nonumber \\
&  &
\hspace{1.8cm}
-\,e^{i[x\cdot \xi  -(e^t-1) |\xi| ]} H_{-}(\frac{1}{2} ;1;2ie^t |\xi| )
H_{+}\big(\frac{1}{2} ;1;2i|\xi|  \big)\Big\}
\hat \varphi _1 (\xi) d\xi \,.
\end{eqnarray*}
In the notations of \cite{B-E} the last functions are   $H_{-}(\alpha ;\gamma ;z) = e^{i \alpha  \pi } \Psi (\alpha ;\gamma ;z) $
and $H_{+}(\alpha ;\gamma ;z) = e^{i \alpha  \pi } \Psi (\gamma -\alpha ;\gamma ;-z) $, where function $ \Psi ( a ;c ; z)  $ 
is defined in \cite[Sec.6.5]{B-E}. Here $\hat \varphi   (\xi) $ is a Fourier transform of $ \varphi   (x) $. 
\smallskip

The typical $L^p-L^q$ decay estimates 
obtained in \cite{Galstian} by dyadic decomposition of the phase space 
 contain some derivative loss.
More precisely, it is proved that for the solution $u = u(x,t)$ to the Cauchy problem (\ref{0.5new_intr}) with $n \ge 2$, 
$\varphi _0(x) \in C_0^\infty ({\mathbb R}^n)$ and $\varphi _1(x) =0$ 
for all large $t\ge T >0$,  the following estimate is satisfied
\begin{equation}
\label{1.16}
\|u(x,t)\|_{L^q({\mathbb R}^n)} 
  \le  
C (1+ e^t )^{ 
- \frac{1}{2} (n-1) (\frac{1}{p}-\frac{1}{q})}
\|\varphi _0\|_{W^N_p( {\mathbb R}^n)},
\end{equation}
 where $1<p\le 2$, $\frac{1}{p}+\frac{1}{q}=1$, and 
$\frac
{1}{2} (n+1) (\frac{1}{p}-\frac{1}{q})\le N < \frac
{1}{2} (n+1) (\frac{1}{p}-\frac{1}{q})+1$ and $W^N_p( {\mathbb R}^n) $ is the Sobolev space.
In particular, the derivative loss, $N$, is positive, unless $p=q=2$. 
This derivative loss phenomenon exists for the classical wave equation as well.
 Indeed, it is well-known  (see, e.g., \cite{Littman_1963,Littman,Peral}) that for the Cauchy problem
$ u_{tt}-\bigtriangleup u=0$, \,$u(x,0)=\varphi (x)$, \,$u_t(x,0)= 0 $,
the estimate
$ \| u(x,t)  \| _{ L^q  ({\mathbb R}_x^n)  }  \le C  \| \varphi  (x) \|_{L^q({\mathbb R}_x^n)}   $
fails to fulfill even for small positive  $t$ unless $q=2$. The obstacle  is created by the 
distinguishing feature of the (different from translation) Fourier integral operators of order zero, which   compose a resolving operator.
\smallskip

According to Theorem~1~\cite{Galstian}, for the solution $u = u(x,t)$ to the Cauchy problem (\ref{0.5new_intr}) with  $n \ge 2$, 
$\varphi _0(x) =0$ and $\varphi _1(x) \in C_0^\infty ({\mathbb R}^n)$ 
for all large $t\ge T >0$ and for any small $\varepsilon >0$, the following estimate is satisfied
\[
\|u(x,t)\|_{L^q({\mathbb R}^n)} 
 \le 
C _\varepsilon(1+t)(1+ e^t )^{ r_0- n(\frac{1}{p}-\frac{1}{q})}
\|\varphi _1\|_{W^N_p( {\mathbb R}^n)},\,
\]
where $1<p\le 2$, $\frac{1}{p}+\frac{1}{q}=1$, 
$r_0 = \max \{\varepsilon; \frac{(n+1)}{2} (\frac{1}{p}-\frac{1}{q})-\frac{1}{q}\}$,  
$\frac{n+1}{2}(\frac{1}{p}-\frac{1}{q})-\frac{1}{q}\le N < \frac{n+1}{2}(\frac{1}{p}-\frac{1}{q})+\frac{1}{p}$.
\medskip

The nonlinear equations (\ref{1.3}) and (\ref{K_G_semilinear})
are those we would like to solve, but the linear problem is a natural first step.
Exceptionally efficient tool for the studying nonlinear equations is a   fundamental 
solution of the associate linear operator.  
\smallskip

The fundamental solutions for the operator of the equation (\ref{0.5new_intr}) are constructed in  \cite{Yag_Galst_Potsdam}
and the representations of the solutions of the Cauchy problem 
\begin{equation}
\label{G_Y_CP}
u_{tt}   - e^{2t}\bigtriangleup u =f(x,t), \qquad u(x,0)= \varphi_0 (x), \quad u_t(x,0)= \varphi_1 (x)\,,
\end{equation}
are given in the terms of the solutions of wave equation in Minkowski spacetime. Then   in \cite{Yag_Galst_Potsdam}
for $n\geq 2$  the following decay estimate
\begin{eqnarray*} 
\hspace{-0.7cm}
\| (-\bigtriangleup )^{-s} u(x,t) \|_{ { L}^{  q} ({\mathbb R}^n)  } \!\! & \!\!\le &  \!\!\!\!
C   e^{t (  2s-n(\frac{1}{p}-\frac{1}{q})) } \int_{ 0}^{t} (1+  t -b  )\|  f(x, b)  \|_{ { L}^{p} ({\mathbb R}^n)  } \,db \nonumber \\
&  &
+ 
 C(e^ t-1)^{2s-n(\frac{1}{p}-\frac{1}{q}) } \left\{ \| \varphi_0  (x) \|_{ { L}^{p} ({\mathbb R}^n)  } 
+ \|\varphi_1  (x)  \|_{ { L}^{p} ({\mathbb R}^n)  }(1+ t ) (1-e^{-t}) \right\} 
\end{eqnarray*}
is proven, provided that $s \ge 0$, $1<p \leq 2$, $\frac{1}{p}+ \frac{1}{q}=1$, $\frac{1}{2} (n+1)\left( \frac{1}{p} - \frac{1}{q}\right) \leq 
2s \leq n \left( \frac{1}{p} - \frac{1}{q}\right) < 2s+1$. Moreover, this 
estimate  is fulfilled  for $n=1$ and $s=0$ as well as if $\varphi_0  (x)=0 $ and $\varphi_1  (x)=0 $. 
Case of $n=1$,  $f (x,t) =0 $, and non-vanishing $\varphi_1  (x)  $ and $\varphi_1  (x)  $ also is discussed in Section~8~\cite{Yag_Galst_Potsdam}.
\smallskip

In the construction of the fundamental solutions for the operator (\ref{K_G_linear}) 
 we follow the approach proposed in \cite{YagTricomi} that allows us to represent  
the fundamental solutions as some integral of the family of the fundamental solutions of the Cauchy problem for the wave equation without source term.
The kernel of that integral contains the Gauss's hypergeometric function. In that way, many properties of the wave equation can be extended to the
hyperbolic equations with the time dependent speed of propagation. That approach  was successfully applied in \cite{YagTricomi_GE,YagTricomi_JMAA} by the first author to
investigate the semilinear Tricomi-type equations.
\smallskip

Thus, in the present paper we consider  Klein-Gordon  operator in de~Sitter model of the universe, that is 
\[
{\mathcal S}:= \partial_t^2  - e^{-2t}\bigtriangleup  + M^2 \,,
\] 
where $M \geq 0$ is the reduced mass, $x \in {\mathbb R}^n$, $t \in {\mathbb R}$. We look for the  fundamental solution  $E=E(x,t;x_0,t_0)$,
\[
E_{tt} - e^{-2t}\Delta E  -M^2 E = \delta (x-x_0,t-t_0),
\]
with a support in the ``forward light cone'' $D_+ (x_0,t_0) $, $x_0 \in {\mathbb R}^n$, $t_0 \in {\mathbb R}$,
and for the  fundamental solution with a support in the ``backward light cone'' $D_- (x_0,t_0) $, $x_0 \in {\mathbb R}^n$, $t_0 \in {\mathbb R}$,
defined as follows
\begin{eqnarray}
\label{D+}
D_\pm (x_0,t_0) 
& := &
\Big\{ (x,t)  \in {\mathbb R}^{n+1}  \, ; \, 
|x -x_0 | \leq \pm( e^{-t_0} - e^{-t })
\,\Big\} \,.
\end{eqnarray}
In fact,  any intersection of  $ D_- (x_0,t_0) $ with the hyperplane $t=const <t_0$ determines the so-called dependence domain
for the point $(x_0,t_0) $, while the  intersection of  $ D_+ (x_0,t_0) $ 
with the hyperplane $t=const >t_0$ is the so-called  domain of influence of the point $(x_0,t_0) $. The equation (\ref{K_G_linear}) is
non-invariant    with respect to time inversion. Moreover, the dependence domain is wider than any given ball if time $const>t_0$ is sufficiently large, while 
the domain of influence is permanently, for all time  $const< t_0$, in the   ball of the radius $e^{t_0} $.
\smallskip

Define for $t_0 \in {\mathbb R}$ in the domain  $D_+ (x_0,t_0)\cup D_- (x_0,t_0) $  the function   
\begin{eqnarray} 
\label{E}
\hspace*{-0.5cm} &  &
E(x,t;x_0,t_0)\\
&  = &
(4e^{-t_0-t })^{iM} \Big((e^{-t }+e^{-t_0})^2 - (x - x_0)^2\Big)^{-\frac{1}{2}-iM    } 
F\Big(\frac{1}{2}+iM   ,\frac{1}{2}+iM  ;1; 
\frac{ ( e^{-t_0}-e^{-t })^2 -(x- x_0 )^2 }{( e^{-t_0}+e^{-t })^2 -(x- x_0 )^2 } \Big) , \nonumber 
\end{eqnarray} 
where $F\big(a, b;c; \zeta \big) $ is the hypergeometric function (See, e.g. \cite{B-E}.).
Let $E(x,t;0,b)$ be  function  (\ref{E}),
and set
\[ 
{\mathcal E}_+(x,t;0,t_0) 
 := 
\cases{ E(x,t;0,t_0) \quad \mbox{\rm in}  \,\, D_+ (0,t_0), \cr
0 \hspace*{2.0cm} \mbox{\rm elsewhere}} \,, \qquad
{\mathcal E}_{-}(x,t;0,t_0) 
 :=  
\cases{ E(x,t;0,t_0) \quad \mbox{\rm in}  \,\, D_- (0,t_0), \cr
0 \hspace*{2.0cm} \mbox{\rm elsewhere}} \,.
\]
Since  function $E=E(x,t;0,t_0)$ is smooth in  $D_{\pm} (0,t_0) $ and is locally integrable,  
it follows that ${\mathcal E}_+(x,t;0,$ $t_0) $ and ${\mathcal E}_-(x,t;0,$ $t_0) $ are  
distributions  whose supports are in   $D_{+}(0,t_0)$  and 
 $D_{-} (0,t_0) $, respectively.  The next theorem  gives our first result.

\begin{theorem}
\label{T1}
Suppose that $n=1$. The distributions ${\mathcal E}_+(x,t;0,t_0) $ and  ${\mathcal E}_-(x,t;0,t_0) $  are the fundamental solutions for the operator
${\mathcal S}= \partial_t^2  - e^{-2t}\bigtriangleup  + M^2$ relative to the point $(0,t_0)$, that is
\[
{\mathcal S} {\mathcal E}_{\pm}(x,t;0,t_0) = \delta (x,t-t_0), 
\]
or
\[
\frac{\partial^2}{\partial t^2}{\mathcal E}_{\pm}(x,t;0,t_0) - e^{-2t} \frac{\partial^2}{\partial x^2}{\mathcal E}_{\pm}(x,t;0,t_0)+ M^2{\mathcal E}_{\pm}(x,t;0,t_0)= \delta (x,t-t_0).
\]
\end{theorem}
\smallskip

To motivate one construction for the higher dimensional case $n\geq 2$ we follow the approach suggested in \cite{YagTricomi} and represent fundamental solution  
${\mathcal E}_+(x,t;0,t_0) $  as follows 
\begin{eqnarray*}
{\mathcal E}_+(x,t;0,t_0) 
& = & 
 \int_{ e^{-t} - e^{-t_0}}^  { e^{-t_0}- e^{-t}} {\mathcal E}^{string} (x,r )   (4e^{-t_0-t })^{iM} \Big((e^{-t_0 }+e^{-t})^2 - r^2\Big)^{-\frac{1}{2}-iM    }  \nonumber \\
&  &
\hspace{2.5cm} \times F\Big(\frac{1}{2}+iM   ,\frac{1}{2}+iM  ;1; 
\frac{ ( e^{-t_0}-e^{-t })^2 -r^2 }{( e^{-t_0}+e^{-t })^2 -r^2 } \Big) \, dr , \quad t>t_0,
\end{eqnarray*}
where the distribution ${\mathcal E}^{string} (x,t )  $ is the fundamental solution of the Cauchy problem for the string equation:
\[
\frac{\partial^2 }{\partial t^2}  {\mathcal E}^{string}   -   \frac{\partial^2 }{\partial x^2}  {\mathcal E}^{string} =0, \qquad  {\mathcal E}^{string}(x,0 )= \delta (x)
, \,\,\, {\mathcal E}^{string}_t  (x,0 )=0\,.
\] 
Hence, ${\mathcal E}^{string} (x,t ) = \frac{1}{2}\{ \delta (x+t )+ \delta (x-t )\}  $. The integral makes sense in the topology of the space of distributions.
The fundamental solution ${\mathcal E}_-(x,t;0,t_0) $ for $t<t_0$ admits a similar representation.
\smallskip

We appeal to
the wave equation in  Minkowski spacetime to obtain in the next theorem very similar representations of the fundamental solutions of the 
higher dimensional equation in  de~Sitter spacetime with $n\geq 2$.
\begin{theorem}
If $x \in {\mathbb R}^n$, $n\geq 2$, then for the operator $
{\mathcal S}=  \partial_t ^2  - e^{-2t}\bigtriangleup +M^2$ 
the  fundamental solution \, ${\mathcal E}_{+,n}(x,t;x_0,t_0) $ $(= {\mathcal E}_{+,n}(x-x_0,t;0,t_0))$ with a 
support in the  forward cone  $D_+ (x_0,t_0) $, $x_0 \in {\mathbb R}^n$, $t_0 \in {\mathbb R}$, \, 
supp$\,{\mathcal E}_{+,n} \subseteq D_+ (x_0,t_0)$, is given by the following integral $(t>t_0)$
\begin{eqnarray} 
\label{E+}
{\mathcal E}_{+,n}(x-x_0,t;0,t_0) 
& = &
2   
  \int_{ 0}^{ e^{-t_0}- e^{-t} }   {\mathcal E}^w (x-x_0,r )  (4e^{-t_0-t })^{iM} \Big((e^{-t_0 }+e^{-t})^2 - r^2\Big)^{-\frac{1}{2}-iM    }  \nonumber \\
&  &
\hspace{2.5cm} \times F\Big(\frac{1}{2}+iM   ,\frac{1}{2}+iM  ;1; 
\frac{ ( e^{-t_0}-e^{-t })^2 -r^2 }{( e^{-t_0}+e^{-t })^2 -r^2 } \Big) \, dr .
\end{eqnarray} 
Here the function 
${\mathcal E}^w(x,t;b)$   
is a fundamental solution to the Cauchy problem for the  wave equation
\[
{\mathcal E}^w_{ tt} -   \bigtriangleup {\mathcal E}^w  =  0 \,, \quad {\mathcal E}^w(x,0)=\delta (x)\,, \quad {\mathcal E}^w_{t}(x,0)= 0\,.
\]
The  fundamental solution ${\mathcal E}_{-,n}(x,t;x_0,t_0) $ $(= {\mathcal E}_{-,n}(x-x_0,t;0,t_0))$ with a 
support in the  backward cone  $D_- (x_0,t_0) $, $x_0 \in {\mathbb R}^n$, $t_0 \in {\mathbb R}$, \, 
supp$\,{\mathcal E}_{-,n} \subseteq  D_- (x_0,t_0)$, is given by the following integral $(t<t_0)$
\begin{eqnarray} 
\label{E-}
{\mathcal E}_{+,n}(x-x_0,t;0,t_0) 
& = &
-2   
  \int_{ e^{-t_0}- e^{-t} } ^{ 0}  {\mathcal E}^w (x-x_0,r )  (4e^{-t_0-t })^{iM} \Big((e^{-t_0 }+e^{-t})^2 - r^2\Big)^{-\frac{1}{2}-iM    }  \nonumber \\
&  &
\hspace{2.5cm} \times F\Big(\frac{1}{2}+iM   ,\frac{1}{2}+iM  ;1; 
\frac{ ( e^{-t_0}-e^{-t })^2 -r^2 }{( e^{-t_0}+e^{-t })^2 -r^2 } \Big) \, dr .
\end{eqnarray} 
\end{theorem}

In particular, the formula (\ref{E+}) shows that Huygens's Principle is not valid for the waves propagating in de~Sitter model of the universe. Fields satisfying a wave equation in de~Sitter model of the universe can be accompanied by {\it tails} propagating inside the light cone. This phenomenon
will be discussed in the spirit of \cite{Sonego-Faraoni}  in the forthcoming paper.
\smallskip

Next we use Theorem~\ref{T1}  to solve the Cauchy problem for the one-dimensional  equation 
\begin{equation}
\label{Int_3} 
u_{tt} - e^{-2t}u_{xx} +M^2 u = f(x,t)\,, \qquad t > 0\,, \quad x \in {\mathbb R} \,, 
\end{equation}
with vanishing initial data,
\begin{equation}
\label{Int_4} 
u(x,0)= u_t(x,0)= 0\,.
\end{equation}
\begin{theorem}
\label{T1.1}
Assume that the function $f$ is continuous along with its all second order derivatives, 
and that for every fixed $t$ it has a compact support, supp$f(\cdot ,t) \subset {\mathbb R}$.
Then the function $u=u(x,t)$ defined by  
\begin{eqnarray*}
u(x,t)  
& = &
 \int_{ 0}^{t} db \int_{ x - (e^{-b}- e^{-t})}^{x+ e^{-b}- e^{-t}}  dy\, f(y,b)  
(4e^{-b-t })^{iM} \Big((e^{-t }+e^{-b})^2 - (x - y)^2\Big)^{-\frac{1}{2}-iM    } \\
&  &
\hspace{2.5cm}\times F\Big(\frac{1}{2}+iM   ,\frac{1}{2}+iM  ;1; 
\frac{ ( e^{-b}-e^{-t })^2 -(x- y )^2 }{( e^{-b}+e^{-t })^2 -(x- y )^2 } \Big)   \nonumber 
\end{eqnarray*}
is a $C^2$-solution to the Cauchy problem for the equation (\ref{Int_3}) with  vanishing initial data,
(\ref{Int_4}).
\end{theorem}
\smallskip

The representation of the solution of the Cauchy problem for the one-dimensional case ($n=1$) 
of the equation (\ref{K_G_linear}) without source term is given by the next theorem. 
\begin{theorem}
\label{T1.3}
The solution $u=u (x,t)$ of the Cauchy problem 
\begin{equation}
\label{oneDphy01}
u_{tt} - e^{-2t}u_{xx} +M^2 u =0\, ,\qquad u(x,0) = \varphi_0  (x) \,, \qquad u_t(x,0) =\varphi_1  (x)\,  ,
\end{equation}
with \,$\varphi_0    ,  \varphi_1  \in C_0^\infty ({\mathbb R})$ can be represented as follows
\begin{eqnarray*}
u(x,t)   
&  =   &
\frac{1}{2} e ^{\frac{t}{2}}  \Big[ 
\varphi_0   (x+ 1-e^{-t})  
+     \varphi_0   (x -  1 +e^{-t})  \Big]  
+ \int_{ 0}^{1-  e^{-t}} \big[ 
\varphi_0   (x - z)  
+     \varphi_0   (x  + z)  \big] K_0(z,t)\,  dz \\
&  &    
+ \,\,\int_{0}^{1-  e^{-t}} \,\Big[      \varphi_1    (x- z)  +   \varphi _1   (x + z)    \Big] K_1(z,t) dz \,,
\end{eqnarray*}
where the kernels  $K_0(z,t)   $    and $K_1(z,t)   $ are defined by
\begin{eqnarray*}
K_0(z,t)
& := & 
- \left[  \frac{\partial }{\partial b}   E(z,t;0,b) \right]_{b=0} \\
&  = &   
 (4e^{-t })^{iM} \big((e^{-t }+1)^2 - z^2\big)^{ -iM    } \frac{1}{ [(1-e^{ -t} )^2 -  z^2]\sqrt{(e^{-t }+1)^2 - z^2} }\\
&   &
\times  \Bigg[  \big(  e^{-t} -1 - iM(e^{ -2t} -      1 -  z^2) \big) 
F \Big(\frac{1}{2}+iM   ,\frac{1}{2}+iM  ;1; \frac{ ( 1-e^{-t })^2 -z^2 }{( 1+e^{-t })^2 -z^2 }\Big) \\
&  &
\hspace{1cm}  +   \big( 1-e^{-2 t}+  z^2 \big)\Big( \frac{1}{2}-iM\Big)
F \Big(-\frac{1}{2}+iM   ,\frac{1}{2}+iM  ;1; \frac{ ( 1-e^{-t })^2 -z^2 }{( 1+e^{-t })^2 -z^2 }\Big) \Bigg]
\end{eqnarray*} 
and
\begin{eqnarray*} 
K_1(z,t) 
& := &
  E(z ,t;0,0) \\
& = &
 (4e^{ -t })^{iM} \big((1+e^{-t })^2 -   z  ^2\big)^{-\frac{1}{2}-iM    } 
F\left(\frac{1}{2}+iM   ,\frac{1}{2}+iM  ;1; 
\frac{ ( 1-e^{-t })^2 -z^2 }{( 1+e^{-t })^2 -z^2 } \right), \, 0\leq z\leq  1-e^{-t}, 
 \end{eqnarray*}  
respectively.
\end{theorem}

The kernels $K_0(z,t)  $ and $K_1(z,t)  $ play leading roles in the derivation of
$L^p-L^q$ estimates. Their main properties are listed and proved in Section~\ref{S11}.

\smallskip

Next we turn to the higher-dimensional equation with $n\geq 2$.

\begin{theorem}
\label{T1.5}
If $n$ is odd, $n=2m+1$,  $m  \in {\mathbb N}$,  then the solution $u= u(x,t)$ to the Cauchy problem   
\begin{equation}
\label{1.26}
u_{tt} - e^{-2t}\Delta u +M^2 u= f ,\quad u(x,0)= 0  , \quad u_t(x,0)=0,
\end{equation}
with \, $ f \in C^\infty ({\mathbb R}^{n+1})$\, and with the vanishing
initial data is given by the next expression
\begin{eqnarray} 
\label{1.27}
\hspace{-0.3cm} u(x,t) 
& = &
2\int_{ 0}^{t} db
  \int_{ 0}^{ e^{-b}- e^{-t}} dr_1 \,  \left(  \frac{\partial }{\partial r} 
\Big(  \frac{1}{r} \frac{\partial }{\partial r}\Big)^{\frac{n-3}{2} } 
\frac{r^{n-2}}{\omega_{n-1} c_0^{(n)} }  \!\!\int_{S^{n-1} } f(x+ry,b)\, dS_y  
\right)_{r=r_1}   \nonumber \\
  &  & 
\times  (4e^{-b-t})^{iM}
\left( (e^{-t}  + e^{-b} )^2 - r_1^2   \right)^{-\frac{1}{2}-iM}
F\left(\frac{1}{2}+iM,\frac{1}{2}+iM;1; 
\frac{ (e^{-b}- e^{-t})^2-r_1^2}
{  (e^{-b}+ e^{-t})^2-r_1^2} \right)\!\! ,
\end{eqnarray} 
where $c_0^{(n)} =1\cdot 3\cdot \ldots \cdot (n-2 )$.
Constant $\omega_{n-1} $ is the area of the unit sphere $S^{n-1} \subset {\mathbb R}^n$. 

If $n$ is even, $n=2m$,  $m  \in {\mathbb N}$,  then the solution $u= u(x,t)$  is given by the next expression 
\begin{eqnarray}
\label{1.28} 
\hspace{-0.3cm} u(x,t) 
& = &
2\int_{ 0}^{t} db
  \int_{ 0}^{ e^{-b}- e^{-t}} dr_1 \,  \left( \frac{\partial }{\partial r} 
\Big( \frac{1}{r} \frac{\partial }{\partial r}\Big)^{\frac{n-2}{2} } 
\frac{2r^{n-1}}{\omega_{n-1} c_0^{(n)} }  \!\!\int_{B_1^{n}(0)} \frac{f(x+ry,b) }{\sqrt{1-|y|^2}} \, dV_y 
\right)_{r=r_1}   \nonumber \\
  &  & 
\times  (4e^{-b-t})^{iM}
\left( (e^{-t}  + e^{-b} )^2 - r_1^2   \right)^{-\frac{1}{2}-iM}
F\left(\frac{1}{2}+iM,\frac{1}{2}+iM;1; 
\frac{ (e^{-b}- e^{-t})^2-r_1^2}
{  (e^{-b}+ e^{-t})^2-r_1^2} \right)\!\! .
\end{eqnarray} 
Here $B_1^{n}(0) :=\{|y|\leq 1\} $ is the unit ball in ${\mathbb R}^n$, while $c_0^{(n)} =1\cdot 3\cdot \ldots \cdot (n-1)$.
\end{theorem}
Thus, in both cases, of even and odd $n$, one can write 
\begin{eqnarray}
\label{1.29} 
u(x,t) 
&  =  &
2   \int_{ 0}^{t} db
  \int_{ 0}^{ e^{-b}- e^{-t}} dr  \,  v(x,r ;b)  (4e^{-b-t})^{iM}
\left( (e^{-t}  + e^{-b} )^2 - r^2   \right)^{-\frac{1}{2}-iM}  \nonumber \\
&  &
\hspace{2.5cm} \times  
F\left(\frac{1}{2}+iM,\frac{1}{2}+iM;1; 
\frac{ (e^{-b}- e^{-t})^2-r^2}
{  (e^{-b}+ e^{-t})^2-r^2} \right) ,
\end{eqnarray}
where the function 
$v(x,t;b)$   
is a solution to the Cauchy problem for the  wave equation
\[
v_{tt} -   \bigtriangleup v  =  0 \,, \quad v(x,0;b)=f(x,b)\,, \quad v_t(x,0;b)= 0\,.
\]

The next theorem gives representation of the solutions of  equation (\ref{K_G_linear}) with the initial data prescribed at $t=0$.  
\begin{theorem}
\label{T1.6}
The solution $u=u (x,t)$ to the Cauchy problem  
\begin{equation}
\label{1.30CP}
u_{tt}-  e^{-2t} \bigtriangleup u +M^2 u =0\,, \quad u(x,0)= \varphi_0 (x)\, , \quad u_t(x,0)=\varphi_1 (x)\,, 
\end{equation}
with \, $ \varphi_0 $,  $ \varphi_1 \in C_0^\infty ({\mathbb R}^n) $, $n\geq 2$, can be represented as follows:
\begin{eqnarray}
\label{1.30}
u(x,t) 
& = &
 e ^{\frac{t}{2}} v_{\varphi_0}  (x, \phi (t))
+ \, 2\int_{ 0}^{1} v_{\varphi_0}  (x, \phi (t)s) K_0(\phi (t)s,t)\phi (t)\,  ds  \nonumber \\
& &
+\, 2\int_{0}^1   v_{\varphi _1 } (x, \phi (t) s) 
  K_1(\phi (t)s,t) \phi (t)\, ds 
, \quad x \in {\mathbb R}^n, \,\, t>0\,,  
\end{eqnarray}
$\phi (t):= 1-e^{-t} $, by means of the kernels  $K_0$ and $K_1$ have been defined in Theorem~\ref{T1.3}. 
 Here for $\varphi \in C_0^\infty ({\mathbb R}^n)$ and for $x \in {\mathbb R}^n$, $n=2m+1$, $m \in {\mathbb N}$,  
\begin{eqnarray*}
 v_\varphi  (x, \phi (t) s) :=  
\Bigg( \frac{\partial}{\partial r} \Big( \frac{1}{r} \frac{\partial }{\partial r}\Big)^{\frac{n-3}{2} } 
\frac{r^{n-2}}{\omega_{n-1} c_0^{(n)} } \int_{S^{n-1}  } 
\varphi (x+ry)\, dS_y   \Bigg)_{r=\phi (t) s}   
\end{eqnarray*}
while for $x \in {\mathbb R}^n$, $n=2m$,  $m \in {\mathbb N}$,   
\[
v_\varphi  (x, \phi (t) s) :=  \Bigg( \frac{\partial }{\partial r}  
\Big( \frac{1}{r} \frac{\partial }{\partial r}\Big)^{\frac{n-2}{2} } 
\frac{2r^{n-1}}{\omega_{n-1} c_0^{(n)}} \int_{B_1^{n}(0)}  \frac{1}{\sqrt{1-|y|^2}}\varphi (x+ry)\, dV_y 
 \Bigg)_{r=s \phi (t)}  \,.
\]
The function $v_\varphi  (x, \phi (t) s)$  coincides with the value $v(x, \phi (t) s) $ 
of the solution $v(x,t)$ of the Cauchy problem
\[
v_{tt}-  \bigtriangleup v =0, \quad v(x,0)= \varphi (x), \quad v_t(x,0)=0\,.
\]
\end{theorem}

\smallskip

As a consequence  of the above theorems  we obtain in Sections~\ref{S12}-\ref{S13}
for $n\geq 2$  and for the particles with ``large'' mass $m$, $m \geq n/2$, 
that is, with nonnegative  curved mass $M \geq 0$, 
the following $L^p-L^q $ estimate
\begin{eqnarray}
\label{1.32} 
\hspace{-0.7cm}  \| (-\bigtriangleup )^{-s} u(x,t) \|_{ { L}^{  q} ({\mathbb R}^n)  }  
\!\! & \!\!  \le \!\! &   \!\! 
C\int_{ 0}^{t} \|  f(x, b)  \|_{ { L}^{p} ({\mathbb R}^n)  }   
 e^{-b}\left( e^{-b}- e^{-t}  \right)^{ 1+  2s-n(\frac{1}{p}-\frac{1}{q}) } \left(1+  t-  b    \right)  db  \nonumber \\
&  &
+ C  (1+ t  )(1- e^{-t}) ^{ 2s-n(\frac{1}{p}-\frac{1}{q}) }\Big\{ e ^{\frac{t}{2}}  \| \varphi_0  (x) \|_{ { L}^{p} ({\mathbb R}^n)  }
+ (1- e^{-t}) \|\varphi_1  \|_{ { L}^{p} ({\mathbb R}^n)  } 
 \Big\} 
\end{eqnarray}
provided that $s \ge 0$, $1<p \leq 2$, $\frac{1}{p}+ \frac{1}{q}=1$, $\frac{1}{2} (n+1)\left( \frac{1}{p} - \frac{1}{q}\right) \leq 
2s \leq n \left( \frac{1}{p} - \frac{1}{q}\right) < 2s+ 1 $. Moreover, according to Theorem~\ref{T10.1}
the estimate (\ref{1.32}) is valid   for $n=1$ and $s=0$ as well as if $\varphi_0  (x)=0 $ and $\varphi_1  (x)=0 $. Case of $n=1$,  $f (x,t) =0 $, and non-vanishing $\varphi_1  (x)  $ and $\varphi_1  (x)  $ is discussed in Section~\ref{S11}.
The case of particles with small mass $m<n/2$ will be discussed in the forthcoming paper.
\smallskip

The paper is organized as follows. In Section~\ref{S2} we construct the Riemann function 
of the operator of (\ref{K_G_linear}) in the characteristic coordinates for the case of $n=1$.  That Riemann function used in Section~\ref{S2a}
to prove Theorem~\ref{T1}. 
Then in Section~\ref{S3} we apply the fundamental solutions to solve the 
Cauchy problem with the source term and with the vanishing initial data  given at $t=0$. 
More precisely, we give a representation formula for the solutions. In that section we also prove several basic properties
of the function $E(x,t;y,b)$.  In Sections~\ref{S5}-\ref{S6} we use formulas of  
Section~\ref{S3} to derive and to complete the list of   representation formulas 
for the solutions of the Cauchy problem for the case of one-dimensional spatial variable. The higher-dimensional
equation with the source term is considered in  Section~\ref{S7}, where we derive a representation formula for the solutions of 
the Cauchy problem with the source term and with the vanishing initial data  given at $t=0$. In the same section this formula is used to 
derive the fundamental solutions of the operator and to complete the proof of Theorem~\ref{T1.6}.
Then in Sections~\ref{S10}-\ref{S13} we establish the $L^p-L^q$ decay estimates. 
Applications of all these results to the nonlinear equations will be done in the forthcoming paper.

\section{The Riemann Function}
\label{S2}

\setcounter{equation}{0}

In the characteristic coordinates $l$ and $m$, 
\begin{equation}
\label{Gel1.8}
l =  x+ e^{-t}, \qquad
m =  x- e^{-t}, 
\end{equation}
one has
\begin{eqnarray*}
\frac{\partial^2  }{\partial t ^2}
& = &
\frac{1}{2}(l-m)\left( \frac{\partial }{\partial l  } -  \frac{\partial }{\partial m  } \right)
+ \frac{1}{4}(l-m)^2 \left( \frac{\partial ^2}{\partial l  ^2} 
   -2 \frac{\partial^2 }{\partial l  \, \partial m } +  \frac{\partial ^2}{\partial m  ^2}  \right)  \\
\end{eqnarray*}
and
\begin{eqnarray*}
&  &
 e^{-2t}  \frac{\partial  ^2}{\partial x  ^2} 
= \frac{1}{4}(l-m)^2 \left( \frac{\partial ^2}{\partial l  ^2}+  2\frac{\partial^2 }{\partial l  \, \partial m }+  \frac{\partial^2 }{\partial m  ^2} \right) .
\end{eqnarray*}
Then the operator\, ${\mathcal S}  $ \, of the equation (\ref{Int_3}) 
 reads
\begin{equation}
\label{Slm}
{\mathcal S} := \frac{\partial^2  }{\partial t ^2} -  e^{-2t}  \frac{\partial^2   }{\partial x^2 } +  M^2
 =  
- (l-m)^{2 } \Bigg\{ \frac{\partial^2  }{\partial l \, \partial m} - \frac{1}{2(l-m)} \Big( \frac{\partial }{\partial l} -
\frac{\partial }{\partial m} \Big)  - \frac{1}{(l-m)^2}M^2 \Bigg\}\,.
\end{equation}
In particular, in the new variables the equation
\[
\left( \frac{\partial^2  }{\partial t ^2} -  e^{-2t}  \frac{\partial^2   }{\partial x^2 } + M^2\right) u=0 
 \quad \mbox{\rm implies}   \quad
 \Bigg\{ \frac{\partial^2  }{\partial l \, \partial m} - \frac{1}{2(l-m)} \Big( \frac{\partial }{\partial l} -
\frac{\partial }{\partial m} \Big) \Bigg\}u -\frac{1}{(l-m)^{2 }} M^2u=0\,.
\]
We need the following lemma  with $\displaystyle \gamma = \frac{1}{2}+  iM$. 
\begin{lemma}
The function
\[
V(l,m;a,b)  = 
(l-b)^{-\gamma } (a-m)^{-\gamma } 
F\Big(\gamma ,\gamma ;1; \frac{(l-a)(m-b)}{(l-b)(m-a)} \Big)
\]
solves the equation
\begin{equation}
\label{eqV}
\Bigg\{ \frac{\partial^2  }{\partial l \, \partial m} - \frac{1}{ (l-m)} \gamma  \Big( \frac{\partial }{\partial l} -
\frac{\partial }{\partial m} \Big) \Bigg\}  V (l,m;a,b) = 0\,.
\end{equation}
\end{lemma}

\noindent
{\bf Proof.} We denote the argument of the hypergeometric function by $z$, and evaluate its derivatives,
\[
z:= \frac{(l-a)(m-b)}{(l-b)(m-a)}\,, \qquad \frac{\partial }{\partial l}z = \frac{(a-b)(b-m)}{(l-b)^2(a-m)}\,, \quad 
\frac{\partial }{\partial m}z= -\frac{(a-b)(a-l)}{(b-l)(a-m)^2} \,.
\]
Further, we obtain
\[ 
 \frac{\partial }{\partial l} V(l,m;a,b)  
 = 
(a-m)^{-\gamma }   (l-b)^{-\gamma -1} \Bigg\{ -\gamma   
F\Big(\gamma ,\gamma ;1; z \Big) 
+ 
\frac{(a-b)(b-m)}{(l-b)(a-m)}  F'_z\Big(\gamma ,\gamma ;1; z \Big)  \Bigg\} 
\,.
\] 
Next
\[
 \frac{\partial }{\partial m} V(l,m;a,b)  
  =   
(a-m)^{-\gamma }   (l-b)^{-\gamma } \Bigg\{ \gamma  (a-m)^{-1} 
F\Big(\gamma ,\gamma ;1; z \Big)  
 -\frac{(a-b)(a-l)}{(b-l)(a-m)^2}F_z'\Big(\gamma ,\gamma ;1; z \Big) \Bigg\} \,.
\] 
Then 
\begin{eqnarray*}
 \left( \frac{\partial }{\partial l} -  \frac{\partial }{\partial m} \right) V(l,m;a,b) 
& = &
(a-m)^{-\gamma }   (l-b)^{-\gamma  } \left(-\gamma \right) 
F\big(\gamma ,\gamma ;1; z \big)\Big\{  \frac{1}{(l-b)} + \frac{1}{(a-m)} \Big\}\\
&  &
+ (a-m)^{-\gamma -1}   (l-b)^{-\gamma  -1} F_z'\big(\gamma ,\gamma ;1; z \big)
(a-b)  \left\{ \frac{(b-m)}{(l-b)} - \frac{(a-l)}{(a-m)}\right\} \,.
\end{eqnarray*} 
Furthermore,
\begin{eqnarray*}
  \frac{\partial^2  }{\partial l\,\partial m} V(l,m;a,b) 
& = &
(a-m)^{-\gamma -1}  \Bigg[    -\gamma    (l-b)^{-\gamma -1} \Bigg\{  \gamma    
F\Big(\gamma ,\gamma ;1; z \Big)  
-\frac{(a-b)(a-l)}{(b-l)(a-m)}F_z'\Big(\gamma ,\gamma ;1; z \Big) \Bigg\} \Bigg] \\
&  &
+ \,(a-m)^{-\gamma -1} \Bigg[    (l-b)^{-\gamma }   \frac{\partial }{\partial l}\Bigg\{ \gamma  
F\Big(\gamma ,\gamma ;1; z \Big)  
-\frac{(a-b)(a-l)}{(b-l)(a-m)}F_z'\Big(\gamma ,\gamma ;1; z \Big) \Bigg\} \Bigg] .
\end{eqnarray*}
We calculate
\begin{eqnarray*}
&  &
\frac{\partial }{\partial l}\Bigg\{\gamma   
F\Big(\gamma ,\gamma ;1;z \Big) 
-\frac{(a-b)(a-l)}{(b-l)(a-m)}F_z'\Big(\gamma ,\gamma ;1; z\Big) \Bigg\} \\
& =&
 \gamma   
F'_z\Big(\gamma ,\gamma ;1; z \Big) \frac{(a-b)(b-m)}{(l-b)^2(a-m)}
- \frac{(a-b)^2 }{(b-l)^2(a-m)}  
F_z'\Big(\gamma ,\gamma ;1; z \Big)  \\
&  &
-
\frac{(a-b)(a-l)}{(b-l)(a-m)} 
F_{zz}{}'{}'\Big(\gamma ,\gamma ;1; z \Big) \frac{(a-b)(b-m)}{(l-b)^2(a-m)}  \,.
\end{eqnarray*}
Here $ \frac{\partial }{\partial l}\frac{(a-b)(a-l)}{(b-l)(a-m)} =  \frac{(a-b)^2 }{(b-l)^2(a-m)} $. 
Finally,
\begin{eqnarray*}
\frac{\partial^2  }{\partial l \, \partial m}  V(l,m;a,b) 
& = &
(a-m)^{-\gamma -1}  \Bigg[   -\gamma   (l-b)^{-\gamma -1} \Bigg\{  \gamma   
F\Big(\gamma ,\gamma ;1; z \Big)  
-\frac{(a-b)(a-l)}{(b-l)(a-m)}F_z'\Big(\gamma ,\gamma ;1; z \Big) \Bigg\} \Bigg] \\
&  &
+ \,(a-m)^{-\gamma -1}      (l-b)^{-\gamma }  \Bigg\{   \gamma   
F'_z\Big(\gamma ,\gamma ;1; z  \Big) \frac{(a-b)(b-m)}{(l-b)^2(a-m)}
 \\
&  &
 - \frac{(a-b)^2 }{(b-l)^2(a-m)}  
F_z'\Big(\gamma ,\gamma ;1;  z \Big) -
\frac{(a-b)(a-l)}{(b-l)(a-m)} 
F_{zz}{}'{}'\Big(\gamma ,\gamma ;1; z  \Big) \frac{(a-b)(b-m)}{(l-b)^2(a-m)}   \Bigg\}   .
\end{eqnarray*}
The coefficients of the   derivatives of the hypergeometric function, $ F_{zz}{}'{}'$, $ F_{z }{}' $, and of $ F $ in the expression for  
$ \displaystyle \frac{\partial^2  }{\partial l \, \partial m}  V(l,m;a,b) $ are
\begin{eqnarray*}
&  &
(a-m)^{-\gamma -3}      (l-b)^{-\gamma -3}    (a-b)^2(a-l) (b-m) \,,\\
&  &
(a-m)^{-\gamma -2}      (l-b)^{-\gamma -2} (a-b)\Big\{    \gamma   
   (b-m-a+l)  -  (a-b)     \Big\} \,,\\
&  &
(a-m)^{-\gamma -1}    \left( -\gamma \right)  (l-b)^{-\gamma -1}  \{ - \left(-\gamma \right)   \}
\, = \,
- (a-m)^{-\gamma -1}   (l-b)^{-\gamma -1}   \gamma ^2 \,, 
\end{eqnarray*}
respectively.
The coefficients of $ F_{z }{}' $  and   $ F  $  in the expression for 
$ \displaystyle  \frac{1}{ (l-m)} \gamma \left( \frac{\partial }{\partial l} -  \frac{\partial }{\partial m} \right) V(l,m;a,b)$ are 
\[
\frac{1}{ (l-m)} \gamma (a-m)^{-\gamma -1}   (l-b)^{-\gamma  -1} 
(a-b)  \left\{ \frac{(b-m)}{(l-b)} - \frac{(a-l)}{(a-m)}\right\}   
\]
and
\[
- \,\frac{1}{ (l-m)}  \gamma  (a-m)^{ -\gamma }   (l-b)^{-\gamma  } \gamma  \Big\{  \frac{1}{(l-b)} + \frac{1}{(a-m)} \Big\} \,.
\]
Now we turn to the equation (\ref{eqV}). The coefficients of  $ F  $  and  $F_{z }{}' $  in that equation are 
\[
 \frac{1}{ (m-l)}(a-b)  (a-m)^{-\gamma -1}  (l-b)^{-\gamma -1} ( - \gamma ^2)   , 
\]
and
\begin{eqnarray*}  
&  &
\frac{1}{ (m-l)}(a-m)^{-\gamma -1}      (l-b)^{-\gamma -1} (a-b) \\
&  &
\times \Bigg[   \gamma    \frac{(m-l)}{(a-m)(l-b)} 
   (b-m-a+l)     + \gamma  \left\{ \frac{(b-m)}{(l-b)} - \frac{(a-l)}{(a-m)}\right\}  
 -   \frac{(m-l)}{(a-m)(l-b)} (a-b) \Bigg]  .
\end{eqnarray*}
The first two terms in the brackets can be written as follows
\begin{eqnarray*} 
\gamma \frac{(m-l)}{(a-m)(l-b)} 
   (b-m-a+l)     +  \gamma  \left\{  \frac{(b-m)}{(l-b)} - \frac{(a-l)}{(a-m)}\right\}
  =  -2\gamma z \,,
\end{eqnarray*}
while the last term can be transformed to 
\begin{eqnarray*} 
-   \frac{(m-l)}{(a-m)(l-b)} (a-b) 
  =  1-z \,.
\end{eqnarray*}
Thus, the coefficient of $ F_{z }{}' $  in the equation  (\ref{eqV}) is 
\begin{eqnarray*}
&  &
\frac{1}{ (m-l)}(a-m)^{-\gamma -1}      (l-b)^{-\gamma -1} (a-b)\left[1   - (1+2  \gamma)  z    \right]   \,.
\end{eqnarray*}
Finally, the coefficient of $ F_{z z}{}' {}'$  in the equation  (\ref{eqV}) is 
\begin{eqnarray*}
 &  &
\frac{1}{ (m-l)}(a-m)^{-\gamma -1}      (l-b)^{-\gamma -1}    (a-b)\Big[\frac{ (a-b)(a-l) (b-m)(m-l)  }{(a-m)^{2} (l-b)^{2} }\Big]\,,
\end{eqnarray*}
where
\begin{eqnarray*}
 &  &
\frac{ (a-b)(a-l) (b-m)(m-l)  }{(a-m)^{2} (l-b)^{2} }= z(1-z)\,.
\end{eqnarray*}
Hence, the left-hand side of (\ref{eqV}) reads
\begin{eqnarray*}
 &  &
\Bigg\{ \frac{\partial^2  }{\partial l \, \partial m} - \frac{1}{ (l-m)} \gamma  \Big( \frac{\partial }{\partial l} -
\frac{\partial }{\partial m} \Big) \Bigg\}  V (l,m;a,b) \\
& = &
\frac{1}{ (m-l)}(a-m)^{-\gamma -1}      (l-b)^{-\gamma -1}    (a-b)\Big[ z(1-z)F_{z z} + \big( 1   - (1+2  \gamma)  z \big) F_{z  } -
\gamma ^2 F \Big] = 0\,,
\end{eqnarray*}
and vanishes, since $F$ solves  the Gauss hypergeometic equation with $c=1$, $a= \gamma  $, and $b= \gamma $. Lemma is proven. \hfill $\square$

\begin{lemma}
\label{L1.2}
For $\gamma \in {\mathbb C}$ such that $ F\big(\gamma ,\gamma ;1; z \big)$ is well defined, the function 
\begin{eqnarray*}
E(l,m;a,b) 
& := &
(a-b)^{\gamma -\frac{1}{2}}(l-m)^{\gamma -\frac{1}{2}}V(l,m;a,b) \\
& = &
(a-b)^{\gamma -\frac{1}{2}}(l-m)^{\gamma -\frac{1}{2}}(l-b)^{-\gamma } (a-m)^{-\gamma } 
F\Big(\gamma ,\gamma ;1; \frac{(l-a)(m-b)}{(l-b)(m-a)} \Big)
\end{eqnarray*}
solves the equation
\begin{equation}
\Bigg\{ \frac{\partial^2  }{\partial l \, \partial m} - \frac{1}{  2(l-m)}   \Big( \frac{\partial }{\partial l} -
\frac{\partial }{\partial m} \Big) \Bigg\}  E (l,m;a,b) + \frac{1}{(l-m)^2}\left( \frac{1}{2}- \gamma \right)^2 E (l,m;a,b)= 0\,.
\end{equation}
\end{lemma}
\medskip

\noindent{\bf Proof.} Indeed, straightforward calculations show 
 \begin{eqnarray*}
&  &
(a-b)^{-\gamma +\frac{1}{2}}\Bigg\{ \frac{\partial^2  }{\partial l \, \partial m} - \frac{1}{  2(l-m)}   \Big( \frac{\partial }{\partial l} -
\frac{\partial }{\partial m} \Big)  + \frac{1}{(l-m)^2}\left( \frac{1}{2}- \gamma \right)^2 \Bigg\}  E \\
& = &
(l- m)^{\gamma -\frac{1}{2}} \left[ V_{lm}  -  \frac{1}{   (l-m)} \gamma   \left( V_{l} - V_{m} \right)  \right] =0 \,.
\end{eqnarray*}
Lemma is proven. \hfill $\square$
\medskip

Consider now the operator 
\[
{\mathcal S}_{ch}^*  := \frac{\partial^2  \,}{\partial l \, \partial m} + \frac{1}{2(l-m)} \Big( \frac{\partial \, }{\partial l} -
\frac{\partial \,}{\partial m} \Big)    - \frac{1}{(l-m)^{2 }} (M^2 +1)\,,
\]
which is a formally adjoint to the operator
\[
{\mathcal S}_{ch}  := \frac{\partial^2  \,}{\partial l \, \partial m} - \frac{1}{2(l-m)} \Big( \frac{\partial \, }{\partial l} -
\frac{\partial \,}{\partial m} \Big)  - \frac{1}{(l-m)^{2 }} M^2\,.
\]
In the next lemma the Riemann function is presented.

\begin{proposition}
\label{PR}
The function 
\begin{eqnarray*}
R(l,m;a,b) 
& = &
(l-m)E(l,m;a,b) \\
& = &
(a-b)^{iM}(l-m)^{1+iM}(l-b)^{-\frac{1}{2}-iM   } (a-m)^{-\frac{1}{2}-iM    } 
F\Big(\frac{1}{2}+iM   ,\frac{1}{2}+iM  ;1; \frac{(l-a)(m-b)}{(l-b)(m-a)} \Big)
\end{eqnarray*}
defined for all $l$, $m$, $a$, $b \in {\mathbb R}$, such that    $ l>m $,   
is a unique solution of the equation ${\mathcal S}_{ch} ^*R=0$, which satisfies the following conditions: \\
$(\mbox{\rm i})$ $\displaystyle{ R_l =  \frac{1}{2 (l-m)}R \quad }$ along the line $m=b$;\\
$(\mbox{\rm ii})$ $\displaystyle{ R_m =  - \frac{1}{2 (l-m)}R \quad }$ along the line $l=a$; \\
$(\mbox{\rm iii})$ $R(a,b;a,b)=  1$.
\end{proposition}
\medskip

\noindent
{\bf Proof.} 
Indeed, if we denote $\gamma  = \frac{1}{2}+iM $, then for the Riemann function we have
\[
R(l,m;a,b) 
  =  
(a-b)^{\gamma -\frac{1}{2}}(l-m)^{\gamma +\frac{1}{2}}V(l,m;a,b) 
  =  
(l-m)E(l,m;a,b)\,.
\]
The operators ${\mathcal S}_{ch} $ and ${\mathcal S}_{ch} ^*$ can be written as follows: 
\begin{eqnarray*}
{\mathcal S}_{ch} 
&  = &
\frac{\partial^2  \,}{\partial l \, \partial m} - \frac{1}{2(l-m)} \Big( \frac{\partial \, }{\partial l} -
\frac{\partial \,}{\partial m} \Big)  + \frac{1}{(l-m)^{2 }} \left( \gamma -\frac{1}{2}\right)^2\,,\\
{\mathcal S}_{ch}^* 
&  = &
\frac{\partial^2  \,}{\partial l \, \partial m} + \frac{1}{2(l-m)} \Big( \frac{\partial \, }{\partial l} -
\frac{\partial \,}{\partial m} \Big)    - \frac{1}{(l-m)^{2 }} \left(1- \left( \gamma -\frac{1}{2}\right)^2\right) \,.
\end{eqnarray*}
The direct calculations  show that, if function $u$ solves the equation  ${\mathcal S}_{ch} u=0$, then the function $v=(l-m)u$ 
solves  the equation  ${\mathcal S}_{ch} ^*v=0$, and vice versa. Then Lemma~\ref{L1.2} completes the proof.
Lemma is proven. \hfill $\square$

\section{Proof of Theorem~\ref{T1}}
\label{S2a}

\setcounter{equation}{0}

Next we use Riemann function $R(l,m;a,b)$ and function $ E(x,t;x_0,t_0)$ defined by (\ref{E}) to 
complete the proof of Theorem~\ref{T1}, which gives the  fundamental solution with a support 
in the forward cone $D_+ (x_0,t_0)$, $x_0 \in {\mathbb R}^n$, $t_0 \in {\mathbb R}$,
and  the  fundamental solution with a support in the backward cone $D_- (x_0,t_0)$, $x_0 \in {\mathbb R}^n$, $t_0 \in {\mathbb R}$,
defined by (\ref{D+}) with plus and minus, respectively. 
\medskip

We present a proof   for   ${\mathcal E}_+(x,t;0,b) $ since for  
${\mathcal E}_-(x,t;0,b) $ it is similar.  
First, we note that the operator ${\mathcal S}$ is formally self-adjoint, ${\mathcal S}={\mathcal S}^*$. We must show that
\[
<{\mathcal E}_+, {\mathcal S} \varphi > = \varphi (0,b)\,, \qquad 
\mbox{\rm for every} \quad \varphi \in C_0^\infty ({\mathbb R}^2)\,.
\]
Since $E(x,t;0,b) $ is locally integrable in ${\mathbb R}^{2}$, 
this is equivalent to showing that
\begin{equation}
\label{Gel3.2}
\int \!\! \int_{{\mathbb R}^2} {\mathcal E}_+(x,t;0,b) {\mathcal S} \varphi (x,t)\, dx \, dt = \varphi (0,b), \quad 
\mbox{\rm for every} \quad \varphi \in C_0^\infty ({\mathbb R}^2).
\end{equation}
In the mean time $ {D(x,t)}/{D(l,m)} = (l-m)^{-1 } $ is the Jacobian of the transformation (\ref{Gel1.8}).
Hence the integral in the left-hand side of (\ref{Gel3.2}) is equal to   
\begin{eqnarray*}
&  &
\int \!\! \int_{{\mathbb R}^2} {\mathcal E}_+(x,t;0,b) {\mathcal S} \varphi (x,t)\, dx \, dt  
  =   
\int_b^\infty  dt \int_{e^{-t}-e^{-b}}^{e^{-b}-e^{-t}} E (x,t;0,b) {\mathcal S} \varphi (x,t)\, dx \\
& =  &
-  \int_{-e^{-b}}^{\infty }  
\int_{ -\infty }^{e^{-b}} R (l,m;e^{-b},-e^{-b})     \, dl \, dm   
 \Bigg\{ \frac{\partial^2  }{\partial l \, \partial m} - \frac{1}{2 (l-m)} \Big( \frac{\partial }{\partial l} -
\frac{\partial }{\partial m} \Big)  - \frac{1}{(l-m)^2 } M^2
\Bigg\}  \varphi .
\end{eqnarray*} 
We consider the first term of the right hand side, and integrate it by parts
\begin{eqnarray*}
&  &
\int_{-e^{-b}}^{\infty }  dm 
\int_{ -\infty }^{e^{-b}}   dl \, \,  R (l,m;e^{-b},-e^{-b})    
  \frac{\partial^2  }{\partial l \, \partial m}   \varphi \\
& = &
\int_{-e^{-b}}^{\infty }  R (e^{-b},m;e^{-b},-e^{-b})       
  \frac{\partial \varphi  }{  \partial m}   \Bigg|_{l= e^{-b}} \,   dm 
-  \Bigg[ -
\int_{ -\infty }^{e^{-b}}   \, dl \left( \frac{\partial}{\partial l } R (l,-e^{-b};e^{-b},-e^{-b})\right)     
   \varphi  \Bigg|_{m= -e^{-b}} \\
&  &
\hspace{2cm} - \int_{-e^{-b}}^{\infty }  dm 
\int_{ -\infty }^{e^{-b}}   \, dl \left( \frac{\partial^2 }{\partial l \,\partial m} R (l,m;e^{-b},-e^{-b})\right)      
  \varphi \Bigg] \\
& = &
-   \varphi ( e^{-b},    -e^{-b}) 
- \int_{-e^{-b}}^{\infty }   dm \, \left( \frac{\partial  }{  \partial m}  R (e^{-b},m;e^{-b},-e^{-b})  \right)      
   \varphi  (e^{-b},m)  \\
  &  &
+
\int_{ -\infty }^{e^{-b}}   \, dl \left( \frac{\partial}{\partial l } R (l,-e^{-b};e^{-b},-e^{-b})\right)   \,   
   \varphi  (l,-e^{-b}) + \int_{-e^{-b}}^{\infty }  dm 
\int_{ -\infty }^{e^{-b}}   \, dl \left( \frac{\partial^2 }{\partial l \,\partial m} R (l,m;e^{-b},-e^{-b})\right)      
  \varphi  .
\end{eqnarray*} 
Then, for the second term we obtain
\begin{eqnarray*}
&  &- \int_{-e^{-b}}^{\infty }  dm 
\int_{ -\infty }^{e^{-b}} \, dl \, R (l,m;e^{-b},-e^{-b})       
\frac{1}{2 (l-m)} \Big( \frac{\partial }{\partial l} -
\frac{\partial }{\partial m} \Big)   \varphi \\
&  =   &
- 
\int_{ -\infty }^{e^{-b}} R (l,-e^{-b};e^{-b},-e^{-b})      
\frac{1}{2 (l+e^{-b})}   \varphi   (l,-e^{-b})  \, dl\\
& &
- \int_{-e^{-b}}^{\infty }  R (e^{-b},m;e^{-b},-e^{-b})       
\frac{1}{2 (e^{-b}-m)}   \varphi ( e^{-b},m) \,  dm \\
&   &
- \int_{-e^{-b}}^{\infty }  dm 
\int_{ -\infty }^{e^{-b}}  dl  \,  \frac{1}{ (l-m)^2}  R (l,m;e^{-b},-e^{-b})  \varphi   (l,m) \,   
\\
&  &
+ \int_{-e^{-b}}^{\infty }  dm 
\int_{ -\infty }^{e^{-b}}  \, dl \,     
\Big[ \frac{1}{2 (l-m)} \Big( \frac{\partial }{\partial l} -
\frac{\partial }{\partial m} \Big)  R (l,m;e^{-b},-e^{-b})  \Big] \varphi   (l,m)\,.
\end{eqnarray*} 
Using  properties of the Riemann function we derive
\begin{eqnarray*}
&  &
\int_{-e^{-b}}^{\infty }  dm 
\int_{ -\infty }^{e^{-b}} R (l,m;e^{-b},-e^{-b})     \, dl \,   
 \Bigg\{ \frac{\partial^2  }{\partial l \, \partial m} - \frac{1}{2 (l-m)} \Big( \frac{\partial }{\partial l} -
\frac{\partial }{\partial m} \Big)  - \frac{1}{(l-m)^2 } M^2
\Bigg\}  \varphi \\
& = &
-   \varphi ( e^{-b},    -e^{-b}) 
- \int_{-e^{-b}}^{\infty }  \left( \frac{\partial  }{  \partial m}  R (e^{-b},m;e^{-b},-e^{-b})  \right)      
   \varphi  (e^{-b},m) \,   dm \\
  &  &
+
\int_{ -\infty }^{e^{-b}}    \left( \frac{\partial}{\partial l } R (l,-e^{-b};e^{-b},-e^{-b})\right)   \,   
   \varphi  (l,-e^{-b}) \, dl\\
 &  & - 
\int_{ -\infty }^{e^{-b}} R (l,-e^{-b};e^{-b},-e^{-b})      
\frac{1}{2 (l+e^{-b})}   \varphi   (l,-e^{-b})  \, dl \\
& &
- \int_{-e^{-b}}^{\infty }  R (e^{-b},m;e^{-b},-e^{-b})       
\frac{1}{2 (e^{-b}-m)}   \varphi ( e^{-b},m) \,  dm \\
& = &
-   \varphi ( e^{-b},    -e^{-b}) \,.
\end{eqnarray*} 
Theorem is proven. \hfill $\square$

\section{Application to the Cauchy Problem: Source Term and $n=1$}
\label{S3}
\setcounter{equation}{0}

Consider now the Cauchy problem for the equation (\ref{Int_3}) 
with vanishing initial data (\ref{Int_4}).
For every  $(x,t) \in D_+ (0,b)$ one has $   e^{-t}- e^{-b}  \le x \le  e^{-b}- e^{-t}$, so that 
\begin{eqnarray*}
E(x,t;0,b) 
& = &
(4e^{-b-t })^{iM} \Big((e^{-t }+e^{-b})^2 -  x  ^2\Big)^{-\frac{1}{2}-iM    } 
F\Big(\frac{1}{2}+iM   ,\frac{1}{2}+iM  ;1; 
\frac{ ( e^{-b}-e^{-t })^2 - x  ^2 }{( e^{-b}+e^{-t })^2 - x  ^2 } \Big) .
\end{eqnarray*}
The coefficient of the  equation (\ref{K_G_linear}) is independent of $x$, therefore ${\mathcal E}_+ (x,t;$ $y,b)$ $= {\mathcal E}_+ (x-y,t;0,b)$.
Using the fundamental solution from Theorem \ref{T1} one can write the convolution 
\begin{equation}
\label{convolution}
u(x,t)
 = 
\int_{ -\infty }^{\infty } \! \int_{ -\infty }^{\infty } {\mathcal E}_+ (x,t;y,b)f(y,b)\, db\, dy 
 = 
\int_{ 0}^{t}\! db \!\int_{ -\infty }^{\infty } {\mathcal E}_+ (x-y,t;0,b) f(y,b) \,dy , 
\end{equation}
which is well-defined since supp$f \subset \{ t \geq 0 \} $. 
 Then according to the definition of the distribution ${\mathcal E}_+ $ 
we obtain the statement of the Theorem~\ref{T1.1}. Thus,  Theorem~\ref{T1.1} is proven.

\begin{remark} 
The argument of the hypergeometric function is nonnegative and bounded, 
\[
0 \leq  \frac{ ( e^{-b}-e^{-t })^2 - z  ^2 }{( e^{-b}+e^{-t })^2 - z  ^2 }  < 1 \quad \mbox{  for all} 
\quad b\in (0,t),\,\, z \in (  e^{-t}- e^{-b} ,e^{-b}- e^{-t})\,.
 \]
\end{remark}
 The following corollary is a manifestation of the  time-speed transformation principle introduced in \cite{YagTricomi}. It implies  
the existence of an operator transforming the solutions of the  Cauchy problem for the string equation
to the solutions of the  Cauchy problem for the inhomogeneous equation with time-dependent speed of propagation.
One may think of this transformation as  a ``two-stage'' Duhamel's principal,
but unlike  the last one, it reduces the equation with the time-dependent speed of propagation to 
the one with the speed of propagation independent of time. 

\begin{corollary}
\label{C2}
The solution $u= u (x,t)$ of the Cauchy problem (\ref{Int_3})-(\ref{Int_4})  can be represented by (\ref{1.29}), 
where the functions $v (x,t;\tau ) := \frac{1}{2} ( f(x+t,\tau ) +   f(x - t,\tau ))$,  $\tau \in [0,\infty )$,
form a one-parameter family  
of  solutions to the Cauchy problem for the string equation, that is,
$v_{tt} - v_{xx} =0$,\, $ v (x,0;\tau ) = f(x,\tau )$, \, $v_t (x,0;\tau ) = 0$.
\end{corollary}
\medskip

\noindent
{\bf Proof.} From the convolution  (\ref{convolution}) we derive
\begin{eqnarray*} 
u(x,t)  
& = &
 \int_{ 0}^{t} db \int_{  - (e^{-b}- e^{-t}) }^{ e^{-b}- e^{-t}} dz \, f(z+x,b)  
(4e^{-b-t })^{iM} \Big((e^{-t }+e^{-b})^2 -  z  ^2\Big)^{-\frac{1}{2}-iM    } \\
&  &
\times F\Big(\frac{1}{2}+iM   ,\frac{1}{2}+iM  ;1; 
\frac{ ( e^{-b}-e^{-t })^2 - z  ^2 }{( e^{-b}+e^{-t })^2 - z  ^2 } \Big)  \\
& = &
2\int_{ 0}^{t} db \int_{ 0  }^{ e^{-b}- e^{-t}} dz \, \frac{1}{2} \{f(x+z ,b) +  f(x-z ,b)\}
(4e^{-b-t })^{iM} \Big((e^{-t }+e^{-b})^2 -  z  ^2\Big)^{-\frac{1}{2}-iM    } \\
&  &
\times F\Big(\frac{1}{2}+iM   ,\frac{1}{2}+iM  ;1; 
\frac{ ( e^{-b}-e^{-t })^2 - z  ^2 }{( e^{-b}+e^{-t })^2 - z  ^2 } \Big) \,.
\end{eqnarray*}
Corollary is proven. \hfill $\square$
\bigskip

\noindent
{\bf Some Properties of the Function  ${\mathbf E(x,t;y,b)}$.}
\smallskip

\noindent In this section we collect some elementary auxiliary formulas in order 
to make the proofs of main theorems more transparent. 

\begin{proposition}
\label{P_E}
Let $E(x,t;x_0,t_0) $ be function defined by (\ref{E}). One has
\begin{eqnarray}
\label{E_0}
E(x,t;y,b) 
& = &  E(y,b;x,t),\,\, \\
\label{E_1}
E(x,t;y,b) 
 =   
E(x-y,t;0,b)  & , & 
E(x,t;0,b) 
  =   
E(-x,t;0,b),\\
\label{E_2a}
E(x ,t;0,-\ln (x+e^{-t })) 
& = &
\frac{1}{2}\frac{1}{\sqrt{e^{-t }}\sqrt{x+e^{-t }}},\\ 
\label{E_3a}
\frac{\partial }{\partial b}\Big(  e^{-b} E(e^{-b}- e^{-t},t;0,b) \Big) 
& = &
- \frac{1}{4}e^{t/2}e^{-b/2}, \\
\label{E_4}
\frac{\partial }{\partial b} \Big( b e^{-b}E(-e^{-b}+ e^{-t},t;0,b) \Big) 
&= &
\frac{\partial }{\partial b} \Big( b e^{-b}E( e^{-b}- e^{-t},t;0,b) \Big) 
 = 
\frac{1}{4}e^{t/2}e^{-b/2}(2-b),\\ 
\label{E_5} 
\lim_{y\rightarrow  x+ e^{-b}- e^{-t}} \frac{\partial }{\partial x}E(x-y,t;0,b)  
& = &
-\frac{1}{16}(1+4M^2)e^{t/2} e^{b/2}(e^t-e^b),\\
\label{E_6}
\lim_{y\rightarrow  x -  e^{-b}+ e^{-t}} \frac{\partial }{\partial x}E(x-y,t;0,b) 
& = &
\frac{1}{16}(1+4M^2)e^{t/2} e^{b/2}(e^t-e^b), \\ 
\label{E_7}
\left[  \frac{\partial }{\partial b}   E(x,t;0,b) \right]_{b= -\ln(x+ e^{-t})} 
 & = &
\frac{1}{16} e^t \frac{ 4+e^tx(1+4M^2) }{\sqrt{1+e^tx}} , 
\end{eqnarray}
\begin{eqnarray} 
\label{E_8} 
\hspace{-0.5cm}\Bigg[\frac{\partial }{\partial b} E(x,t;0,b) \Bigg]_{b=0}
\!\! &\!\!  = \!\! &
 -(4e^{-t })^{iM} \big((e^{-t }+1)^2 - x^2\big)^{ -iM    } \frac{1}{ [(1-e^{ -t} )^2 -  x^2]\sqrt{(e^{-t }+1)^2 - x^2} } \nonumber \\
&   &
\times  \Bigg[  \big(  e^{-t} -1 - iM(e^{ -2t} -      1 -  x^2) \big) 
F \Big(\frac{1}{2}+iM   ,\frac{1}{2}+iM  ;1; \frac{ ( 1-e^{-t })^2 -x^2 }{( 1+e^{-t })^2 -x^2 }\Big)  \nonumber \\
&  &
\hspace{0.5cm}  +   \big( 1-e^{-2 t}+  x^2 \big)\Big( \frac{1}{2}-iM\Big)
F \Big(-\frac{1}{2}+iM   ,\frac{1}{2}+iM  ;1; \frac{ ( 1-e^{-t })^2 -x^2 }{( 1+e^{-t })^2 -x^2 }\Big) \Bigg].
\end{eqnarray}
\end{proposition}
\medskip

\noindent
{\bf Proof.} The properties (\ref{E_0}),(\ref{E_1}),   and (\ref{E_2a}) are evident. 
To prove (\ref{E_3a}) and (\ref{E_4}) we write
\[
E(e^{-b}- e^{-t},t;0,b)  
  =   
(4e^{-b-t })^{-\frac{1}{2}} ,
\] 
that implies both (\ref{E_3a}) and (\ref{E_4}). 
To prove (\ref{E_5}) we denote 
\begin{eqnarray*}
z:= \frac{ ( e^{-b}-e^{-t })^2 -(x- y )^2 }{( e^{-b}+e^{-t })^2 -(x- y )^2 } \, ,
\end{eqnarray*}
so that,
\begin{eqnarray*}
\frac{\partial z }{\partial x} 
& = &
 -2(x- y ) \frac{   4 e^{-b}e^{-t } }{[( e^{-b}+e^{-t })^2 -(x- y )^2 ]^2}\,, \qquad \frac{\partial z }{\partial b} 
= -\frac{4  e^{-b}e^{-t } \left(-e^{-2 t}+e^{-2 b}+  (x-y)^2\right) }{\left[(e^{- b}+e^{- t})^2- (x-y)^2\right]^2}\,.
\end{eqnarray*}
Then we obtain
\begin{eqnarray}   
\label{paE}
\frac{\partial }{\partial x} E(x-y,t;0,b) 
& = &
(4e^{-b-t })^{iM}  \Bigg[ (-2(x - y))(-\frac{1}{2}-iM  )\Big((e^{-t }+e^{-b})^2 - (x - y)^2\Big)^{-\frac{3}{2}-iM    } \nonumber \\
&  &
\times F\Big(\frac{1}{2}+iM   ,\frac{1}{2}+iM  ;1; z \Big)  \nonumber \\
 &  &
+ \Big((e^{-t }+e^{-b})^2 - (x - y)^2\Big)^{-\frac{1}{2}-iM    } 
\Big( \frac{\partial z }{\partial x}   \Big) F_z\Big(\frac{1}{2}+iM   ,\frac{1}{2}+iM  ;1; z \Big)\Bigg] 
\,.
\end{eqnarray} 
It is easily seen that
\[
\lim_{y\to x+ e^{-b}- e^{-t}} z= 0, \qquad  
 \lim_{y\to x+ e^{-b}- e^{-t}} \frac{\partial z}{\partial x}  =   \frac{ 1}{ 2  }( e^{ t} - e^{b }),
\]
while according to (20)  \cite[Sec.2.8 v.1]{B-E}, we have
\begin{eqnarray} 
\label{Fderiv} 
\lim_{y\to x+ e^{-b}- e^{-t}}  \partial _z   F \left(\frac{1}{2}+iM,\frac{1}{2}+iM;1;z\right)  
& = &
\frac{\Gamma^2 (\frac{3}{2}+iM)}{\Gamma^2 (\frac{1}{2}+iM) } 
  =  
\left(\frac{1}{2}+iM \right)^2\,.
\end{eqnarray} 
Consequently, from (\ref{paE}) we obtain
\begin{eqnarray*}   
&  &
\lim_{y\to x+ e^{-b}- e^{-t}} \frac{\partial }{\partial x} E(x-y,t;0,b) \\
& = &
(4e^{-b-t })^{iM} \Bigg[   2(e^{-b}- e^{-t}) \left(-\frac{1}{2}-iM  \right)\Big( 4e^{-t } e^{-b}   \Big)^{-\frac{3}{2}-iM    } 
+ \Big(4e^{-t } e^{-b}\Big)^{-\frac{1}{2}-iM    } 
\frac{ 1}{ 2  }( e^{ t} - e^{b })   \left(\frac{1}{2}+iM \right)^2\Bigg] \\
& = &
-  e^{ \frac{t}{2} } e^{\frac{b}{2} }  \frac{ 1}{ 16  }( e^{ t} - e^{b })  \Big[  1+ 4  M^2  \Big] 
\,.
\end{eqnarray*} 
The proof of (\ref{E_6}) is similar. To prove (\ref{E_7}) 
we write
\begin{eqnarray}
\label{AAA} 
\hspace{-0.5cm}& &
\frac{\partial }{\partial b} E(x,t;0,b) \nonumber \\
\hspace{-0.5cm}& = &
(-iM)    (4e^{-b-t })^{iM}  \Big((e^{-t }+e^{-b})^2 - x^2\Big)^{-\frac{1}{2}-iM    }
F\Big(\frac{1}{2}+iM   ,\frac{1}{2}+iM  ;1; 
\frac{ ( e^{-b}-e^{-t })^2 -x^2 }{( e^{-b}+e^{-t })^2 -x^2 } \Big)  \nonumber  \\
\hspace{-0.5cm}&   &
-\Big( \frac{1}{2}+iM\Big) (4e^{-b-t })^{iM}  (-e^{-b}) 2(e^{-t }+e^{-b})    \Big((e^{-t }+e^{-b})^2 - x^2\Big)^{-\frac{3}{2}-iM    }  \nonumber  \\
\hspace{-0.5cm}&  &
\times F\Big(\frac{1}{2}+iM   ,\frac{1}{2}+iM  ;1; 
\frac{ ( e^{-b}-e^{-t })^2 -x^2 }{( e^{-b}+e^{-t })^2 -x^2 } \Big)  \nonumber  \\
\hspace{-0.5cm}&   &
 +   (4e^{-b-t })^{iM} \Big((e^{-t }+e^{-b})^2 - x^2\Big)^{-\frac{1}{2}-iM    } 
\frac{\partial }{\partial b}\Bigg(F\Big(\frac{1}{2}+iM   ,\frac{1}{2}+iM  ;1; 
\frac{ ( e^{-b}-e^{-t })^2 -x^2 }{( e^{-b}+e^{-t })^2 -x^2 } \Big)  \Bigg).
\end{eqnarray}
On the other hand (\ref{Fderiv}) and (\ref{AAA}) imply
\begin{eqnarray*}  
\hspace{-0.5cm} &  &
\left[ \frac{\partial }{\partial b} E(x,t;0,b) \right]_{b= -\ln (x+e^{-t })} \\
\hspace{-0.5cm} & = &
 (4e^{-b-t })^{iM}  \Big(4e^{-b-t }\Big)^{-\frac{1}{2}-iM    }  \Bigg[(-iM)  
+ 2\Big( \frac{1}{2}+iM\Big)   e^{-b} (e^{-t }+e^{-b})    \Big(4e^{-b-t }\Big)^{-1    }   \\
\hspace{-0.5cm} &   &
\hspace{5cm}  +  
\frac{\partial }{\partial b}\Bigg(F\Big(\frac{1}{2}+iM   ,\frac{1}{2}+iM  ;1; 
\frac{ ( e^{-b}-e^{-t })^2 -x^2 }{( e^{-b}+e^{-t })^2 -x^2 } \Big)  \Bigg)\Bigg]_{b= -\ln (x+e^{-t })}\\
\hspace{-0.5cm} & = &
\frac{1}{2}e^{ \frac{1}{2}b+\frac{1}{2}t }  \Bigg[(-iM)  
+ \frac{1}{2}\Big( \frac{1}{2}+iM\Big)   (e^{-t }+e^{-b})      e^{ t }   
-\frac{4  e^{-b}e^{-t } \left( e^{-2 b}-e^{-2 t}+  x^2 \right) }
{\left[(e^{- b}+e^{- t})^2- x^2\right]^2}\Big( \frac{1}{2}+iM\Big)^2\Bigg]_{b= -\ln (x+e^{-t })}\\
\hspace{-0.5cm} & = &
\frac{1}{2}\frac{1}{\sqrt{1+xe^{t }}}e^{  t }  \Bigg[ -iM   
+ \frac{1}{2}\Big( \frac{1}{2}+iM\Big)   (2+x e^{ t } )       
-\frac{   x e^{t }}{ 2   }\Big( \frac{1}{2}+iM\Big)^2\Bigg]\\
\hspace{-0.5cm} & = &
\frac{1}{2}\frac{1}{\sqrt{1+xe^{t }}}e^{  t }  \frac{4+x e^{ t }(1+4M^2)}{8}\,.
\end{eqnarray*} 
Thus, (\ref{E_7})  is proven. 
To prove (\ref{E_8})  
 we  appeal to (23)\cite[v.1, Sec.2.8]{B-E}, that reads with \,$a=b=\frac{1}{2}+iM$, $c=1$,
\begin{eqnarray*} 
F_z\Big(\frac{1}{2}+iM   ,\frac{1}{2}+iM  ;1; z\Big)& = & 
\frac{1}{z(1-z)}\Bigg\{    \Big[ -  \Big( \frac{1}{2}-iM\Big)
+ \Big( \frac{1}{2}+iM\Big)z\Big] F \Big(\frac{1}{2}+iM   ,\frac{1}{2}+iM  ;1; z\Big)\\
&  &
\hspace{2cm} + \Big( \frac{1}{2}-iM\Big)F \Big(-\frac{1}{2}+iM   ,\frac{1}{2}+iM  ;1; z\Big)\Bigg\}\,.
\end{eqnarray*}
Then we plug the last relation in (\ref{AAA}) and obtain  
\begin{eqnarray*} 
&  &
\Bigg[\frac{\partial }{\partial b} E(x,t;0,b) \Bigg]_{b=0}\nonumber \\
& = &
-iM (4e^{-t })^{iM}  \Big((e^{-t }+1)^2 - x^2\Big)^{-\frac{1}{2}-iM    }
F\Big(\frac{1}{2}+iM   ,\frac{1}{2}+iM  ;1; 
\zeta_0 \Big)  \nonumber  \\
&   &
-\Big( \frac{1}{2}+iM\Big) (4e^{-t })^{iM}  (-1) 2(e^{-t }+1)    \Big((e^{-t }+1)^2 - x^2\Big)^{-\frac{3}{2}-iM    }  
 F\Big(\frac{1}{2}+iM   ,\frac{1}{2}+iM  ;1; 
\zeta_0  \Big)  \nonumber  \\
&   &
 +   (4e^{-t })^{iM} \Big((e^{-t }+1)^2 - x^2\Big)^{-\frac{1}{2}-iM    } 
\left[\frac{\partial }{\partial b}\left(F\Big(\frac{1}{2}+iM   ,\frac{1}{2}+iM  ;1; 
\zeta\Big)  \right) \right]_{b=0}\,,
\end{eqnarray*}
where  we have denoted
\[
\zeta := \frac{ ( e^{-b}-e^{-t })^2 -x^2 }{( e^{-b}+e^{-t })^2 -x^2 } \,, 
\qquad \zeta_0 := \frac{ ( 1-e^{-t })^2 -x^2 }{( 1+e^{-t })^2 -x^2 }\,, \qquad  
\zeta_0(1-\zeta_0)=\frac{ 4e^{-t }[( 1-e^{-t })^2 -x^2 ]}{[( 1+e^{-t })^2 -x^2]^2 }\,,
\]
with
\[
 \frac{\partial }{\partial b}\zeta 
  =  -\frac{4  e^{-b}e^{-t } \left( e^{-2 b}-e^{-2 t}+  x^2\right) }{\left[(e^{- b}+e^{- t})^2- x^2\right]^2} 
\,,\qquad \frac{\partial \zeta}{\partial b} \Big|_{b=0}  
  =    
 -\frac{4  e^{-t } \left( 1-e^{-2 t}+  x^2\right) }{\left[(1+e^{- t})^2- x^2\right]^2}  \,.
\]
Hence, due to (20) \cite[v.1, Sec.2.8]{B-E},  we obtain
\begin{eqnarray*}
&  &
\Bigg[\frac{\partial }{\partial b}F\Big(\frac{1}{2}+iM   ,\frac{1}{2}+iM  ;1; 
\frac{ ( e^{-b}-e^{-t })^2 -x^2 }{( e^{-b}+e^{-t })^2 -x^2 } \Big)   \Bigg]_{b=0}\\
& = &
 -\frac{  1-e^{-2 t}+  x^2 }{ ( 1-e^{-t })^2 -x^2}\Bigg\{    \Big[ -  \Big( \frac{1}{2}-iM\Big)
+ \Big( \frac{1}{2}+iM\Big)\zeta_0\Big] F \Big(\frac{1}{2}+iM   ,\frac{1}{2}+iM  ;1; \zeta_0\Big)\\
&  &
\hspace{6cm} + \Big( \frac{1}{2}-iM\Big)F \Big(-\frac{1}{2}+iM   ,\frac{1}{2}+iM  ;1; \zeta_0\Big)\Bigg\}\,.
\end{eqnarray*}
Hence,
\begin{eqnarray*} 
&  &
\Bigg[\frac{\partial }{\partial b} E(x,t;0,b) \Bigg]_{b=0}\nonumber \\
& = &
-iM (4e^{-t })^{iM}  \Big((e^{-t }+1)^2 - x^2\Big)^{-\frac{1}{2}-iM    }
F\Big(\frac{1}{2}+iM   ,\frac{1}{2}+iM  ;1; \zeta_0 \Big)  \nonumber  \\
&   &
-\Big( \frac{1}{2}+iM\Big) (4e^{-t })^{iM}  (-1) 2(e^{-t }+1)    \Big((e^{-t }+1)^2 - x^2\Big)^{-\frac{3}{2}-iM    }  \nonumber   F\Big(\frac{1}{2}+iM   ,\frac{1}{2}+iM  ;1; \zeta_0\Big)  \nonumber  \\
&   &
 +   (4e^{-t })^{iM} \Big((e^{-t }+1)^2 - x^2\Big)^{-\frac{1}{2}-iM    } \\
&  &
\times \Bigg[  -\frac{   1-e^{-2 t}+  x^2  }{  ( 1-e^{-t })^2 -x^2  }\Bigg\{    \Big[ -  \Big( \frac{1}{2}-iM\Big)
+ \Big( \frac{1}{2}+iM\Big)\zeta_0\Big] F \Big(\frac{1}{2}+iM   ,\frac{1}{2}+iM  ;1; \zeta_0\Big)\\
&  &
+ \Big( \frac{1}{2}-iM\Big)F \Big(-\frac{1}{2}+iM   ,\frac{1}{2}+iM  ;1; \zeta_0\Big)\Bigg\} \Bigg] \\
& = &
(4e^{-t })^{iM} \Big((e^{-t }+1)^2 - x^2\Big)^{-\frac{1}{2}-iM    }\Bigg[ F \Big(\frac{1}{2}+iM   ,\frac{1}{2}+iM  ;1; \zeta_0\Big)  \Big\{ -iM   + ( 1+2iM )   \frac{e^{-t }+1}{(e^{-t }+1)^2 - x^2}  \\
&   &
    - \frac{ 1-e^{-2 t} +  x^2 }{ ( 1-e^{-t })^2 -x^2 }\Big[ -  \Big( \frac{1}{2}-iM\Big)
+ \Big( \frac{1}{2}+iM\Big)\frac{ ( 1-e^{-t })^2 -x^2 }{( 1+e^{-t })^2 -x^2 }\Big] \Big\}\\
&  &
   -  \frac{ 1-e^{-2 t}+  x^2 }{ ( 1-e^{-t })^2 -x^2}\Big( \frac{1}{2}-iM\Big)F \Big(-\frac{1}{2}+iM   ,\frac{1}{2}+iM  ;1; \zeta_0\Big) \Bigg] .
\end{eqnarray*}
Finally,
\begin{eqnarray*} 
\Bigg[\frac{\partial }{\partial b} E(x,t;0,b) \Bigg]_{b=0} 
& = &
(4e^{-t })^{iM} \big((e^{-t }+1)^2 - x^2\big)^{-\frac{1}{2}-iM    } \\
&   &
\times  \Bigg[\frac{e^t- iM +e^{2 t}  \big(  iM(1 +  x^2) -1  \big)}{ e^t  \big(2+e^t  (x^2-1 ) \big)-1} 
F \Big(\frac{1}{2}+iM   ,\frac{1}{2}+iM  ;1; \zeta_0\Big) \\
&  &
\hspace{1cm}   -  \frac{ 1-e^{-2 t}+  x^2 }{ ( 1-e^{-t })^2 -x^2}\Big( \frac{1}{2}-iM\Big)
F \Big(-\frac{1}{2}+iM   ,\frac{1}{2}+iM  ;1; \zeta_0\Big) \Bigg]\,.
\end{eqnarray*}
The formula (\ref{E_8}) and, consequently, the proposition are proven. \hfill $\square$

\section{The Cauchy Problem:  Second Datum and $n=1$}
\label{S5}

\setcounter{equation}{0}

In this section we prove Theorem~\ref{T1.3} in the case of  $\varphi_0  (x)=0$. More precisely, we have to prove that 
the solution $u (x,t)$ of the Cauchy problem  (\ref{oneDphy01}) with  $\varphi_0  (x)=0$ and  $\varphi_1  (x)=\varphi (x)$
can be represented as follows
\begin{equation}
\label{3.1n}
u(x,t)  
 =   
 \int_{0}^{1 - e^{-t}} \,\Big[      \varphi   (x+ z)  +   \varphi   (x - z)    \Big] K_1(z,t) dz
 =   
   \int_{0}^{  1} \,\Big[      \varphi   (x+ \phi (t) s )  +   \varphi   (x - \phi (t) s)    \Big] K_1(\phi (t) s,t) \phi (t) ds ,
\end{equation}
where $\phi (t) =  1 - e^{-t}$. 
The proof of the theorem is splitted into several steps.

\begin{proposition}
\label{C3a}
The solution $u = u (x,t)$ of the Cauchy problem    (\ref{oneDphy01}) with  $\varphi_0  (x)=0$ and  $\varphi_1  (x)=\varphi (x)$
can be represented as follows
\begin{eqnarray}
\label{5.2}
u(x,t)
& = &
  \int_{ 0}^{t}   db \,  \Big[  \frac{1}{4}e^{t/2}e^{-b/2}(2-b) +\frac{1}{16} be^{-3b/2}e^{t/2}(e^b-e^t)(1+4M^2)\Big] \nonumber \\
&  &
\hspace{1cm} \times \Big[ \varphi  (x+ e^{-b}- e^{-t} ) +   \varphi  (x - e^{-b}+ e^{-t} ) \Big]\nonumber \\
&  &
+ \int_{ 0}^{t} db \int_{ x - (e^{-b}- e^{-t})}^{x+ e^{-b}- e^{-t}} dy \, \varphi (y)   b \Big[e^{-2b}\,
\left( \frac{\partial}{\partial y} \right)^2 E(x-y,t;0,b)  
- M^2 E(x-y,t;0,b) \Big] \, \,.
\end{eqnarray}
\end{proposition}
 \medskip

 \noindent
 {\bf Proof.} 
We look for the solution $u=u(x,t)$ of the form $
u(x,t)=w(x,t)+ t  \varphi (x)$. 
Then     (\ref{oneDphy01})  implies 
\begin{eqnarray*}
&  &
w_{tt} - e^{-2t}w_{xx} +M^2w = t e^{-2t}\varphi^{(2)} (x)-M^2t\varphi(x) , \qquad  w(x,0) = 0,\quad w_t(x,0)= 0 \,.
\end{eqnarray*}
We set $f(x,t)=  t e^{-2t}\varphi^{(2)} (x) -M^2t\varphi(x)$ and due to Theorem~\ref{T1.1} obtain  
\[
w(x,t)  =\widetilde{ w(x,t)} - M^2\int_{ 0}^{t} b \, db \int_{ x - (e^{-b}- e^{-t})}^{x+ e^{-b}- e^{-t}} dy \, \varphi (y)  
E(x-y,t;0,b) ,
\]
where
\begin{eqnarray*}
\widetilde{ w(x,t)}    
&  :=  &
 \int_{ 0}^{t} b e^{-2b}\, db \int_{ x - (e^{-b}- e^{-t})}^{x+ e^{-b}- e^{-t}} dy \, \varphi^{(2)} (y)  
E(x-y,t;0,b)  \,.
\end{eqnarray*}
Then we integrate by parts:
\begin{eqnarray*}
\widetilde{ w(x,t)}    
\!\! &\!\!  = \!\!  &\!\! 
 \int_{ 0}^{t} b e^{-2b}\, db \Bigg[  \varphi^{(1)} (x+ e^{-b}- e^{-t})  
E(-e^{-b}+ e^{-t},t;0,b)-  \varphi^{(1)} (x - e^{-b}+ e^{-t} ) E(e^{-b}- e^{-t},t;0,b) \Bigg]\\
&  &
-  \int_{ 0}^{t} b e^{-2b}\, db \int_{ x - (e^{-b}- e^{-t})}^{x+ e^{-b}- e^{-t}} dy \, \varphi^{(1)} (y)  
\frac{\partial}{\partial y} E(x-y,t;0,b) \,.
\end{eqnarray*}
But
\[
 \varphi^{(1)} (x+ e^{-b}- e^{-t} ) = -  e^{b} \frac{\partial }{\partial b}  \varphi  (x+ e^{-b}- e^{-t} ) \quad \mbox{\rm and} \quad
 \varphi^{(1)} (x - e^{-b}+ e^{-t} )=    e^{b} \frac{\partial }{\partial b}  \varphi  (x - e^{-b}+ e^{-t} )\,.
\]
Then, \, $ E(e^{-b}- e^{-t},t;0,b)= E(-e^{-b}+ e^{-t},t;0,b)$\, due to (\ref{E_1}), and we obtain 
\begin{eqnarray*}
\widetilde{ w(x,t)} 
& = &
 \int_{ 0}^{t} b e^{-2b}\, db \Bigg[   -  e^{b} \frac{\partial }{\partial b}  \varphi  (x+ e^{-b}- e^{-t} ) 
-  e^{b} \frac{\partial }{\partial b}  \varphi  (x - e^{-b}+ e^{-t} ) \Bigg] E(e^{-b}- e^{-t},t;0,b) \\
&  &
-  \int_{ 0}^{t} b e^{-2b}\, db \int_{ x - (e^{-b}- e^{-t})}^{x+ e^{-b}- e^{-t}} dy \, \varphi^{(1)} (y)  
\frac{\partial}{\partial y} E(x-y,t;0,b) \,.
\end{eqnarray*}
One more integration by parts leads to
\begin{eqnarray*}
\widetilde{ w(x,t)} 
& = &
-  2t e^{-t}      \varphi  (x )  E(0,t;0,t)   \\
&  &
+ \int_{ 0}^{t} \, db \Big(   \varphi  (x+ e^{-b}- e^{-t} ) 
+ \varphi  (x - e^{-b}+ e^{-t} ) \Big)\frac{\partial }{\partial b} \Big( b e^{-b}E(e^{-b}- e^{-t},t;0,b) \Big)\\
&  &
-  \int_{ 0}^{t} b e^{-2b}\, db \int_{ x - (e^{-b}- e^{-t})}^{x+ e^{-b}- e^{-t}} dy \, \varphi^{(1)} (y)  
\frac{\partial}{\partial y} E(x-y,t;0,b) \,.
\end{eqnarray*}
Since  \, $ E( 0 ,t;0,t)=e^{t}/2 $ \,  we use (\ref{E_4})  of Proposition \ref{P_E} 
to derive the next representation 
\begin{eqnarray*}
\widetilde{ w(x,t)} + t      \varphi  (x )   
& = &
  \int_{ 0}^{t}     \frac{1}{4}e^{t/2}e^{-b/2}(2-b)\Big( \varphi  (x+ e^{-b}- e^{-t} ) +   \varphi  (x - e^{-b}+ e^{-t} ) \Big)\, db \\
&  &
-  \int_{ 0}^{t} b e^{-2b}\, db \int_{ x - (e^{-b}- e^{-t})}^{x+ e^{-b}- e^{-t}} dy \, \varphi^{(1)} (y)  
\frac{\partial}{\partial y} E(x-y,t;0,b) \,.
\end{eqnarray*}
The integration by parts in the second term leads to
\begin{eqnarray*}
\widetilde{ w(x,t)} + t      \varphi  (x )  
& = &
  \int_{ 0}^{t}     \frac{1}{4}e^{t/2}e^{-b/2}(2-b)\Big[ \varphi  (x+ e^{-b}- e^{-t} ) +   \varphi  (x - e^{-b}+ e^{-t} ) \Big]\, db \\
&  &
- \int_{ 0}^{t} b e^{-2b}\, db \, \varphi (x+ e^{-b}- e^{-t})  
\left[ \frac{\partial}{\partial y} E(x-y,t;0,b) \right]_{y=x+ e^{-b}- e^{-t}}\\
&  &
+ \int_{ 0}^{t} b e^{-2b}\, db  \, \varphi (x - e^{-b}+ e^{-t} )  
\left[ \frac{\partial}{\partial y} E(x-y,t;0,b)\right]_{ y=x - e^{-b}+ e^{-t} }\\
&  &
+ \int_{ 0}^{t} b e^{-2b}\, db \int_{ x - (e^{-b}- e^{-t})}^{x+ e^{-b}- e^{-t}} dy \, \varphi (y)  
\left( \frac{\partial}{\partial y} \right)^2 E(x-y,t;0,b) 
\,.
\end{eqnarray*}
The application of  (\ref{E_5}) and  (\ref{E_6}) of Proposition \ref{P_E}  
and \, $\frac{\partial }{\partial y} E(x-y,t;0,b)=-\frac{\partial }{\partial x} E(x-y,t;0,b)  $ imply 
\begin{eqnarray*}
&  &
\widetilde{ w(x,t)} + t      \varphi  (x )  \\
& = &
  \int_{ 0}^{t}    \Big[  \frac{1}{4}e^{t/2}e^{-b/2}(2-b) +\frac{1}{16} be^{-3b/2}e^{t/2}(e^b-e^t)(1+4M^2)\Big]\Big[ \varphi  (x+ e^{-b}- e^{-t} ) +   \varphi  (x - e^{-b}+ e^{-t} ) \Big]\, db \\
&  &
+ \int_{ 0}^{t} b e^{-2b}\, db \int_{ x -  e^{-b}+ e^{-t} }^{x+ e^{-b}- e^{-t}} dy \, \varphi (y)  
\left( \frac{\partial}{\partial y} \right)^2 E(x-y,t;0,b) \,.
\end{eqnarray*}
Thus, for the function $u(x,t)=  w(x,t) + t      \varphi  (x ) $ we have obtained
\begin{eqnarray}
&  &
u(x,t) \nonumber \\
& = &
  \int_{ 0}^{t}    \Big[  \frac{1}{4}e^{t/2}e^{-b/2}(2-b) +\frac{1}{16} be^{-3b/2}e^{t/2}(e^b-e^t)(1+4M^2)\Big]\Big[ \varphi  (x+ e^{-b}- e^{-t} ) +   \varphi  (x - e^{-b}+ e^{-t} ) \Big]\, db \nonumber \\
&  &
+ \int_{ 0}^{t} b e^{-2b}\, db \int_{ x -  e^{-b}+ e^{-t}}^{x+ e^{-b}- e^{-t}} dy \, \varphi (y)  
\left( \frac{\partial}{\partial y} \right)^2 E(x-y,t;0,b) \nonumber \\
&  &
- M^2\int_{ 0}^{t} b \, db \int_{ x -  e^{-b}+ e^{-t}}^{x+ e^{-b}- e^{-t}} dy \, \varphi (y)  
E(x-y,t;0,b)  \, \,.
\end{eqnarray}
The proposition is proven.
\hfill $\square$

\begin{corollary}
\label{C5.3}
The solution $u= u (x,t)$ of the Cauchy problem  (\ref{oneDphy01}) with  $\varphi_0  (x)=0$ and  $\varphi_1  (x)=\varphi (x)$
can be represented as follows
\begin{eqnarray}
&  &
u(x,t) \nonumber \\
& = &
  \int_{ 0}^{t}    \Big[  \frac{1}{4}e^{t/2}e^{-b/2}(2-b) +\frac{1}{16} be^{-3b/2}e^{t/2}(e^b-e^t)(1+4M^2)\Big]\Big[ \varphi  (x+ e^{-b}- e^{-t} ) +   \varphi  (x - e^{-b}+ e^{-t} ) \Big]\, db \nonumber \\
&  &
+ \int_{ 0}^{t} db \int_{ 0 }^{e^{-b}- e^{-t}}  dz \, \Big[ \varphi (x-z) +  \varphi (x+z)\Big]
b \left[e^{-2b}\, \left( \frac{\partial}{\partial z} \right)^2 E(z,t;0,b) - M^2E(z,t;0,b)\right]\,,
\end{eqnarray}
as well as by (\ref{3.1n}),
where 
\begin{eqnarray}
\label{3.3}
K_1(z,t) 
&  =  &
  \Big[  \frac{1}{4}e^{t/2}\big(2+\ln (z+e^{-t })\big) 
+\frac{1}{16}(1+4M^2) e^{3t/2}z\ln (z+e^{-t })\Big]\frac{1}{\sqrt{z+e^{-t }}} \nonumber \\
&  &
+ \int_{   0}^{ -\ln (z+e^{-t })}
 b \left[e^{-2b}\left( \frac{\partial}{\partial z} \right)^2 E(z,t;0,b) - M^2E(z,t;0,b) \right] db\,.
 \end{eqnarray}
\end{corollary}
\medskip

\noindent
{\bf Proof of Corollary.} By means of  the statement (\ref{5.2}) of  Proposition~\ref{C3a} and (\ref{E_1}) we obtain
\begin{eqnarray*}
u(x,t)
& = &
  \int_{ 0}^{t}  db   \, \Big[  \frac{1}{4}e^{t/2}e^{-b/2}(2-b) +\frac{1}{16} be^{-3b/2}e^{t/2}(e^b-e^t)(1+4M^2)\Big]\\
&  &
\hspace{1cm} \times \Big[ \varphi  (x+ e^{-b}- e^{-t} ) +   \varphi  (x - e^{-b}+ e^{-t} ) \Big]  \\
&  &
+ \int_{ 0}^{t} db \int_{ 0 }^{-(e^{-b}- e^{-t})} (-1) dz \, \varphi (x-z)  
\left[b e^{-2b}\, \left( \frac{\partial}{\partial z} \right)^2 E(z,t;0,b) - M^2bE(z,t;0,b)\right]  \\
&  &
+ \int_{ 0}^{t} db \int_{ - (e^{-b}- e^{-t})}^{0} dz \, \varphi (x+z)  
\left[b e^{-2b}\, \left( \frac{\partial}{\partial z} \right)^2 E(z,t;0,b) - M^2bE(z,t;0,b)\right]\,.
\end{eqnarray*}
To   prove (\ref{3.1n}) with $K_1(z,t) $ defined by (\ref{3.3}) we apply   (\ref{E_1})  and write
\begin{eqnarray*}
u(x,t)
& = &
  \int_{ 0}^{t}   db  \, \Big[  \frac{1}{4}e^{t/2}e^{-b/2}(2-b) +\frac{1}{16} be^{-3b/2}e^{t/2}(e^b-e^t)(1+4M^2)\Big]\\
&  &
\hspace{1cm} \times \Big[ \varphi  (x+ e^{-b}- e^{-t} ) +   \varphi  (x - e^{-b}+ e^{-t} ) \Big]\nonumber \\
&  &
+ \int_{ 0}^{t} \, db \int_{   0}^{  e^{-b}- e^{-t}} dz \, \Big[\varphi (x+z)+  \varphi (x-z) \Big]
\left[ b e^{-2b}\left( \frac{\partial}{\partial z} \right)^2 E(z,t;0,b) - M^2bE(z,t;0,b) \right].
\end{eqnarray*}
Next we make change $z= e^{-b }- e^{-t }$, \, $dz= -e^{-b }db$,\, $db= -(z+e^{-t })^{-1}dz$, and $b=-\ln (z+e^{-t })$ in the first integral:
\begin{eqnarray*}
&  &
  \int_{ 0}^{t}    \Big[  \frac{1}{4}e^{t/2}e^{-b/2}(2-b) +\frac{1}{16} be^{-3b/2}e^{t/2}(e^b-e^t)(1+4M^2)\Big]\Big[ \varphi  (x+ e^{-b}- e^{-t} ) +   \varphi  (x - e^{-b}+ e^{-t} ) \Big]\, db  \\
& = &
 \int_{0}^{ 1-e^{-t}}    \Big[  \frac{1}{4}e^{t/2}\big(2+\ln (z+e^{-t })\big) 
+\frac{1}{16}(1+4M^2) e^{3t/2}z\ln (z+e^{-t })\Big]\frac{1}{\sqrt{z+e^{-t }}}\\
&  &
\hspace{2cm}\times \Big[ \varphi  (x+ z ) +   \varphi  (x - z ) \Big]\,dz \,.
\end{eqnarray*}
Then
\begin{eqnarray*}
u(x,t)
& = &
  \int_{0}^{ 1-e^{-t}}    \Big[  \frac{1}{4}e^{t/2}\big(2+\ln (z+e^{-t })\big) 
+\frac{1}{16}(1+4M^2) e^{3t/2}z\ln (z+e^{-t })\Big]\frac{1}{\sqrt{z+e^{-t }}}\\
&  &
\hspace{2cm}\times \Big[ \varphi  (x+ z ) +   \varphi  (x - z ) \Big]\,dz \\
&  &
+ \int_{ 0}^{t} \, db \int_{   0}^{  e^{-b}- e^{-t}} dz \, \Big[\varphi (x+z)+  \varphi (x-z) \Big]
 b \left[e^{-2b}\left( \frac{\partial}{\partial z} \right)^2 E(z,t;0,b) - M^2E(z,t;0,b) \right] .
\end{eqnarray*}
In the last integral we change the order of integration,
\begin{eqnarray*}
&  &
u(x,t) \\
& = &
\int_{0}^{ 1-e^{-t}}    \Big[  \frac{1}{4}e^{t/2}\big(2+\ln (z+e^{-t })\big) 
+\frac{1}{16}(1+4M^2) e^{3t/2}z\ln (z+e^{-t })\Big]\frac{1}{\sqrt{z+e^{-t }}} \Big[ \varphi  (x+ z ) +   \varphi  (x - z ) \Big]\,dz \\
&  &
+ \int_{ 0}^{1-e^{-t}} \, dz  \Big[\varphi (x+z)+  \varphi (x-z) \Big]\int_{   0}^{ -\ln (z+e^{-t })} db \,
 b \left[e^{-2b}\left( \frac{\partial}{\partial z} \right)^2 E(z,t;0,b) - M^2E(z,t;0,b) \right],
\end{eqnarray*}
and obtain (\ref{3.1n}), where  $K_1(z,t)   $ is defined by (\ref{3.3}). 
Corollary is proven.
\hfill $\square$
\medskip

\noindent
{\bf Proof of Theorem \ref{T1.3} with $\varphi _0=0$.}  
The next lemma completes the proof of Theorem~\ref{T1.3}.
\begin{lemma} 
The kernel   $K_1(z,t)   $ defined by (\ref{3.3}) coincides with one given in Theorem~\ref{T1.3}. 
\end{lemma}
\medskip

\noindent
{\bf Proof.} We have due to Lemma~\ref{L1.2}, (\ref{E_1}), and by integration by parts 
\begin{eqnarray*}
&  &
\int_{   0}^{ -\ln (z+e^{-t })} 
 b \left[e^{-2b}\left( \frac{\partial}{\partial z} \right)^2 E(z,t;0,b) - M^2 E(z,t;0,b) \right] \, db 
  =  
\int_{   0}^{ -\ln (z+e^{-t })}  b  \Big( \frac{\partial }{\partial b} \Big)^2 E(z ,t;0,b)  db\\
& = &
(-\ln (z+e^{-t })) \Big[ \frac{\partial }{\partial b}  E(z ,t;0,b) \Big]_{b= -\ln (z+e^{-t })}
-  E(z ,t;0,-\ln (z+e^{-t }))  +  E(z ,t;0,0) \,.
 \end{eqnarray*} 
On the other hand,  (\ref{E_2a}) and (\ref{E_7}) of Proposition \ref{P_E} 
imply 
 \begin{eqnarray*}
&  &
\int_{   0}^{ -\ln (z+e^{-t })} db \,
 b \left[e^{-2b}\left( \frac{\partial}{\partial z} \right)^2 E(z,t;0,b) - M^2 E(z,t;0,b) \right] \\
& = &
-\ln (z+e^{-t }) \frac{1}{16}\frac{4+e^tz(1+4M^2)}{\sqrt{e^{-t }}\sqrt{z+e^{-t }}}
-  \frac{1}{2}\frac{1}{\sqrt{e^{-t }}\sqrt{z+e^{-t }}}  +  E(z ,t;0,0) .
 \end{eqnarray*}
Thus, for the kernel   $K_1(z,t)   $ defined by (\ref{3.3}) we have
 \begin{eqnarray*}
K_1(z,t) 
& = &
   \Big[  \frac{1}{4}e^{t/2}\big(2+\ln (z+e^{-t })\big) 
+\frac{1}{16}(1+4M^2) e^{3t/2}z\ln (z+e^{-t })\Big]\frac{1}{\sqrt{z+e^{-t }}}\\
&  &
 -\ln (z+e^{-t }) \frac{1}{16}\frac{4+e^tz(1+4M^2)}{\sqrt{e^{-t }}\sqrt{z+e^{-t }}}
-  \frac{1}{2}\frac{1}{\sqrt{e^{-t }}\sqrt{z+e^{-t }}}  +  E(z ,t;0,0) \\ 
& = &
 e^{t/2}  \Big[  \frac{1}{4}\big(2+\ln (z+e^{-t })\big) 
+\frac{1}{16}(1+4M^2) e^{t}z\ln (z+e^{-t })\\
&  &
 - \ln (z+e^{-t }) \frac{1}{16}(4+e^tz(1+4M^2)) 
-  \frac{1}{2}   \Big]\frac{1}{\sqrt{z+e^{-t }}} +  E(z ,t;0,0) \\ 
& = &
E(z ,t;0,0) \,.
 \end{eqnarray*}
The last line coincides with $K_1(z,t)$ of Theorem~\ref{T1.3}. Lemma is proven. \hfill $\square$

\section{The Cauchy Problem:  First Datum and $n=1$}
\label{S6}
\setcounter{equation}{0}

In this section we prove Theorem~\ref{T1.3} in the case of  $\varphi_1  (x)=0$. Thus, we have to prove the representation given by   Theorem~\ref{T1.3} for
the solution $u=u (x,t)$ of the Cauchy problem (\ref{oneDphy01}) with $\varphi_1  (x)=0$, that is equivalent to 
\[
u(x,t)  
  =    
\frac{1}{2} e ^{ \frac{t}{2}} \Big[     \varphi_0   (x+ 1- e^{-t})   
+   \varphi_0   (x -  1+ e^{-t} )\Big] 
+ \, \int_{ 0}^{1} \big[ 
\varphi_0   (x - \phi (t)s)  
+     \varphi_0   (x  + \phi (t)s)  \big] K_0(\phi (t)s,t)\phi (t)\,  ds ,
\]
where $\phi (t) =  1- e^{-t}$. 
The proof of this case consists of the several steps.

\begin{proposition} 
The solution $u= u (x,t)$ of the Cauchy problem (\ref{oneDphy01}) 
can be represented as follows
\begin{eqnarray}
\label{5.1}
&  &
u(x,t)  \\
& = &
\frac{1}{2} e ^{ \frac{t}{2}} \Big[     \varphi_0   (x+ 1- e^{-t})   
+   \varphi_0   (x -  1+ e^{-t} )\Big]  \nonumber \\
&  &
-  \int_{ 0}^{t}  \Big[ \frac{1}{4}e^{t/2}e^{-b/2} + \frac{1}{16}(1+4M^2)e^{t/2} e^{-3b/2}(e^t-e^b)\Big]  \Big[    \varphi_0   (x+ e^{-b}- e^{-t})   
+   \varphi_0   (x -  e^{-b}+ e^{-t}))\Big] \, db  \nonumber\\
&  &
+ \int_{ 0}^{t} db \int_{ x - (e^{-b}- e^{-t})}^{x+ e^{-b}- e^{-t}} dy \,  \varphi_{0}   (y) \left[ e^{-2b}  
\left( \frac{\partial }{\partial y} \right)^2 E(x-y,t;0,b) -
  M^2  E(x-y,t;0,b)  \right] \,.\nonumber
\end{eqnarray}
\end{proposition}
\medskip

\noindent
{\bf Proof.} We set $u(x,t)= w(x,t)+ \varphi_0 (x)$, then 
\[
w_{tt} - e^{-2t}w_{xx} +M^2w =e^{-2t}\varphi_{0,xx} -M^2 \varphi_0 (x)\, ,\qquad w(x,0) = 0 \,, \qquad w_t(x,0) =0\,  .
\]
Next we plug $f(x,t)= e^{-2t}\varphi_{0,xx}(x) -M^2 \varphi_0 (x) $ in the formula given by Theorem~\ref{T1.1} 
 and obtain
\begin{eqnarray}
\label{5.2a}
w(x,t)  
& = &
  \widetilde{ w(x,t)}-  \int_{ 0}^{t} db \int_{ x -  e^{-b}+ e^{-t} }^{x+ e^{-b}- e^{-t}} dy \,  M^2 \varphi_0 (y)   
E(x-y,t;0,b)  ,
\end{eqnarray}
where we have denoted
\begin{eqnarray*}
\widetilde{ w(x,t)} 
& := &
 \int_{ 0}^{t} db \int_{ x -  e^{-b}+ e^{-t} }^{x+ e^{-b}- e^{-t}} dy \,  e^{-2b}\varphi_{0}^{(2)} (y)   
E(x-y,t;0,b)   \,.
\end{eqnarray*}
Next we integrate by parts and apply (\ref{E_1}):
\begin{eqnarray*}
\widetilde{ w(x,t)} 
& = &
 \int_{ 0}^{t}   e^{-2b} \Big[ \varphi_{0}^{(1)} (x+ e^{-b}- e^{-t})   
- \varphi_{0}^{(1)} (x -  e^{-b}+ e^{-t} ) \Big]  
E(e^{-b}- e^{-t},t;0,b)  \, db  \\
&  &
- \int_{ 0}^{t} db \int_{ x -  e^{-b}+ e^{-t} }^{x+ e^{-b}- e^{-t}} dy \,  e^{-2b}\varphi_{0}^{(1)}  (y)   
\frac{\partial }{\partial y}E(x-y,t;0,b)  \,.
\end{eqnarray*}
On the other hand,
\[
 \varphi_0 ^{(1)} (x+ e^{-b}- e^{-t})  = -  e^{b} \frac{\partial }{\partial b}  \varphi_0   (x+ e^{-b}- e^{-t})  \quad
 \mbox{\rm and} \quad 
 \varphi_0 ^{(1)} (x -  e^{-b}+ e^{-t} )=    e^{b} \frac{\partial }{\partial b}  \varphi_0   (x -  e^{-b}+ e^{-t} )
\]
imply
\begin{eqnarray*}
\widetilde{ w(x,t)} 
& = &
- \int_{ 0}^{t}   e^{-b} \Big[    \frac{\partial }{\partial b}  \varphi_0   (x+ e^{-b}- e^{-t})   
+   \frac{\partial }{\partial b}  \varphi_0   (x -  e^{-b}+ e^{-t} )\Big]  
E(e^{-b}- e^{-t},t;0,b)  \, db  \\
&  &
- \int_{ 0}^{t} db \int_{ x -  e^{-b}+ e^{-t} }^{x+ e^{-b}- e^{-t}} dy \,  e^{-2b}\varphi_{0}^{(1)}  (y)   
\frac{\partial }{\partial y}E(x-y,t;0,b)  \,.
\end{eqnarray*}
One more integration by parts leads to
\begin{eqnarray*}
\widetilde{ w(x,t)}  + \varphi_0   (x)  
& = &
\frac{1}{2} e ^{ \frac{t}{2}} \Big[     \varphi_0   (x+ 1- e^{-t})   
+   \varphi_0   (x -  1+ e^{-t} )\Big]   \\
&  &
+ \int_{ 0}^{t}  \Big[    \varphi_0   (x+ e^{-b}- e^{-t})   
+   \varphi_0   (x -  e^{-b}+ e^{-t} )\Big]  
 \frac{\partial }{\partial b}\Big(  e^{-b} E(e^{-b}- e^{-t},t;0,b) \Big) \, db  \\
&  &
- \int_{ 0}^{t} db \int_{ x -  e^{-b}+ e^{-t} }^{x+ e^{-b}- e^{-t}} dy \,  e^{-2b}\varphi_{0}^{(1)}  (y)   
\frac{\partial }{\partial y}E(x-y,t;0,b) \,,
\end{eqnarray*}
where  $E( 0 ,t;0,t)=  \frac{1}{2} e^{t}$
  and $ E(1- e^{-t},t;0,0) =  \frac{1}{2} e ^{ \frac{t}{2}} $
have been used.
Next we apply  (\ref{E_3a})  of Proposition~\ref{P_E} 
and an integration by parts to obtain
\begin{eqnarray*}
\widetilde{ w(x,t)}  + \varphi_0   (x) 
& = &
\frac{1}{2} e ^{ \frac{t}{2}} \Big[     \varphi_0   (x+ 1- e^{-t})   
+   \varphi_0   (x -  1+ e^{-t} )\Big]   \\
&  &
-  \int_{ 0}^{t} \frac{1}{4}e^{t/2}e^{-b/2} \Big[    \varphi_0   (x+ e^{-b}- e^{-t})   
+   \varphi_0   (x -  e^{-b} + e^{-t} )\Big] \, db  \\
&  &
+ \int_{ 0}^{t} db \,  e^{-2b}
\Big[ \varphi_{0}   (y)   \frac{\partial }{\partial x}E(x-y,t;0,b) \Big]_{y=x -  e^{-b} + e^{-t} }^{y=x+ e^{-b}- e^{-t}}\\
&  &
+ \int_{ 0}^{t} db \int_{ x -  e^{-b} + e^{-t} }^{x+ e^{-b}- e^{-t}} dy \,  \varphi_{0}   (y) e^{-2b}  
\left( \frac{\partial }{\partial y} \right)^2 E(x-y,t;0,b) \,.
\end{eqnarray*}
We have due to (\ref{E_5}) and (\ref{E_6}) of Proposition \ref{P_E} 
\begin{eqnarray*}
&  &
\widetilde{ w(x,t)}  + \varphi_0   (x) \\
& = &
\frac{1}{2} e ^{ \frac{t}{2}} \Big[     \varphi_0   (x+ 1- e^{-t})   
+   \varphi_0   (x -  1+ e^{-t} )\Big]   \\
&  &
-  \int_{ 0}^{t}  \Big[ \frac{1}{4}e^{t/2}e^{-b/2} + \frac{1}{16}(1+4M^2)e^{t/2} e^{-3b/2}(e^t-e^b)\Big]  \Big[    \varphi_0   (x+ e^{-b}- e^{-t})   
+   \varphi_0   (x -  e^{-b}+ e^{-t}))\Big] \, db  \\
&  &
+ \int_{ 0}^{t} db \int_{ x -  e^{-b}+ e^{-t} }^{x+ e^{-b}- e^{-t}} dy \,  \varphi_{0}   (y) e^{-2b}  
\left( \frac{\partial }{\partial y} \right)^2 E(x-y,t;0,b) \,.
\end{eqnarray*}
Then the last equation together with (\ref{5.2a})  proves the desired   representation. Proposition is proven. \hfill $\square$
\smallskip

\noindent
{\bf Completion of the proof of Theorem~\ref{T1.3}.}  We make change $z= e^{-b }- e^{-t }$,\, $dz= -e^{-b }db$, and $b=-\ln (z+e^{t })$ in 
the second  term of the representation given by the previous proposition:
\begin{eqnarray*} 
&   &  
\int_{ 0}^{t}  \Big[ \frac{1}{4}e^{t/2}e^{-b/2} + \frac{1}{16}(1+4M^2)e^{t/2} e^{-3b/2}(e^t-e^b)\Big]  \Big[    \varphi_0   (x+ e^{-b}- e^{-t})   
+   \varphi_0   (x -  e^{-b}+ e^{-t}))\Big] \, db \\
& = &  
\int^{1-e^{-t }}_{0}  \Big[ \frac{1}{4}e^{t/2}\sqrt{z+ e^{-t}} + \frac{1}{16}(1+4M^2)e^{t/2} 
(\sqrt{z+ e^{-t}})^3\frac{ze^t}{z+ e^{-t}} \Big]  \Big[    \varphi_0   (x+ z)   
+   \varphi_0   (x -  z)\Big] \frac{1}{z+ e^{-t}}  \, dz  \\
& = &  
\int^{1-e^{-t }}_{0} e^{t } \Big[ \frac{1}{4} + \frac{1}{16}(1+4M^2) 
ze^t \Big]  \frac{1}{\sqrt{e^{ t}z+ 1}}  \Big[    \varphi_0   (x+ z)   
+   \varphi_0   (x -  z)\Big] \, dz  \,.
\end{eqnarray*}
Next we apply (\ref{E_1}) to the last term of that representation, and then we change the order of integration: 
\begin{eqnarray*}
&  &
\int_{ 0}^{t} db \int_{ x - (e^{-b}- e^{-t})}^{x+ e^{-b}- e^{-t}} dy \,  \varphi_{0}   (y) \left[ e^{-2b}  
\left( \frac{\partial }{\partial y} \right)^2 E(x-y,t;0,b) -
  M^2  E(x-y,t;0,b)  \right] \\
& = &
\int_{ 0}^{t} db \int_{ 0}^{ e^{-b}- e^{-t}} dz \, \Big[ \varphi_{0}   (x+z) + \varphi_{0}   (x-z)\Big]\left[ e^{-2b}  
\left( \frac{\partial }{\partial z} \right)^2 E(z,t;0,b) -
  M^2  E(z,t;0,b)  \right]  \\
& = &
\int_{ 0}^{1-e^{-t}} dz \Big[ \varphi_{0}   (x+z) + \varphi_{0}   (x-z)\Big] \int_{ 0}^{ -\ln(z+ e^{-t})} db \, \left[ e^{-2b}  
\left( \frac{\partial }{\partial z} \right)^2 E(z,t;0,b) -
  M^2  E(z,t;0,b)  \right]   \,.
\end{eqnarray*}
On the other hand, since the function $ E(z,t;0,b)$ solves Klein-Gordon equation,  
the last integral is equal to
\begin{eqnarray*}
 &  &
\int_{ 0}^{1-e^{-t}} dz \Big[ \varphi_{0}   (x+z) + \varphi_{0}   (x-z)\Big] \int_{ 0}^{ -\ln(z+ e^{-t})} db \, \left[ 
\left( \frac{\partial }{\partial t} \right)^2 E(z,t;0,b) \right] \\
 & = &
\int_{ 0}^{1-e^{-t}} dz \Big[ \varphi_{0}   (x+z) + \varphi_{0}   (x-z)\Big]  \left[ 
 \frac{\partial }{\partial b}   E(z,t;0,b) \right]_{b=0}^{b= -\ln(z+ e^{-t})} \,.
 \end{eqnarray*}
Application of (\ref{E_7}) and (\ref{E_8}) completes the proof.
Theorem~\ref{T1.3} is proven.
\hfill $\square$

\section{n-Dimensional Case, $n\geq 2$}
\label{S7}
\setcounter{equation}{0}

{\bf Proof of Theorem~\ref{T1.5}}.  
Let us consider the case of $x \in {\mathbb R}^n$, where $n=2m+1$,\, $m \in {\mathbb N}$.
First for the given function $u=u(x,t)$ we   define the spherical means of 
$u$ about point $x$:
\begin{eqnarray*}
I_u(x,r,t) 
 & = &
\frac{1}{\omega_{n-1} } \int_{S^{n-1}  } u(x+ry,t)\, dS_y \,,   
\end{eqnarray*}
where $\omega_{n-1} $ denotes the area of the unit sphere $S^{n-1} \subset {\mathbb R}^n$. Then we define
an operator $ \Omega _r $ by 
\[
\Omega _r ( u) (x,t) := \Big( \frac{1}{r} \frac{\partial }{\partial r}\Big)^{m-1} r^{2m-1}I_u(x,r,t) \,.
\]
One can show that there are constants $c_j^{(n)} $, $j=0,\ldots,m-1$, where $n=2m+1$, with
$c_0^{(n)} =1\cdot 3\cdot 5\cdots (n-2)$, \, 
such that
\[
\Big( \frac{1}{r} \frac{\partial }{\partial r}\Big)^{m-1} r^{2m-1} \varphi (r) 
= r \sum_{j=0}^{m-1} c_j^{(n)} r^{j} \frac{\partial^j }{\partial r^j} \varphi (r) \,.
\]
One can recover the functions according to
\begin{eqnarray}
\label{7.1}
u(x,t) 
& = &
\lim_{r \to 0}  I_u (x,r,t) = \lim_{r \to 0}  \frac{1}{c_0^{(n)}r} \Omega _r ( u) (x,t) \,,\\
\label{7.2}
u(x,0) 
& = &
\lim_{r \to 0}  \frac{1}{c_0^{(n)}r} \Omega _r ( u) (x,0) \,,
\quad u_t(x,0) = \lim_{r \to 0}  \frac{1}{c_0^{(n)}r} \Omega_r ( \partial_t u) (x,0)  \,. 
\end{eqnarray}
It is well known that $\Delta_x \Omega_r  h =  
\frac{\partial^2 }{\partial\, r^2}  \Omega_r  h $  for every function $h \in C^2({\mathbb R}^n)$.
Therefore we arrive at the following mixed problem for the  function $v(x,r,t) := \Omega_r  (u )(x,r,t) $: 
\begin{eqnarray*}
\cases{
v_{tt} (x,r,t) - e^{-2t} v_{rr}  (x,r,t) +M^2 v(x,r,t) = F (x,r,t)  
\quad  \mbox{\rm for all}\quad   t \ge 0\,, 
\,\, r \ge 0\,, \,\, x \in {\mathbb R}^n\,,\cr
v(x,0,t) = 0 
\quad \mbox{\rm for all}\quad t \ge 0\,, \quad x \in {\mathbb R}^n,\cr
v(x,r,0) = 0 \,, \quad v_t(x,r,0) = 0
\quad \mbox{\rm for all}
\quad r \ge 0\,, \quad x \in {\mathbb R}^n\,,\cr
F (x,r,t)  := \Omega_r  (f ) (x,t) \,, \quad  F (x,0,t)  =0 \,, 
\quad \mbox{\rm for all}\quad x \in {\mathbb R}^n\,.}
\end{eqnarray*}
It must be noted here that the spherical mean $I_u$ defined for $r>0$ has an 
extension as even function for $r<0 $ and hence $\Omega_r  (u ) $ has a natural extension as an
odd function. That allows replacing the mixed problem with the Cauchy problem. Namely, 
let functions $\widetilde v$ and $\widetilde F$ be the continuations of the functions $v$ and $F$, respectively, by 
\[
\hspace*{-0.7cm} 
\widetilde v(x,r,t)   = \cases{\,v (x,r,t) , \,\, if \,\,r \ge 0 \cr  - v (x,- r,t) , \,\,if \,\,r \le 0 }\,, \quad 
\widetilde F(x,r,t)   = \cases{\,F (x,r,t) , \,\,if \,\,r \ge 0 \cr  - F (x,- r,t) , \,\,if \,\,r \le 0 } \,.
\]
Then $\widetilde v$ solves the Cauchy problem
\begin{eqnarray*}
\cases{
\widetilde v_{tt} (x,r,t) - e^{-2t} \widetilde v_{rr}  (x,r,t) +M^2 \widetilde v  (x,r,t)= \widetilde F (x,r,t)  
\quad \mbox{\rm for all}\quad t \ge 0\,, 
\quad r \in {\mathbb R} \,, \quad x \in {\mathbb R}^n\,,\cr
\widetilde v(x,r,0) = 0 \,, \quad \widetilde v_t(x,r,0) = 0
\quad \mbox{\rm for all}
\quad r \in {\mathbb R} \,, \quad x \in {\mathbb R}^n .}
\end{eqnarray*} 
Hence, according to Theorem \ref{T1.1}, one has the representation
\begin{eqnarray*}
 \widetilde v (x,r,t)  
& = &
\int_{ 0}^{t} db \int_{ r - (e^{-b}- e^{-t} )}^{r+ e^{-b}- e^{-t} } dr_1 \,  \widetilde F (x,r_1,b)   (4e^{-b-t})^{iM}
\left( (e^{-t}  + e^{-b} )^2 - (r-r_1)^2   \right)^{-\frac{1}{2}-iM}\\
&  &
\hspace*{3cm}\times F\left(\frac{1}{2}+iM,\frac{1}{2}+iM;1; 
\frac{ (e^{-b}  - e^{-t} )^2 - (r-r_1)^2 }
{  (e^{-b}  + e^{-t} )^2 - (r-r_1)^2 } \right) . 
\end{eqnarray*}
Since $u(x,t) = \lim_{r \to 0} \big( \widetilde v (x,r,t)/(c_0^{(n)}r)\big) $,
we consider the case of $r < t$ in the above representation to obtain:
\begin{eqnarray*}
u(x,t) 
& = &
\frac{1}{c_0^{(n)} }   \int_{ 0}^{t} db
  \int_{ 0}^{ e^{-b}- e^{-t}} dr_1 \,  \lim_{r \to 0}
\frac{1}{ r}  \left\{\widetilde F (x,r+r_1,b) +   \widetilde F (x,r-r_1,b)\right\} \\
  &  & 
\times  (4e^{-b-t})^{iM}
\left( (e^{-t}  + e^{-b} )^2 - r_1^2   \right)^{-\frac{1}{2}-iM}
F\left(\frac{1}{2}+iM,\frac{1}{2}+iM;1; 
\frac{ (e^{-b}- e^{-t})^2-r_1^2}
{  (e^{-b}+ e^{-t})^2-r_1^2} \right) \,.
\end{eqnarray*}
Then by definition of the function $\widetilde F $, we replace 
$\lim_{r \to 0} \frac{1}{ r} \Big\{  \widetilde F (x,r- r_1,b)  + \widetilde F (x,r+ r_1,b) \Big\}  $
with \\
$2\Big( \frac{\partial }{\partial r} F (x,r ,b) \Big)_{r=r_1}  $ in the last formula. 
The definitions of $F (x,r ,t)$  and of the operator $\Omega _r $  yield:
\begin{eqnarray*}
u(x,t) 
& = &
\frac{2}{c_0^{(n)} }   \int_{ 0}^{t} db
  \int_{ 0}^{ e^{-b}- e^{-t}} dr_1 \,  \left( \frac{\partial }{\partial r} 
\Big( \frac{1}{r} \frac{\partial }{\partial r}\Big)^{m-1} r^{2m-1}I_f(x,r,t) 
\right)_{r=r_1} \\
  &  & 
\times  (4e^{-b-t})^{iM}
\left( (e^{-t}  + e^{-b} )^2 - r_1^2   \right)^{-\frac{1}{2}-iM}
F\left(\frac{1}{2}+iM,\frac{1}{2}+iM;1; 
\frac{ (e^{-b}- e^{-t})^2-r_1^2}
{  (e^{-b}+ e^{-t})^2-r_1^2} \right) \,,
\end{eqnarray*}
where $x \in {\mathbb R}^n$, $n=2m+1$, $m \in {\mathbb N}$. Thus, the solution to the Cauchy problem 
is given by (\ref{1.27}). We employ the method of descent to complete the proof for 
the case of even $n$, $n=2m$, $m \in {\mathbb N}$. 
Theorem~\ref{T1.5} is proven. 
\hfill $\square$

\noindent
{\bf Proof of (\ref{E+}) and (\ref{E-}).} We set $f(x,b)=\delta (x)\delta (t-t_0) $ in  (\ref{1.27}) and (\ref{1.28}),
and we obtain  (\ref{E+}) and (\ref{E-}),  
where if $n$ is odd, then  
\[
E^{w} (x,t) :=\frac{1}{\omega_{n-1} 1\cdot 3\cdot 5\ldots \cdot (n-2) } \frac{\partial }{\partial t} 
\Big( \frac{1}{t} \frac{\partial }{\partial t}\Big)^{\frac{n-3}{2} } 
\frac{1}{t} \delta (|x|-t) \, ,
\] 
while for $n$ even  we have
\[
 E^{w} (x,t) :=\frac{2}{\omega_{n-1} 1\cdot 3\cdot 5\ldots \cdot (n-1) } \frac{\partial }{\partial t} 
\Big( \frac{1}{t} \frac{\partial }{\partial t}\Big)^{\frac{n-2}{2} } 
\frac{1}{\sqrt{t^2-|x|^2}}\chi _{B_t(x)} \, .
\] 
Here $\chi _{B_t(x)} $ denotes  the characteristic function of the ball $B_t(x):= \{x \in {\mathbb R}^{n};\, |x| $ $\leq t \}  $.
Constant $\omega_{n-1} $ is the area of the unit sphere $S^{n-1} \subset {\mathbb R}^n$. The distribution $\delta (|x|-t) $ is defined  by 
$ <\delta (|\cdot |-t) , f (\cdot ) > = \int_{|x|=t}f(x)\, dx $ for   $f \in C^\infty _0({\mathbb R}^{n}) $.
\medskip

\noindent
{\bf Proof of Theorem~\ref{T1.6}}.  First we consider  case of  $\varphi_0  (x)=0$. More precisely, we have to prove that 
the solution $u (x,t)$ of the Cauchy problem (\ref{1.30CP}) with  $\varphi_0  (x)=0$
can be represented by (\ref{1.30}) with  $\varphi_0  (x)=0$. The next lemma will be used in both cases.

\begin{lemma}
\label{L8.2} 
Consider the mixed problem
\begin{eqnarray*}
\cases{ 
v_{tt} - e^{-2t}   v_{rr} +M^2 v=0 \quad \mbox{\rm for all}\quad t \ge 0\,, 
\quad r \ge 0\,,  \quad \cr  
v(r,0)= \tau _0 (r) \,, \quad v_t(r,0)= \tau _1 (r) \quad \mbox{\rm for all} 
\quad r \ge 0\,, \cr
v(0,t)=0 \quad \mbox{\rm for all}\quad t \ge 0\,,}
\end{eqnarray*}
and denote by  $\widetilde \tau _0 (r) $ and $\widetilde \tau _1 (r) $  the 
continuations of the functions $\tau _0 (r) $ and $\tau _1 (r) $  for negative $r$ as odd functions: 
\, $\widetilde \tau _0 (-r) =- \tau _0 (r)$ and  $\widetilde \tau _1 (-r) =- \tau _1 (r)$  for all $r \geq 0$, respectively. 
Then solution $v(r,t) $ to the mixed problem is given by the restriction of (\ref{3.1n}) to
$r \ge 0$:
\begin{eqnarray*} 
  v(r,t) 
& = &
\frac{1}{2} e ^{\frac{t}{2}}  \Big[ 
 \widetilde \tau _0 (r+ 1- e^{-t})  
+  \widetilde  \tau _0  (r - 1+ e^{-t})  \Big]   
+ \, \int_{ 0}^{1} \big[ 
 \widetilde \tau _0  (r - \phi (t)s)  
+   \widetilde \tau _0  (r  + \phi (t)s)  \big] K_0(\phi (t)s,t)\phi (t)\,  ds \\
&  &
+ \int_0^1\Big[  \widetilde \tau_1 \Big( r +  \phi (t) s \Big) +  \widetilde \tau_1 \Big( r -  \phi (t) s \Big) \Big] K_1(\phi (t)s,t) \phi (t)\,ds \,, 
\end{eqnarray*}
where $K_0(z,t)$ and  $K_1(z,t)$ are defined in Theorem~\ref{T1.3} and $\phi (t)= 1- e^{-t}$.
\end{lemma}
\medskip

\noindent
{\bf Proof.} This lemma is a direct consequence of Theorem~\ref{T1.3}. \hfill $\square$

\medskip

Now let us  consider the case of $x \in {\mathbb R}^n$, where $n=2m+1$.
First for the given function $u=u(x,t)$ we   define the spherical means of 
$u$ about  point $x$. 
One can recover the functions by means of (\ref{7.1}), (\ref{7.2}), and
\begin{eqnarray*}
\varphi_i(x) 
& = &
\lim_{r \to 0}  I_{\varphi_i} (x,r) = \lim_{r \to 0}  \frac{1}{c_0^{(n)}r} \Omega_r ( \varphi_i) (x) \,, \quad i=0,1   \,.
\end{eqnarray*}
Then 
we arrive at the following mixed problem 
\begin{eqnarray*}
\cases{
v_{tt} (x,r,t) - e^{-2t} v_{rr}  (x,r,t) +M^2v(x,r,t) = 0  
\quad  \mbox{\rm for all}\quad   t \ge 0\,, 
\,\, r \ge 0\,, \,\, x \in {\mathbb R}^n\,,\cr
v(x,0,t) = 0 
\quad \mbox{\rm for all}\quad t \ge 0\,, \quad x \in {\mathbb R}^n\,,\cr
v(x,r,0) = 0 \,, \quad v_t(x,r,0) = \Phi _1(x,r)
\quad \mbox{\rm for all}
\quad r \ge 0\,, \quad x \in {\mathbb R}^n\,,}
\end{eqnarray*}
with the  unknown function $v(x,r,t) := \Omega_r  (u )(x,r,t) $, where
\begin{eqnarray}
\label{8.2}
&  &
\Phi _i (x,r)  := \Omega_r  (\varphi_i  ) (x) 
= \Big( \frac{1}{r} \frac{\partial }{\partial r}\Big)^{m-1} r^{2m-1}\frac{1}{\omega_{n-1} } \int_{S^{n-1}  } 
\varphi_i(x+ry)\, dS_y
\,, \\
\label{8.3}
&  &
\Phi _i (x,0)  =0 \,, \quad i=0,1, \quad
\quad \mbox{\rm for all}\quad x \in {\mathbb R}^n\,.
\end{eqnarray} 
Then, according to Lemma \ref{L8.2} and to \,  
$u(x,t) = \lim_{r \to 0} \big(  v (x,r,t)/(c_0^{(n)}r)\big) $,
 we obtain:
\begin{eqnarray*}
u(x,t) 
 & = &
\frac{  1}{c_0^{(n)}} \lim_{r \to 0} \frac{  1}{r} \int_0^1\Big[ \widetilde \Phi _1 \big(x, r +  \phi (t) s \big) +   \widetilde \Phi _1 \big(x, r -  \phi (t) s \big) \Big] K_1(\phi (t)s,t) \phi (t)\,ds  \,.
\end{eqnarray*}
The last limit is equal to 
\begin{eqnarray*}
&   &
2 \int_0^1\left( \frac{ \partial }{\partial r} \Phi _1 (x,  r  )   
  \right)_{r= \phi (t)  s}  K_1(\phi (t)s,t) \phi (t)\,ds  \\
& = &
 2\int_0^1 \left( \frac{\partial}{\partial r} \Big( \frac{1}{r} \frac{\partial }{\partial r}\Big)^{\frac{n-3}{2} } 
\frac{r^{n-2}}{\omega_{n-1} } \int_{S^{n-1}  } 
\varphi_1 (x+ry)\, dS_y \right)_{r=\phi (t) s}    K_1(\phi (t)s,t) \phi (t)\,ds  \,.
\end{eqnarray*}
Thus, Theorem~\ref{T1.6} in the case of  $\varphi_0  (x)=0$ is proven.
\bigskip

Now we turn to  the case of  $\varphi_1  (x)=0$. 
Thus, we arrive at the following mixed problem
\begin{eqnarray*}
\cases{
v_{tt} (x,r,t) - e^{-2t} v_{rr}  (x,r,t) +M^2 v (x,r,t) = 0  
\quad  \mbox{\rm for all}\quad   t \ge 0\,, 
\,\, r \ge 0\,, \,\, x \in {\mathbb R}^n\,,\cr
v(x,r,0) = \Phi _0(x,r) \,, \quad v_t(x,r,0) = 0
\quad \mbox{\rm for all}
\quad r \ge 0\,, \quad x \in {\mathbb R}^n\,,\cr
v(x,0,t) = 0 
\quad \mbox{\rm for all}\quad t \ge 0\,, \quad x \in {\mathbb R}^n\,,}
\end{eqnarray*}
with the  unknown function $v(x,r,t) := \Omega_r  (u )(x,r,t) $ defined by (\ref{8.2}), (\ref{8.3}). 
Then, according to Lemma \ref{L8.2} and  to \,  
$u(x,t) = \lim_{r \to 0} \big(  v (x,r,t)/(c_0^{(n)}r)\big) $,
 we obtain:
\begin{eqnarray*}
u(x,t) 
 & = &
 \frac{1}{c_0^{(n)}} e ^{\frac{t}{2}} \lim_{r \to 0}\frac{1}{2r} \Big[ 
 \widetilde \Phi _0 (x,r+ e^t- 1)  
+  \widetilde  \Phi _0 (x,r - e^t+ 1)  \Big]   \\
&  &
+ \, \frac{2}{c_0^{(n)}} \int_{ 0}^{1} \lim_{r \to 0}\frac{1}{2r} \big[ 
 \widetilde \Phi _0 (x,r - \phi (t)s)  
+   \widetilde \Phi _0 (x,r  + \phi (t)s)  \big] K_0(\phi (t)s,t)\phi (t)\,  ds\,, \\
& = &
 \frac{1}{c_0^{(n)}} e ^{\frac{t}{2}} \left( \frac{\partial }{\partial r} 
\Phi _0 (x,r)  \right)_{r=\phi (t)}   
+ \, \frac{2}{c_0^{(n)}} \int_{ 0}^{1} \left( \frac{\partial }{\partial r} 
\Phi _0 (x,r)  \right)_{r=\phi (t)s}  K_0(\phi (t)s,t)\phi (t)\,  ds\\
& = &
 e ^{\frac{t}{2}} v_{\varphi_0}  (x, \phi (t))
+ \, 2\int_{ 0}^{1} v_{\varphi_0}  (x, \phi (t)s) K_0(\phi (t)s,t)\phi (t)\,  ds\,.
\end{eqnarray*}
Theorem~\ref{T1.6} is proven. \hfill $\square $

\section{$L^p-L^q$ and $L^q-L^q$ Estimates for the Solutions of  One-dimensional Equation,   $n=1$}
\label{S10}
\setcounter{equation}{0}

Consider now the Cauchy problem for the equation (\ref{Int_3})  with the source term and with vanishing initial  data (\ref{Int_4}).  
In this and next sections we restrict ourselves to the particles with ``large'' mass $m \geq n/2$,
that is, with nonnegative curved mass $M \geq 0$, 
to make presentation more transparent. The case of  $m < n/2$ will be discussed in the forthcoming paper.
\begin{theorem}
\label{T10.1}
For every function $f \in C^2 ({\mathbb R}\times [0,\infty))$ such that $f(\cdot ,t) \in C_0^\infty ({\mathbb R}_x)$ the solution 
$u = u(x,t)$ of the Cauchy problem (\ref{Int_3}), (\ref{Int_4})  satisfies inequality
\begin{eqnarray*}
\| u(x,t) \| _{L^q({\mathbb R}_x)} 
& \leq  &
C_M  e^{t(1-1/\rho )} \int_{ 0}^{t} (1+  t-b  )  (e^{t-b}   - 1)^{1/\rho } (e^{t-b}   + 1)^{-1 }   \|  f(x,b) \| _{L^p({\mathbb R}_x)} \, db  
\end{eqnarray*}
for all $t >0$, where $1 <p<\rho '$, $ \frac{1}{q} = \frac{1}{p} - \frac{1}{\rho '} $, $ \rho <2$, $ \frac{1}{\rho } + \frac{1}{\rho '}= 1$. 
\end{theorem}
\medskip

\noindent
{\bf Proof.} 
Using the fundamental solution from Theorem \ref{T1} one can write the convolution 
\[
\hspace*{-0.5cm} u(x,t)
 = 
\int_{ -\infty }^{\infty } \! \int_{ -\infty }^{\infty } {\mathcal E}_+ (x,t;y,b)f(y,b)\, db\, dy 
 = 
\int_{ 0}^{t}\! db \!\int_{ -\infty }^{\infty } {\mathcal E}_+ (x-y,t;0,b) f(y,b) \,dy  \,.
\]
Due to Young's inequality we have
\begin{eqnarray*}
\| u(x,t) \| _{L^q({\mathbb R}_x)} & \leq  &
c_k \int_{ 0}^{t} db \Bigg( \int_{  - (\phi (t) - \phi (b))}^{ \phi (t) - \phi (b)} 
|E (x,t;0,b)|^\rho  dx \Bigg)^{1/\rho } \|  f(x,b) \| _{L^p({\mathbb R}_x)} ,
\end{eqnarray*}
where $1 <p<\rho '$, $ \frac{1}{q} = \frac{1}{p} - \frac{1}{\rho '} $, $ \frac{1}{\rho } + \frac{1}{\rho '}= 1$, $\phi (t)=1- e^{-t}$. 
The integral in parentheses can be transformed as follows
\begin{eqnarray*}
&  &
 \int_{  - (\phi (t) - \phi (b))}^{ \phi (t) - \phi (b)} 
|E (x,t;0,b)|^\rho  dx \\
& = &
2 e^{-t +t\rho }\int_{ 0}^{ e^{t-b}  - 1} 
\Big((e^{t-b}  + 1)^2 - y^2\Big)^{-\frac{\rho }{2}   } 
\left| F\Big(\frac{1}{2}+iM   ,\frac{1}{2}+iM  ;1; 
\frac{ (e^{t-b}  - 1)^2 - y ^2 }{(e^{t-b}  + 1)^2 - y ^2 } \Big) \right|^\rho  dy  .
\end{eqnarray*}
On the other hand, due to integral representation for the hypergeometric function (1)\cite[v.1, Sec.2.12]{B-E} for $\zeta \in [0,1) $   one has
\begin{eqnarray}
\label{HypComplex}
\left| F\Big(\frac{1}{2}+iM   ,\frac{1}{2}+iM  ;1; 
\zeta \Big) \right|
&  \leq  &
 \frac{1}{|\Gamma \left(\frac{1}{2}+iM\right)|^2}  \pi F\Big(\frac{1}{2}  ,\frac{1}{2} ;1; 
\zeta \Big) \,.
\end{eqnarray} 
Thus,
\[
 \int_{  - (\phi (t) - \phi (b))}^{ \phi (t) - \phi (b)} 
|E (x,t;0,b)|^\rho  dx 
 \leq  
C_M  e^{-t +t\rho }\int_{ 0}^{ e^{t-b}  - 1} \!
\left((e^{t-b}  + 1)^2 - y^2\right)^{-\frac{\rho }{2}   } 
\left| F\Big(\frac{1}{2}  ,\frac{1}{2}  ;1; 
\frac{ (e^{t-b}  - 1)^2 - y ^2 }{(e^{t-b}  + 1)^2 - y ^2 } \Big) \right|^\rho \! dy .
\]
\begin{lemma}$\cite{Yag_Galst_Potsdam}$
\label{L10.2}
For all $z>1$ the following estimate 
\[
\int_{  0}^{ z  - 1} 
((z  + 1)^2  - r^2  )^{-\frac{\rho }{2}}  
F\left(\frac{1}{2},\frac{1}{2};1; 
\frac{ (z  - 1 )^2 -  r^2   }
{ (z  + 1 )^2  - r^2} \right)^\rho  d r  
  \leq  
C (1+ \ln   z )^\rho (z  - 1) (z  + 1)^{-\rho } F\Big(\frac{1}{2},\frac{\rho }{2};\frac{3}{2}; \frac{  (z  - 1)^2  }{ (z  + 1)^2  } \Big)
\]
is fulfilled, provided that $1 <p<\rho '$, $ \frac{1}{q} = \frac{1}{p} - \frac{1}{\rho '} $, $ \frac{1}{\rho } + \frac{1}{\rho '}= 1$. 
In particular, if $\rho < 2$, then 
\begin{eqnarray*}
\int_{  0}^{ z  - 1} 
((z  + 1)^2  - r^2  )^{-\frac{\rho }{2}}  
F\left(\frac{1}{2},\frac{1}{2};1; 
\frac{ (z  - 1 )^2 -  r^2   }
{ (z  + 1 )^2  - r^2} \right)^\rho  d r 
& \leq &
C (1+ \ln   z )^\rho (z  - 1) (z  + 1)^{-\rho } \,.
\end{eqnarray*}
\end{lemma}
\medskip

\medskip

\noindent
{\bf Completion of the proof of Theorem~\ref{T10.1}.} Thus for $\rho <2$ and $ z= e^{t-b}$ we have 
\begin{eqnarray*}
\| u(x,t) \| _{L^q({\mathbb R}_x)} 
& \leq  &
C_M  e^{t(1-1/\rho )} \int_{ 0}^{t} (1+  t-b  )  (e^{t-b}   - 1)^{1/\rho } (e^{t-b}   + 1)^{-1 }   \|  f(x,b) \| _{L^p({\mathbb R}_x)} \, db  \,.
\end{eqnarray*}
The last inequality implies the  estimate of the statement of theorem. Theorem~\ref{T10.1} is proven. \hfill $\square $
\bigskip

\begin{proposition}
\label{P9.3}
The solution $u=u(x,t) $ of the Cauchy problem 
\[
u_{tt} - e^{-2t}u_{xx} +M^2 u=0\, ,\qquad u(x,0) = \varphi_0  (x) \,, \qquad u_t(x,0) =\varphi_1  (x)\,  ,
\]
with $\varphi_0  $, $\varphi_1  \in C_0^\infty ({\mathbb R})$ satisfies the following estimate
\begin{equation}
\label{9.1}   
\| u(x,t) \| _{L^q({\mathbb R}_x)} 
\leq 
C (1+t)\Big(e^{\frac{t }{2}  } \| \varphi_0   (x)  \| _{L^q({\mathbb R}_x)}  
 +
 (e^{ t}  - 1) e^{ -t}\|  \varphi_1 (x) \| _{L^q({\mathbb R}_x)} \Big)
\, \quad \mbox{for all} \,\quad t \in (0,\infty) .
\end{equation}
\end{proposition}
\medskip  

\noindent
{\bf Proof.} First we consider the  equation without source term but with the second datum that is the case of $\varphi _0=0$. 
For the convenience we drop subindex of
$\varphi _1$.
Then we apply the representation given by Theorem~\ref{T1.3} for the solution $u=u (x,t)$ of the Cauchy problem with $\varphi _0=0$,
and obtain
\begin{eqnarray*}
u(x,t)   
&  =   &
\int_{0}^{1-  e^{-t}} \,\Big[      \varphi_1    (x- z)  +   \varphi _1   (x + z)    \Big] K_1(z,t) dz \,,
\end{eqnarray*}
where the kernel   $K_1(z,t)   $ is defined in Theorem~\ref{T1.3}.  Hence, we arrive at inequality 
\[  
\|  u(x,t) \|_{ { L}^{  q} ({\mathbb R})  } 
 \le    
C_M \|  \varphi   (x)   \|_{ { L}^{  q} ({\mathbb R})  }  \int_{0}^{1-  e^{-t}} \, |K_1(r,t)| dr 
\,.
\]
To estimate the  last integral we denote it by $I_1$, 
\begin{eqnarray*}  
I_1 (z)
& :=  &   \int_{0}^{1-  e^{-t}} \, |K_1(r,t)| dr \,,
\end{eqnarray*} 
and  with   $z=  e^{t}> 1$ due to (\ref{HypComplex}) we   write
\begin{eqnarray} 
\label{10.1}   
I_1 (z)
& \leq   &   
C \int_{0}^{z-1} \, \frac{1}{  \sqrt{(1+z)^2-y^2 } }    F\left(\frac{1}{2},\frac{1}{2};1; 
\frac{  y^2 - ( 1-z )^2  }
{   y^2 - ( 1+z )^2 } \right) dy 
\,.
\end{eqnarray} 
Then, according to Lemma \ref{L10.2} (the case of $\rho =1$) we have for that integral the following estimate
\begin{equation}
\label{10.4I1}  
I_1 (e^t)
  \leq  
C (1+ t )(e^t-1)(e^t+1)^{-1}\,.
\end{equation}
Finally, (\ref{10.1}) and  (\ref{10.4I1})  imply the $L^q- L^q$ estimate (\ref{9.1})  for the case of $\varphi _0=0$.
\medskip

Next we consider the equation without source but with the first datum, that is,   the  case of $\varphi _1=0$. 
We apply the representation given by Theorem~\ref{T1.3} for the solution $u=u (x,t)$ of the Cauchy problem with $\varphi _1=0$,
and obtain
\begin{eqnarray*}
u(x,t)   
&  =   &
\frac{1}{2} e ^{\frac{t}{2}}  \Big[ 
\varphi_0   (x+ 1-e^{-t})  
+     \varphi_0   (x -  1 +e^{-t})  \Big]  
+ \int_{ 0}^{1-  e^{-t}} \big[ 
\varphi_0   (x - r)  
+     \varphi_0   (x  + r)  \big] K_0(r,t)\,  dr \,,
\end{eqnarray*}
where the kernel   $K_0(r,t)   $ is defined in Theorem~\ref{T1.3}.
Then we easily obtain the following two estimates:
\begin{eqnarray*}
\|  u(x,t) - \int_{ 0}^{1-  e^{-t}} \big[ 
\varphi_0   (x - r)  
+     \varphi_0   (x  + r)  \big] K_0(r,t)\,  dr \|_{ { L}^{  q} ({\mathbb R})  } 
& \le &   
e ^{\frac{t}{2}} 
\|  \varphi_0   (x) \|_{ { L}^{  q} ({\mathbb R})  }   
\end{eqnarray*}
and
\begin{eqnarray*}
\|  u(x,t) \|_{ { L}^{  q} ({\mathbb R})  } 
& \le &   
e ^{\frac{t}{2}} 
\|  \varphi_0   (x) \|_{ { L}^{  q} ({\mathbb R})  }  
+ 2\|  \varphi_0   (x) \|_{ { L}^{  q} ({\mathbb R})  }  \int_{ 0}^{1-  e^{-t}} \left|   K_0(r,t)  \right|  dr \,.
\end{eqnarray*}
Finally,  the following lemma completes the proof of proposition.
\begin{lemma} 
\label{L10.3}
The kernel $K_0(r,t)$ has an integrable singularity at $r=e^t-1$, more precisely, one has 
\begin{eqnarray*}
\int_{ 0}^{1-e^{-t}} \left|   K_0(r,t)  \right|  dr \leq C(e^t-1)e^{-\frac{1}{2}t} (1+t) \quad \mbox{for all} \quad    t \in[0,\infty)\,. 
\end{eqnarray*} 
\end{lemma}
\medskip

\noindent
{\bf Proof.} 
For the  integral we obtain
\begin{eqnarray*}
 \int_{ 0}^{1-  e^{-t}} \left|   K_0(r,t)  \right|  dr 
& \leq  &
  \int_{ 0}^{z-  1}   \frac{1}{  [(z-1 )^2 -  y^2]\sqrt{ [(z+1)^2 - y^2]} }\\
&   &
\times  \Bigg|     \big(  z -z^{ 2 }  - iM(1 -     z^{ 2 } -  y^2) \big) 
F \Big(\frac{1}{2}+iM   ,\frac{1}{2}+iM  ;1; \frac{ (z -1)^2 -y^2 }{(z +1)^2 -y^2 }\Big) \\
&  &
\hspace{0.7cm}  +    \big( z^{2  }-1+  y^2 \big)\Big( \frac{1}{2}-iM\Big)
F \Big(-\frac{1}{2}+iM   ,\frac{1}{2}+iM  ;1;  \frac{ (z -1)^2 -y^2 }{(z +1)^2 -y^2 }\Big) \Bigg|  dy
\end{eqnarray*}
for all $z:=e^t >1$. We divide the domain of integration into  two zones, 
\begin{eqnarray}
\label{9.10}
Z_1(\varepsilon, z) 
& := &
\left\{ (z,r) \,\Big|\, \frac{ (z-1)^2 -r^2   }{ (z+1)^2 -r^2 } \leq \varepsilon,\,\, 0 \leq r \leq z-1 \right\} ,\\ 
\label{9.11}
Z_2(\varepsilon, z) 
& := &
\left\{ (z,r) \,\Big|\, \varepsilon \leq  \frac{ (z-1)^2 -r^2   }{ (z+1)^2 -r^2 },\,\, 0 \leq r \leq z-1  \right\},
\end{eqnarray}
and   split the integral into conformable two parts,
\begin{eqnarray*}
\int_{ 0}^{e^t-1} \left|   K_0(r,t)  \right|  dr 
& = &
\int_{ (z,r) \in Z_1(\varepsilon, z)   }   \left|   K_0(r,t)  \right|  dr 
+ \int_{ (z,r) \in Z_2(\varepsilon, z)  }  \left|   K_0(r,t)  \right|  dr \,.
\end{eqnarray*} 
In the first zone we have
\begin{eqnarray}
\label{9.8}
\hspace{-0.5cm} 
F\Big(\frac{1}{2}+iM,\frac{1}{2}+iM;1; \frac{ (z-1)^2 -y^2   }{ (z+1)^2 -y^2 }   \Big) 
\!\! & \!\! = \!\! &\!\! 
 1 + \left( \frac{1}{2}+iM \right)^2\frac{ (z-1)^2 -y^2   }{ (z+1)^2 -y^2 }   
+ O\left(\left( \frac{ (z-1)^2 -y^2   }{ (z+1)^2 -y^2 }\right)^2\right),  \\
\label{9.8a}
\hspace{-0.5cm} 
F\Big(-\frac{1}{2}+iM,\frac{1}{2}+iM;1; \frac{ (z-1)^2 -y^2   }{ (z+1)^2 -y^2 }   \Big) 
\!\! & \!\! = \!\! &\!\! 
 1 - \left( \frac{1}{4}+M^2 \right)\frac{ (z-1)^2 -y^2   }{ (z+1)^2 -y^2 }   
+ O\left(\left( \frac{ (z-1)^2 -y^2   }{ (z+1)^2 -y^2 }\right)^2\right)   . 
\end{eqnarray}
We use the last formulas to estimate the term containing the hypergeometric functions:
\begin{eqnarray}
\hspace{-0.5cm} &   &
\Bigg|  \big( z -z^{ 2 }  - iM(1 -    z^{ 2 }-  r^2)    \big) 
F \Big(\frac{1}{2}+iM   ,\frac{1}{2}+iM  ;1; \frac{ ( z-1)^2 -r^2 }{( z+1)^2 -r^2 }\Big) \nonumber \\
\hspace{-0.5cm} &  &
\hspace{1cm}   +    \big(z^{2 }- 1+  r^2  \big) \Big( \frac{1}{2}-iM\Big)
F \Big(-\frac{1}{2}+iM   ,\frac{1}{2}+iM  ;1; 
\frac{ ( z-1)^2 -r^2 }{( z+1)^2 -r^2 }\Big)  \Bigg|  \nonumber \\
\hspace{-0.5cm} & \leq  &
\frac{1}{2}\big[(z-1)^2 -r^2 \big] \nonumber \\
\hspace{-0.5cm} &  &
+\left| \big( z -z^{ 2 }  - iM(1 -    z^{ 2 }-  r^2)   \big) 
\left( \frac{1}{2}+iM \right)^2   -  \big(z^{2 }- 1+  r^2  \big) \Big( \frac{1}{2}-iM\Big)
  \left( \frac{1}{4}+M^2 \right) \right|  \frac{ (z-1)^2 -r^2   }{ (z+1)^2 -r^2 }  \nonumber \\
\hspace{-0.5cm} &  &
+ \Big( \left|  z -z^{ 2 }  - iM(1 -    z^{ 2 }-  r^2)      \right| 
+\left|  z^{2 }- 1+  y^2  \right|\Big)O\left(\left( \frac{ (z-1)^2 -r^2   }{ (z+1)^2 -r^2 }\right)^2\right)  \nonumber \\
\hspace{-0.5cm} & = &
\frac{1}{2}\big[(z-1)^2 -r^2 \big] \nonumber \\
\hspace{-0.5cm} &  &
+\frac{1}{8}\frac{ (z-1)^2 -r^2   }{ (z+1)^2 -r^2 }\left| (1-2 i M)  (-1+4 M^2 )  (y^2+z^2-1  )+2 (1+2 i M)^2  (-z^2+z +i M  (y^2+z^2 -1) )\right|  \nonumber  \\
\hspace{-0.5cm} &  &
\label{9.15}
 + \Big( \left|  z -z^{ 2 }  - iM(1 -    z^{ 2 }-  r^2)      \right| 
+\left|  z^{2 }- 1+  y^2  \right|\Big)O\left(\left( \frac{ (z-1)^2 -y^2   }{ (z+1)^2 -y^2 }\right)^2\right) .
\end{eqnarray}
Hence, we have to consider the following three integrals, which can be easily evaluated and estimated,
\begin{eqnarray*}
A_1
&  :=  &
\int_{ (z,r) \in Z_1(\varepsilon, z)  }  \frac{1}{     \sqrt{(z+1)^2-r^2}    }   dr   
  \leq   \mbox{\rm Arctan}\left( \frac{z-1}{2\sqrt{z}}\right) \leq 
\frac{\pi}{2} \,, \\
A_2
& := &
\int_{  (z,r) \in Z_1(\varepsilon, z) }  \frac{z^2}{ ((z+1)^2 -r^2 )\sqrt{(z+1)^2-r^2}    }   dr  
 \leq (z+1)^{-1/2}(z-1)   ,
\end{eqnarray*} 
and
\[
A_3
 := 
\int_{  (z,r) \in Z_1(\varepsilon, z) }  \frac{\left|  z -z^{ 2 }  - iM(1 -    z^{ 2 }-  r^2)     \right| 
+\left|  z^{2 }- 1+  r^2  \right| }{  \sqrt{(z+1)^2-r^2}    } \frac{ (z-1)^2 -r^2   }{ ((z+1)^2 -r^2)^2 }  dr  \leq  
C_M(z+1)^{-1/2}(z-1) 
\]
 for all  $z \in [1,\infty)$. Finally, for the integral over  the first zone we have obtained
\begin{eqnarray*}
&   &
  \int_{  (z,r) \in Z_1(\varepsilon, z) }   dr \frac{1}{  [(z-1 )^2 -  r^2]\sqrt{ [(z+1)^2 - r^2]} }\\
&   &
\times  \Bigg|     \big(  z -z^{ 2 }  - iM(1 -    z^{ 2 }-  r^2) \big) 
F \Big(\frac{1}{2}+iM   ,\frac{1}{2}+iM  ;1; \frac{ (z -1)^2 -r^2 }{(z +1)^2 -r^2 }\Big) \\
&  &
\hspace{1cm}  +    \big( z^{ 2 } -1+  r^2 \big)\Big( \frac{1}{2}-iM\Big)
F \Big(-\frac{1}{2}+iM   ,\frac{1}{2}+iM  ;1;  \frac{ (z -1)^2 -r^2 }{(z +1)^2 -r^2 }\Big) \Bigg| 
\, \leq  
\,C_M (z+1)^{-1/2}(z-1) 
\end{eqnarray*}
for all $z \in [1,\infty)$. In the second zone we have
\begin{equation}
\label{10.4n}
\varepsilon \leq  \frac{ (z-1)^2 -r^2   }{ (z+1)^2 -r^2 } \leq 1 \quad \mbox{\rm and}  \quad 
\frac{ 1  }{ (z-1)^2 -r^2 }  \leq  \frac{ 1   }{ \varepsilon[(z+1)^2 -r^2] }\,.
\end{equation}
According to the formula 15.3.10 of \cite[Ch.15]{B-E}   the hypergeometric functions obey the estimate
\begin{equation}
\label{FGH+xat1}
\hspace{-0.5cm} \left| F\Big(-\frac{1}{2}+iM,\frac{1}{2}+iM;1; x   \Big) \right|  \leq C \,\,  \mbox{\rm and}  \,\,   
\left| F\Big(\frac{1}{2}+iM,\frac{1}{2}+iM;1; x  \Big) \right|  \leq C \big(1-\ln(1-x)) \,\, \, \forall   x \in [\varepsilon ,1) .
\end{equation}
This allows to  estimate  the integral over the second zone:
\begin{eqnarray}
\label{10.4}
&   &
  \int_{  (z,r) \in Z_2(\varepsilon, z) }   dr \frac{1}{  [(z-1 )^2 -  r^2]\sqrt{ (z+1)^2 - r^2} } \\
&   &
\times  \Bigg|     \big(  z -z^{ 2 }  - iM(1 -    z^{ 2 }-  r^2) \big) 
F \Big(\frac{1}{2}+iM   ,\frac{1}{2}+iM  ;1; \frac{ (z -1)^2 -r^2 }{(z +1)^2 -r^2 }\Big) \nonumber \\
&  &
\hspace{0.5cm}  +    \big( z^{ 2 } -1+  r^2 \big)\Big( \frac{1}{2}-iM\Big)
F \Big(-\frac{1}{2}+iM   ,\frac{1}{2}+iM  ;1;  \frac{ (z -1)^2 -r^2 }{(z +1)^2 -r^2 }\Big) \Bigg| 
 \leq  
C_M(z+1)^{-1/2}(z-1)  \nonumber 
\end{eqnarray}
for all $z \in [1,\infty)$. Indeed, for the argument of the hypergeometric functions we have
\begin{equation}
\label{8.13}
\varepsilon \leq \frac{ (z-1)^2 -r^2   }{ (z+1)^2 -r^2 }  = 1- \frac{  4z  } { (z  + 1)^2 -r ^2} < 1,\quad 
\frac{  4z  } { (z  + 1)^2 -r ^2} < 1- \varepsilon \quad \mbox{\rm for all}  \quad (z,r) \in Z_2(\varepsilon, z) \,.
\end{equation}
Hence,
\begin{equation}
\label{10.8}
\left| F\Big(\frac{1}{2}+iM,\frac{1}{2}+iM;1; \frac{ (z-1)^2 -r^2   }{ (z+1)^2 -r^2 }    \Big) \right|   
  \leq   
C \left(1+ \ln  z     \right) 
, \, (z,r) \in Z_2(\varepsilon, z)  . 
\end{equation}
To prove (\ref{10.4}) we estimate the following  integral
\begin{eqnarray*} 
\int_{  (z,r) \in Z_2(\varepsilon, z) }   \frac{ z^{2} }{     ((z-1)^2-r^2 ) \sqrt{(z+1)^2-r^2}    } 
 dr 
& \leq    & 
C_\varepsilon  z^{2}  \int_{  0 }^{z-1}  \frac{1}{ ((z+1)^2-r^2 )^{3/2}   }  dr 
\, \leq   \, 
C_\varepsilon     \frac{ (z -1)}{\sqrt{z} } \,.
\end{eqnarray*} 
Thus, the  lemma is proven.    \hfill $\square$

\section{$L^p-L^q$  Estimates for  Equation with $n=1$ and without Source Term. Some Estimates of  Kernels $K_0$ and $K_1$ }
\label{S11}
\setcounter{equation}{0}

\begin{theorem}
Let $u=u(x,t) $  be a solution of the Cauchy problem 
\[
u_{tt} - e^{-2t}u_{xx} +M^2 u=0\, ,\qquad u(x,0) = \varphi_0  (x) \,, \qquad u_t(x,0) =\varphi_1  (x)\,  ,
\]
with $\varphi_0  $, $\varphi_1  \in C_0^\infty ({\mathbb R})$. 
If $\rho \in (1,2) $, then
\begin{eqnarray*}
\| u(x,t) \| _{L^q({\mathbb R}_x)} 
& \leq &
e ^{\frac{t}{2}} 
\| \varphi_0   (x)  \| _{L^q({\mathbb R}_x)}  
+  C_\rho (1+t)(e^{ t}  - 1)^{\frac{1}{\rho }}e^{t[\frac{1}{2}- \frac{1}{\rho }]} \| \varphi_0   (x)  \| _{L^p({\mathbb R}_x)}  \\
&  &
 +\,
C_\rho (1+  t)  (e^{ t}  - 1)^{\frac{1}{\rho }} e^{ -\frac{t}{\rho }}\|  \varphi_1 (x) \| _{L^p({\mathbb R}_x)} 
\,,
\end{eqnarray*}
for all $ t \in (0,\infty)$. Here $1 <p<\rho '$, $ \frac{1}{q} = \frac{1}{p} - \frac{1}{\rho '} $, $ \frac{1}{\rho } + \frac{1}{\rho '}= 1$. 
 If $\rho =1$, then 
\begin{equation}
\label{10.1LpLq}   
\| u(x,t) \| _{L^q({\mathbb R}_x)} 
\leq 
C (1+t)\Big(e^{\frac{t }{2}  } \| \varphi_0   (x)  \| _{L^q({\mathbb R}_x)}  
 +
 (e^{ t}  - 1) e^{ -t}\|  \varphi_1 (x) \| _{L^q({\mathbb R}_x)} \Big)
\,,
\end{equation}
for all $ t \in (0,\infty)$. 
\end{theorem}
\medskip

\noindent
{\bf Proof.}
For $\rho =1 $ we just apply Proposition~\ref{P9.3}. 
To prove this theorem for  $\rho >1 $  we need some auxiliary estimates for the kernels $K_0$ and $K_1$. We start with the case of $\varphi_0  =0 $,
 where the kernel $K_1$
appears.  The application of Theorem~\ref{T1.3} and Young's inequality lead  to 
\begin{eqnarray*}
\| u(x,t) \| _{L^q({\mathbb R}_x)} 
& \leq  &
2\Bigg( \int_{0}^{ 1 -e^{-t} }
| K_1(x,t)|^\rho  dx \Bigg)^{1/\rho } \|  \varphi_1 (x) \| _{L^p({\mathbb R}_x)} ,
\end{eqnarray*}
where $1 <p<\rho '$, $ \frac{1}{q} = \frac{1}{p} - \frac{1}{\rho '} $, $ \frac{1}{\rho } + \frac{1}{\rho '}= 1$.  Now we have to estimate
the last integral.
\begin{proposition}
We have
\begin{eqnarray*}
\left( \int_{0}^{ 1 -e^{-t}  } 
| K_1(x,t)|^\rho  dx \right)^{1/\rho } 
& \leq  &
 C  (1+ t )  (1- e^{-t} )^{1/\rho }  \quad \mbox{for all} \,\quad t \in (0,\infty) \,.
\end{eqnarray*}
\end{proposition}
\medskip

\noindent
{\bf Proof.} One can write
\[
\left( \int_{0}^{1 -e^{-t}  } 
| K_1(x,t)|^\rho  dx \right)^{1/\rho }  \leq   
C_Me^{ t(1-1/\rho )}  \Bigg( \int_{0}^{e^t- 1} 
\big((e^t+1 )^2 -  y^2\big)^{-\frac{\rho }{2}   }  \Big| F\left(\frac{1}{2}   ,\frac{1}{2}  ;1; 
\frac{ (e^t- 1 )^2 -y^2 }{( e^t+ 1 )^2 -y^2 } \right) \Big|^\rho  dy \Bigg)^{1/\rho }.
\]
Denote $z:= e^t >1 $ and consider the  integral $
\displaystyle \int_{0}^{ z-1 } 
\left| \frac{1}{  \sqrt{(1+z)^2-x^2 } }    F\Big(\frac{1}{2},\frac{1}{2};1; 
\frac{  (z  - 1  )^2 - x^2}
{  (z  + 1)^2-x^2} \Big)  \right|^\rho  dx  $ of the right-hand side. 
Then we apply Lemma \ref{L10.2} and obtain 
\[
\left( \int_{0}^{1 -e^{-t}  } 
| K_1(x,t)|^\rho  dx \right)^{1/\rho } 
  \leq  
 Ce^{ t(1-1/\rho )}  (1+ \ln   e^t )  (e^t  - 1)^{1/\rho } (e^t  + 1)^{-1 } \leq 
 C  (1+ t )  (1- e^{-t} )^{1/\rho } \,.
\]
Proposition is proven. \hfill $\square$
\medskip

Thus, the theorem in the case of $\varphi_0  =0 $ is proven.
\medskip

\noindent
Now we turn to the case of $\varphi_1  =0 $,
 where the kernel $K_0$
appears.  The application of Theorem~\ref{T1.3} leads to 
\begin{eqnarray*}
\| u(x,t) \| _{L^q({\mathbb R}_x)} 
& \leq  &  
e ^{\frac{t}{2}} 
\| \varphi_0   (x)  \| _{L^q({\mathbb R}_x)}  
+  \, \left\| \int_{ 0}^{1 -e^{-t} }\left[ 
\varphi_0   (x - z) 
+     \varphi_0   (x  + z)  \right] K_0(z,t)\, dz  \right\| _{L^q({\mathbb R}_x)} \,.
\end{eqnarray*} 
Similarly to the case of the second datum we arrive at
\begin{eqnarray*}
\| u(x,t) \| _{L^q({\mathbb R}_x)} 
& \leq  &  
e ^{\frac{t}{2}} 
\| \varphi_0   (x)  \| _{L^q({\mathbb R}_x)}  
+  \, 
\| \varphi_0   (x )  \| _{L^p({\mathbb R}_x)} \left( \int_{0}^{ 1 -e^{-t} } | K_0(r,t)|^\rho  dr \right)^{1/\rho }\,.
\end{eqnarray*} 
The next proposition gives an estimate for  the integral of the last inequality. 
\begin{proposition}
\label{P11.4}
Let $1 <p<\rho '$, $ \frac{1}{q} = \frac{1}{p} - \frac{1}{\rho '} $, $ \frac{1}{\rho } + \frac{1}{\rho '}= 1$, and $ \rho \in [1,2)$.  We have
\begin{eqnarray*}
\left( \int_{0}^{  1 -e^{-t} } 
| K_0(r,t)|^\rho  dr \right)^{1/\rho } 
& \leq  &
C_\rho  (1+t)(e^{t}-1)^{  \frac{1}{\rho}  }e^{t(\frac{1}{2} -\frac{1}{\rho } )} \quad \mbox{ for all} \quad t \in (0,\infty)  \,.
\end{eqnarray*}
\end{proposition}
\medskip

\noindent
{\bf Proof.} 
We turn to the integral ($z:= e^t>1$)
\begin{eqnarray*}
\left( \int_{0}^{  1 -e^{-t} } 
| K_0(r,t)|^\rho  dr \right)^{1/\rho }  
\!\!& \!\!= \!\!&\!\!
\Bigg( \int_{0}^{z- 1 }  dy 
\left( \frac{1}{ [(z-1  )^2 -  y^2]\sqrt{(z+1)^2 - y^2} } \right)^\rho\\
&   &
\times \Bigg|  \big(  z -z^{2 } - iM( 1-    z^{2 }   -  y^2) \big) 
F \Big(\frac{1}{2}+iM   ,\frac{1}{2}+iM  ;1; \frac{ (z -1)^2 -y^2  }{(z +1)^2 -y^2  }\Big) \\
&  &
 +   \big(z^{2 } -1+  y^2 \big)\Big( \frac{1}{2}-iM\Big)
F \Big(-\frac{1}{2}+iM   ,\frac{1}{2}+iM  ;1; \frac{(z -1)^2 -y^2  }{(z +1)^2 -y^2 }\Big)  \Bigg|^\rho \Bigg)^{1/\rho }.
\end{eqnarray*}
The formulas (\ref{9.8}) and (\ref{9.8a})  describe the behavior of the  hypergeometric functions in the neighbourhood of zero.
Consider therefore two zones, $Z_1(\varepsilon, z) $ and $ Z_2(\varepsilon, z )$, defined in  (\ref{9.10}) and (\ref{9.11}),
respectively.
We  split integral into two parts:
\begin{eqnarray*}
\int_{ 0}^{1-e^{-t} } \left|   K_0(r,t)  \right|^\rho   dr 
& = &
\int_{ (z,r) \in Z_1(\varepsilon, z)   }   \left|   K_0(r,t)  \right|^\rho   dr 
+ \int_{ (z,r) \in Z_2(\varepsilon, z)  }  \left|   K_0(r,t)  \right|^\rho   dr \,.
\end{eqnarray*} 
In the proof of Lemma~\ref{L10.3} the   relation (\ref{9.15}) was checked in the first zone.  
If $1\leq  z \leq N$ with some constant $N$, then the argument of the hypergeometric functions is bounded,
\begin{equation}
\label{argument}
\frac{ (z-1)^2 -r^2   }{ (z+1)^2 -r^2 } \leq \frac{ (z-1)^2  }{ (z+1)^2 } \leq \frac{ (N-1)^2  }{ (N+1)^2 }  <1 \quad \mbox{for all} \quad r \in (0,z-1),
\end{equation}
and we obtain   with $z= e^{t }$, 
\begin{eqnarray*}
\Bigg( \int_{ 0}^{1-e^{-t} } \left|   K_0(r,t)  \right|^\rho   dr  \Bigg)^{1/\rho}
& \leq  & 
C_{M,N} \Bigg( \int_{0}^{z- 1 } 
\Bigg[ \frac{1}{ [(z-1  )^2 -  y^2]\sqrt{(z+1)^2 - y^2} } \\
&   &
\hspace{0.1cm} \times   
\Bigg\{ \frac{1}{2}[(z -1)^2 -y^2]+ z^2\frac{(z -1)^2 -y^2  }{(z +1)^2 -y^2 }+ z^2\left(\frac{(z -1)^2 -y^2  }{(z +1)^2 -y^2 } \right)^2 \Bigg\} 
 \Bigg]^\rho  dy \Bigg)^{1/\rho }\\
& \leq  & 
C_{M,N} \Bigg( \int_{0}^{z- 1 } 
\Bigg[ \frac{1}{  \sqrt{(z+1)^2 - y^2} }   
\Bigg\{ 1+ z^2\frac{1 }{(z +1)^2 -y^2 }   \Bigg\}  
 \Bigg]^\rho  dy \Bigg)^{1/\rho }\\
& \leq  & 
C_{M,N}   (z-1)^{1/\rho } (z+1)^{-1 }   \,.
\end{eqnarray*}
Thus, we can restrict ourselves to the  case of large $   z \geq N$ in both zones.
Consider therefore for $\rho \in (1,2)$ the following integrals over the first zone
\begin{eqnarray*}
A_4
& := &
\int_{(z,r)\in Z_1(\varepsilon )} \left( \frac{1}{   \sqrt{(z+1)^2-r^2}    } \right)^\rho  dr  
 \,\, \leq    \,\,
\int_{0}^{z-1}   \left(\frac{1}{   \sqrt{(z+1)^2-r^2}    } \right)^\rho dr \\
& \leq  & 
C  (z-1)(z+1)^{-\rho }F\Big(\frac{1}{2},\frac{\rho }{2};\frac{3}{2}; \frac{ (z-1)^2     }{ (z+1)^2   }   \Big)  \\
& \leq  &
C  (z-1)(z+1)^{-\rho }\,, \\
A_5
& := &
\int_{(z,r)\in Z_1(\varepsilon )} \left( \frac{z^2}{[(z+1)^2-r^2  ]   \sqrt{(z+1)^2-r^2}    } \right)^\rho  dr    
 \leq   
C z^{2\rho } (z-1)(z+1)^{-3\rho }F\Big(\frac{1}{2},\frac{3\rho }{2};\frac{3}{2}; \frac{ (z-1)^2     }{ (z+1)^2   }   \Big) .
\end{eqnarray*}
Then, according to  (15.3.6) of Ch.15\cite{A-S}  and \cite{B-E},
\begin{eqnarray}
\label{A-S_15.3.6}
 F  \left( a,b;c;\zeta   \right) 
& = &
\frac{\Gamma (c)\Gamma (c-a-b)}{\Gamma (c-a)\Gamma (c-b)}F  \left( a,b;a+b-c+1;1-\zeta   \right) \\
 &  &
 + (1-z)^{c-a-b}\frac{\Gamma (c)\Gamma (a+b-c)}{\Gamma (a)\Gamma (b)}F  \left( c-a,c-b;c-a-b+1;1-\zeta   \right)  \nonumber
\end{eqnarray}
for all $\zeta \in {\mathbb C} $, \,$ |\arg (1-\zeta )| <\pi $. We use (\ref{A-S_15.3.6}) with  
\[
\zeta = \frac{ (z-1)^2     }{ (z+1)^2   } ,\qquad 1- \zeta =  \frac{ 4z     }{ (z+1)^2   }
\]
to obtain   for $\rho <2$ and large $z \geq N$ the following estimate for the hypergeometric function,
\begin{equation}
\label{8.4}
F\Big(\frac{1}{2},\frac{3\rho }{2};\frac{3}{2}; \frac{ (z-1)^2     }{ (z+1)^2   }   \Big) 
 \leq 
C (z+1)^{-1+\frac{3\rho }{2} }\,.
\end{equation}
Thus,
\begin{eqnarray*}
A_5
& \leq  & C  (z-1) (z+1)^{-1+\frac{\rho }{2} }\,.
\end{eqnarray*}
For the next term we obtain a similar estimate,
\begin{eqnarray*}
A_6
 := 
\int_{(z,r)\in Z_1(\varepsilon )} \left| \frac{z^{2}}{    ((z-1)^2-r^2 ) \sqrt{(z+1)^2-r^2}    } 
\left( \frac{ (z-1)^2 -r^2   }{ (z+1)^2 -r^2 }\right)^2\right|^\rho  dr 
& \leq  & 
C  (z-1) (z+1)^{-1+\frac{\rho }{2} }\,.
\end{eqnarray*} 
Hence,
\begin{eqnarray*}
&  &
\int_{ (z,r) \in Z_1(\varepsilon, z)   }   \left|   K_0(r,t)  \right|^\rho   dr 
\leq  C (z-1) (z+1)^{-1+\frac{\rho }{2} }\,.
\end{eqnarray*} 
In the second zone $Z_2(\varepsilon, z)$  for the argument of the hypergeometric functions we have (\ref{10.4n}),
(\ref{8.13}), and (\ref{10.8}). 
We have to estimate the following  integral
\begin{eqnarray*} 
A_7
& :=   &
\int_{  (z,r) \in Z_2(\varepsilon, z) }  \left| \frac{z^{2} \left(1+ \ln  z     \right) }{     ((z-1)^2-r^2 ) \sqrt{(z+1)^2-r^2}    }   
 \right|^\rho   dr  \,.
\end{eqnarray*} 
We apply (\ref{10.4n}) and (\ref{8.4}) to obtain 
\begin{eqnarray*} 
A_7
& \leq    &
C\left(1+ \ln  z     \right)^{\rho } (z-1)(z+1)^{-1+\frac{\rho }{2} } \,.
\end{eqnarray*} 
Hence,
\begin{eqnarray*}
\int_{ (z,r) \in Z_2(\varepsilon, z)  }  \left|   K_0(r,t)  \right|^\rho   dr \leq C\left(1+ \ln  z     \right)^{\rho } (z-1)(z+1)^{-1+\frac{\rho }{2} } 
\qquad \mbox{\rm for all}  \quad z\geq N \,.
\end{eqnarray*} 
Proposition is proven. \hfill $\square $

\section{$L^p-L^q$ Estimates for the Equation with  Source, $n\geq 2$}
\label{S12}
\setcounter{equation}{0}

For the wave equation the Duhamel's principle allows to reduce the case of source term to the case of the Cauchy problem without
source term and consequently to derive the $L^p-L^q$-decay estimates for the equation. For (\ref{K_G_linear}) the Duhamel's principle
is not applicable straightforward and we have to appeal to the representation formula of Theorem~\ref{T1.5}. In fact, one can regard that
formula as an expansion of the two-stage Duhamel's principle.
In this section we consider the Cauchy problem (\ref{1.26}) for the equation with the source term with zero initial data.   
\begin{theorem} 
\label{T11.1}
Let $u=u(x,t)$ be solution of the Cauchy problem (\ref{1.26}). Then for $n>1$ one has the following decay estimate
\begin{eqnarray*}
&  &
\| (-\bigtriangleup )^{-s} u(x,t) \|_{ { L}^{  q} ({\mathbb R}^n)  } \\
& \le &  
C \int_{ 0}^{t} db\, \|  f(x, b)  \|_{ { L}^{p} ({\mathbb R}^n)  }   
  \int_{ 0}^{ e^{-b}- e^{-t}}  dr \, r^{ 2s-n(\frac{1}{p}-\frac{1}{q}) }  \frac{1}{\sqrt{(e^{-t}  + e^{-b})^2    - r^2}}
F\left(\frac{1}{2},\frac{1}{2};1; 
\frac{ (e^{-b}- e^{-t})^2-r^2}
{  (e^{-b}+ e^{-t})^2-r^2}   \right) 
\end{eqnarray*}
provided that $s \ge 0$, $1<p \leq 2$, $\frac{1}{p}+ \frac{1}{q}=1$, $\frac{1}{2} (n+1)\left( \frac{1}{p} - \frac{1}{q}\right) \leq 
2s \leq n \left( \frac{1}{p} - \frac{1}{q}\right)< 2s+1 $.
\end{theorem}
\medskip

\noindent
{\bf Proof.} In both cases, of even and odd $n$, one can write the representation  (\ref{1.29}). 
Due to the results of \cite{Brenner,Pecher} for the wave equation, we have
\begin{eqnarray*}
&  &
\| (-\bigtriangleup )^{-s} u(x,t) \|_{ { L}^{  q} ({\mathbb R}^n)  } \\
& \le &  
C\int_{ 0}^{t} db
  \int_{ 0}^{ e^{-b}- e^{-t}}  \| (-\bigtriangleup )^{-s} v(x,r ;b)  \|_{ { L}^{  q} ({\mathbb R}^n)  }   \frac{1}{\sqrt{(e^{-t}  + e^{-b})^2    - r^2}}
F\left(\frac{1}{2},\frac{1}{2};1; 
\frac{ (e^{-b}- e^{-t})^2-r^2}
{  (e^{-b}+ e^{-t})^2-r^2}  \right) dr \\ 
& \le &  
C \int_{ 0}^{t} db\, \|  f(x, b)  \|_{ { L}^{p} ({\mathbb R}^n)  }   
  \int_{ 0}^{ e^{-b}- e^{-t}} r^{ 2s-n(\frac{1}{p}-\frac{1}{q}) }  \frac{1}{\sqrt{(e^{-t}  + e^{-b})^2    - r^2}}
F\left(\frac{1}{2},\frac{1}{2};1; 
\frac{ (e^{-b}- e^{-t})^2-r^2}
{  (e^{-b}+ e^{-t})^2-r^2}   \right) dr\,.
\end{eqnarray*}
The theorem is proven. \hfill $\square$

We are going to transform the estimate of the last theorem to more cosy form. 
To this aim we  estimate  for $n(\frac{1}{p}-\frac{1}{q})<2s+1 $ the  last integral of the right hand side. 
If we replace $e^{-b}/ e^{-t} >  1$ with $z := e^{-b}/ e^{-t}>  1$, then the integral will be simplified.
\begin{eqnarray*}
&  &
\int_{ 0}^{ e^{-b}- e^{-t}} r^{ 2s-n(\frac{1}{p}-\frac{1}{q}) }  \frac{1}{\sqrt{(e^{-t}  + e^{-b})^2    - r^2}}
F\left(\frac{1}{2},\frac{1}{2};1; 
\frac{ (e^{-b}- e^{-t})^2-r^2}
{  (e^{-b}+ e^{-t})^2-r^2}   \right) dr \\
& = &
e^{-t[2s-n(\frac{1}{p}-\frac{1}{q})]}\int_{ 0}^{ z- 1} y^{ 2s-n(\frac{1}{p}-\frac{1}{q}) }  \frac{1}{\sqrt{(z  + 1)^2    - y^2}}
F\left(\frac{1}{2},\frac{1}{2};1; 
\frac{ ( z  - 1 )^2-y^2}
{   (z  + 1) ^2-y^2}   \right) dy
\end{eqnarray*}
\begin{lemma}
\label{L12.2}$\cite[Lemma~9.2]{Yag_Galst_Potsdam}$
Assume that $0 \ge 2s-n(\frac{1}{p}-\frac{1}{q})>-1 $. Then 
\[
  \int_{ 0}^{ z- 1}  r^{ 2s-n(\frac{1}{p}-\frac{1}{q}) }  \frac{1}{\sqrt{(z  + 1)^2    - r^2}}
F\left(\frac{1}{2},\frac{1}{2};1; 
\frac{  (z  - 1)^2  -r^2 }
{  (z  + 1)^2  -r^2  }  \right)  \, dr  
 \leq   
C     z^{-1}(z-1)^{ 1+  2s-n(\frac{1}{p}-\frac{1}{q}) } (1+ \ln z ),
\]
for all $z>1$.
\end{lemma}
\medskip

\begin{corollary}
\label{C12.3}
Let $u=u(x,t)$ be solution of the Cauchy problem (\ref{1.26}). Then for $n\geq 2$ one has the following decay estimate  
 \begin{eqnarray}
 \label{12.1}
\| (-\bigtriangleup )^{-s} u(x,t) \|_{ { L}^{  q} ({\mathbb R}^n)  } 
& \le &  
C\int_{ 0}^{t} \|  f(x, b)  \|_{ { L}^{p} ({\mathbb R}^n)  }   
 e^{-b}\left( e^{-b}- e^{-t}  \right)^{ 1+  2s-n(\frac{1}{p}-\frac{1}{q}) } \left(1+  t-  b    \right)\,  db 
\end{eqnarray}
provided that $s \ge 0$, $1<p \leq 2$, $\frac{1}{p}+ \frac{1}{q}=1$, $\frac{1}{2} (n+1)\left( \frac{1}{p} - \frac{1}{q}\right) \leq 
2s \leq n \left( \frac{1}{p} - \frac{1}{q}\right)< 2s  +1 $.
\end{corollary}
\medskip
 
\noindent
{\bf Proof.}  
Indeed, we apply Lemma~\ref{L12.2}  with $z= e^{t -b}$ to the right-hand side of the estimate given by Theorem~\ref{T11.1} :
\begin{eqnarray*}
\| (-\bigtriangleup )^{-s} u(x,t) \|_{ { L}^{  q} ({\mathbb R}^n)  } 
& \le &  
C \int_{ 0}^{t} db\, \|  f(x, b)  \|_{ { L}^{p} ({\mathbb R}^n)  }   
e^{-t[2s-n(\frac{1}{p}-\frac{1}{q})]}    z^{-1}(z-1)^{ 1+  2s-n(\frac{1}{p}-\frac{1}{q}) } (1+ \ln z )\\
& \leq  &
C\int_{ 0}^{t} \|  f(x, b)  \|_{ { L}^{p} ({\mathbb R}^n)  }   
 e^{-b}\left( e^{-b}- e^{-t}  \right)^{ 1+  2s-n(\frac{1}{p}-\frac{1}{q}) } \left(1+  t-  b    \right)\,  db\,.
\end{eqnarray*}
Corollary is proven.
\hfill $\square$

\section{$L^p-L^q$ Estimates for  Equation without Source, $n\geq 2$}
\label{S13}
\setcounter{equation}{0}

The $L^p-L^q$-decay estimates for the energy of the solution of the Cauchy problem for the 
wave equation without source can be proved 
by the representation formula, $L_1-L_ \infty$ and  
$L_2-L_2$ estimates, and interpolation argument. (See, e.g., \cite[Theorem~2.1]{Racke}.) 
There is also a proof of the $L^p-L^q$-decay estimates that is based on the
microlocal consideration and dyadic decomposition of the phase space.  (See, e.g., \cite{Brenner,Pecher}.) 
To avoid the derivative loss and obtain more sharp estimates  we 
appeal to the representation formula provided by Theorem~\ref{T1.6}.

\begin{theorem} 
\label{T13.2}
The solution  $u=u(x,t)$ of the Cauchy problem (\ref{1.30CP})  satisfies the following $L^p-L^q$ estimate
\begin{eqnarray*} 
\| (-\bigtriangleup )^{-s} u(x,t) \|_{ { L}^{  q} ({\mathbb R}^n)  } 
& \leq & 
C  (1+ t  )(1- e^{-t}) ^{ 2s-n(\frac{1}{p}-\frac{1}{q}) }\Big\{ e ^{\frac{t}{2}}  \| \varphi_0  (x) \|_{ { L}^{p} ({\mathbb R}^n)  }
+ (1- e^{-t}) \|\varphi_1  \|_{ { L}^{p} ({\mathbb R}^n)  } 
 \Big\} 
\end{eqnarray*}
for all $t \in (0,\infty)$, provided that $s \ge 0$, $1<p \leq 2$, $\frac{1}{p}+ \frac{1}{q}=1$, $\frac{1}{2} (n+1)\left( \frac{1}{p} - \frac{1}{q}\right) \leq 
2s \leq n \left( \frac{1}{p} - \frac{1}{q}\right)  < 2s +1 $. 

\end{theorem}
\medskip

\noindent
{\bf Proof.} We start with the case of $\varphi _0=0$. Due  to Theorem~\ref{T1.6} for the solution $u=u (x,t)$ of the Cauchy 
problem (\ref{1.30CP}) with  $\varphi _0=0$ and to the results of \cite{Brenner,Pecher} 
we have:
\begin{eqnarray*}
&  &
\| (-\bigtriangleup )^{-s} u(x,t) \|_{ { L}^{  q} ({\mathbb R}^n)  } \\
& \leq &
C \|\varphi_1  \|_{ { L}^{p} ({\mathbb R}^n)  }\int_{0}^{1-e^{-t}}   
 r^{ 2s-n(\frac{1}{p}-\frac{1}{q}) }  \left| K_1(r,t)\right|   \, dr \\ 
& \leq &
C \|\varphi_1  \|_{ { L}^{p} ({\mathbb R}^n)  } e ^{ -t[2s-n(\frac{1}{p}-\frac{1}{q})] }\int_{0}^{e^{t } -1}   
y^{ 2s-n(\frac{1}{p}-\frac{1}{q}) } \big( (e^{t } +1)^2 -y^2 \big)^{-\frac{1}{2}} 
 F\left(\frac{1}{2} ,\frac{1}{2}  ;1; 
\frac{ (e^{t } -1)^2 -y^2 }{(e^{t } +1)^2 -y^2 } \right) \, dy \,.  
\end{eqnarray*}
To continue we apply Lemma~\ref{L12.2} and obtain 
\begin{eqnarray*}
\| (-\bigtriangleup )^{-s} u(x,t) \|_{ { L}^{  q} ({\mathbb R}^n)  } 
& \leq &
C \|\varphi_1  \|_{ { L}^{p} ({\mathbb R}^n)  } 
(1+ t)(1- e^{-t})^{ 1+  2s-n(\frac{1}{p}-\frac{1}{q})  }\,.  
\end{eqnarray*}
Thus, in the case of $\varphi _0=0$ the theorem is proven.
\medskip

Next we turn to the case of $\varphi _1=0$.    Due  to Theorem~\ref{T1.6} for the solution $u=u (x,t)$ of the Cauchy 
problem (\ref{1.30CP}) with  $\varphi _1=0$ and to the results of \cite{Brenner,Pecher} 
we have: 
\begin{eqnarray*} 
\| (-\bigtriangleup )^{-s} u(x,t) \|_{ { L}^{  q} ({\mathbb R}^n)  } 
& \leq & 
C \Bigg( e ^{\frac{t}{2}} (1- e^{-t}) ^{ 2s-n(\frac{1}{p}-\frac{1}{q}) } 
+ \int_{ 0}^{1- e^{-t}}  r^{ 2s-n(\frac{1}{p}-\frac{1}{q}) }
|K_0(r,t)|  \,  dr \Bigg) \| \varphi_0  (x) \|_{ { L}^{p} ({\mathbb R}^n)  }.
\end{eqnarray*}
One can estimate  the last integral  
\begin{eqnarray*} 
&  &
\int_{ 0}^{1- e^{-t}}  r^{ 2s-n(\frac{1}{p}-\frac{1}{q}) }
|K_0(r,t)|  \,  dr \\
& \leq & 
e^{ -t[2s-n(\frac{1}{p}-\frac{1}{q})] }\int_{ 0}^{e^{t }-1}  y^{ 2s-n(\frac{1}{p}-\frac{1}{q}) }
\frac{1}{ [(e^{t }-1)^2 -  y^2]\sqrt{(e^{t }+1)^2 - y^2} }\\
&   &
\times   \Bigg|  \big(  e^{t} - e^{2t} - iM(1 -     e^{ 2t} -  y^2) \big) 
F \Big(\frac{1}{2}+iM   ,\frac{1}{2}+iM  ;1; \frac{ ( e^{t }-1)^2 -y^2 }{( e^{t }+1)^2 -y^2 }\Big) \\
&  &
\hspace{1cm}  +   \big(  e^{2t} -1+  y^2 \big)\Big( \frac{1}{2}-iM\Big)
F \Big(-\frac{1}{2}+iM   ,\frac{1}{2}+iM  ;1; \frac{ ( e^{t }-1)^2 -y^2 }{( e^{t }+1)^2 -y^2 }\Big) \Bigg|  \,  d y \,.
\end{eqnarray*}
The following proposition gives the remaining   estimate for that integral
and completes the proof of the theorem. 
\begin{proposition}
\label{P13.4}
If  $2s-n(\frac{1}{p}-\frac{1}{q})>-1 $, then
\begin{eqnarray*} 
&  &
\int_{ 0}^{z-1}  y^{ 2s-n(\frac{1}{p}-\frac{1}{q}) }
\frac{1}{ [(z-1)^2 -  y^2]\sqrt{(z+1)^2 - y^2} }\\
&   &
\times   \Bigg|  \big( z - z^{2 } - iM(1 -     z^{ 2 } -  y^2) \big) 
F \Big(\frac{1}{2}+iM   ,\frac{1}{2}+iM  ;1; \frac{ ( z-1)^2 -y^2 }{( z+1)^2 -y^2 }\Big) \\
&  &
\hspace{1cm}  +   \big(  z^{2 } -1+  y^2 \big)\Big( \frac{1}{2}-iM\Big)
F \Big(-\frac{1}{2}+iM   ,\frac{1}{2}+iM  ;1; \frac{ ( z-1)^2 -y^2 }{( z+1)^2 -y^2 }\Big) \Bigg|  \,  d y \\
& \leq &
 C z^{ -\frac{1}{2}}(z-1)^{1+ 2s-n(\frac{1}{p}-\frac{1}{q}) }  \left(1+ \ln  z     \right)    \quad \mbox{ for all} \quad z>1.
\end{eqnarray*}
\end{proposition}
\medskip

\noindent
{\bf Proof.} We follow the arguments have been used in the proof of Proposition~\ref{P11.4}. If $1\leq  z \leq N$ with some constant $N$, then the argument  of the hypergeometric functions is bounded (\ref{argument}), and  the integral can be estimated by:
\begin{eqnarray*} 
&  &
\int_{ 0}^{z-1}  y^{ 2s-n(\frac{1}{p}-\frac{1}{q}) }
\frac{1}{ [(z-1)^2 -  y^2]\sqrt{(z+1)^2 - y^2} }\\
&   &
\times   \Bigg|  \big( z - z^{2 } - iM(1 -     z^{ 2 } -  y^2) \big) 
F \Big(\frac{1}{2}+iM   ,\frac{1}{2}+iM  ;1; \frac{ ( z-1)^2 -y^2 }{( z+1)^2 -y^2 }\Big) \\
&  &
\hspace{1cm}  +   \big(  z^{2 } -1+  y^2 \big)\Big( \frac{1}{2}-iM\Big)
F \Big(-\frac{1}{2}+iM   ,\frac{1}{2}+iM  ;1; \frac{ ( z-1)^2 -y^2 }{( z+1)^2 -y^2 }\Big) \Bigg|  \,  d y \\
& \leq &
C_M \int_{ 0}^{z-1}  y^{ 2s-n(\frac{1}{p}-\frac{1}{q}) }
\Bigg[ \frac{1}{  \sqrt{(z+1)^2 - y^2} }   
\Bigg\{ 1+ z^2\frac{1 }{(z +1)^2 -y^2 }   \Bigg\} 
 \Bigg]  \,  d y \\
& \leq &
C_Mz^{-1 }(z-1)^{ 1+2s-n(\frac{1}{p}-\frac{1}{q}) }  \quad \mbox{\rm for all}  \quad z \in [1,N]\,. 
\end{eqnarray*}
Thus, we can restrict ourselves to the  case of large $   z \geq M$ in both zones $Z_1(\varepsilon, z) $ and $ Z_2(\varepsilon, z )$, 
defined in  (\ref{9.10}) and (\ref{9.11}),
respectively. In the first zone we have (\ref{9.15}). 
Consider therefore the following inequalities,
\begin{eqnarray*}
A_{8}
& := &
\int_{ (z,r) \in Z_1(\varepsilon, z)  }  r^{ 2s-n(\frac{1}{p}-\frac{1}{q}) }  \frac{1}{     \sqrt{(z+1)^2-r^2}    }  \, dr \\
& \leq   &
C  z^{-1 }(z-1)^{ 1+2s-n(\frac{1}{p}-\frac{1}{q}) }  \quad \mbox{\rm for all}  \quad z \in [N,\infty)  \,.
\end{eqnarray*}
For $0 \geq a>-1 $ and $z\geq N$ the following integral can be easily estimated: 
\begin{eqnarray*} 
\int_{  0 }^{z-1}  r^{ a}    \frac{1}{ ((z+1)^2-r^2 )^{3/2}   }  dr 
& = &
\int_{  0 }^{z/2}  r^{ a}    \frac{1}{ ((z+1)^2-r^2 )^{3/2}   }  dr  + \int_{  z/2 }^{z-1}  r^{ a}    \frac{1}{ ((z+1)^2-r^2 )^{3/2}   }  dr  \\
& \leq  &
\frac{16}{9} z^{  -3}  \int_{  0 }^{z/2}  r^{ a }   dr  + \frac{z^{ a}}{4^{ a}} \int_{  z/2 }^{z-1}     \frac{1}{ ((z+1)^2-r^2 )^{3/2}   }  dr  \\
& \leq  &
C  z^{ a -3/2}  \quad \mbox{\rm for all}  \quad z \in [N,\infty)  \,.
\end{eqnarray*} 
Hence,
\begin{eqnarray*}
A_{9}
& := &
z^{2}\int_{  (z,r) \in Z_1(\varepsilon, z) } r^{ 2s-n(\frac{1}{p}-\frac{1}{q}) }  
 \frac{1}{   \sqrt{(z+1)^2-r^2}    } \frac{ 1}{ (z+1)^2 -r^2 }  
 dr \\
& \leq   &
z^{2}\int_{ 0}^{z-1} r^{ 2s-n(\frac{1}{p}-\frac{1}{q}) }  
 \frac{1}{   \sqrt{(z+1)^2-r^2}    } \frac{ 1}{ (z+1)^2 -r^2 }  
 dr \\
& \leq   &
 C  z^{ -\frac{1}{2}}(z-1)^{1+ 2s-n(\frac{1}{p}-\frac{1}{q}) }  
  \quad \mbox{\rm for all}  \quad z \in [N,\infty)\,,
\end{eqnarray*} 
and 
\begin{eqnarray*}
A_{10}
& := &
z^2\int_{  (z,r) \in Z_1(\varepsilon, z) }   r^{ 2s-n(\frac{1}{p}-\frac{1}{q}) } \frac{1}{ ((z-1)^2-r^2 ) \sqrt{(z+1)^2-r^2}    } \left(\frac{ (z-1)^2 -r^2   }{ (z+1)^2 -r^2 } \right)^2  dr \\
& \leq  &
z^2\int_{  (z,r) \in Z_1(\varepsilon, z) }  r^{ 2s-n(\frac{1}{p}-\frac{1}{q}) }  \frac{1}{  \sqrt{(z+1)^2-r^2}    } \frac{ 1   }{  (z+1)^2 -r^2  }   dr \\
& \leq  &
 C  z^{ -\frac{1}{2}}(z-1)^{1+ 2s-n(\frac{1}{p}-\frac{1}{q}) }  
  \quad \mbox{\rm for all}  \quad z \in [N,\infty) \,.
\end{eqnarray*}
Finally,
\begin{eqnarray*}
&  &
\int_{  (z,y) \in Z_1(\varepsilon, z) } \, y^{ 2s-n(\frac{1}{p}-\frac{1}{q}) }
\frac{1}{ [(z-1)^2 -  y^2]\sqrt{(z+1)^2 - y^2} }\\
&   &
\times   \Bigg|  \big( z - z^{2 } - iM(1 -     z^{ 2 } -  y^2) \big) 
F \Big(\frac{1}{2}+iM   ,\frac{1}{2}+iM  ;1; \frac{ ( z-1)^2 -y^2 }{( z+1)^2 -y^2 }\Big) \\
&  &
\hspace{1cm}  +   \big(  z^{2 } -1+  y^2 \big)\Big( \frac{1}{2}-iM\Big)
F \Big(-\frac{1}{2}+iM   ,\frac{1}{2}+iM  ;1; \frac{ ( z-1)^2 -y^2 }{( z+1)^2 -y^2 }\Big) \Bigg|  \,  d y \\& \leq  &
 C  z^{ -\frac{1}{2}}(z-1)^{1+ 2s-n(\frac{1}{p}-\frac{1}{q}) }  
  \quad \mbox{\rm for all}  \quad z \in [1,\infty) \,.
\end{eqnarray*}
In the second zone we use (\ref{10.4n}), (\ref{FGH+xat1}),   and  (\ref{10.8}). Thus, we have to estimate the next two integrals:
\begin{eqnarray*} 
A_{11} 
& :=   &
z^{2}\int_{  (z,r) \in Z_2(\varepsilon, z) }   r^{ 2s-n(\frac{1}{p}-\frac{1}{q}) } \frac{1}{     ((z-1)^2-r^2 ) \sqrt{(z+1)^2-r^2}    }  \,dr \,, \\
A_{12} 
& :=   & 
z^{2}\left(1+ \ln  z     \right) \int_{  (z,r) \in Z_2(\varepsilon, z) }   r^{ 2s-n(\frac{1}{p}-\frac{1}{q}) }  
\frac{1}{     ((z-1)^2-r^2 ) \sqrt{(z+1)^2-r^2}    } \,dr \,.
\end{eqnarray*}
We apply (\ref{10.4n}) to $A_{11}$ and obtain
\begin{eqnarray*} 
A_{11} 
 \leq   
C_\varepsilon z^{2}\int_{  (z,r) \in Z_2(\varepsilon, z) }   r^{ 2s-n(\frac{1}{p}-\frac{1}{q}) }   
\frac{ 1   }{ [(z+1)^2 -r^2] }\frac{1}{    \sqrt{(z+1)^2-r^2}    }  \, dr
 \leq  
C_\varepsilon z^{ -\frac{1}{2}}(z-1)^{1+ 2s-n(\frac{1}{p}-\frac{1}{q}) }    
\end{eqnarray*} 
for all $ z \in [1,\infty)$, while
\begin{eqnarray*} 
A_{12} 
& \leq    & 
C_\varepsilon z^{ -\frac{1}{2}}(z-1)^{1+ 2s-n(\frac{1}{p}-\frac{1}{q}) }  \left(1+ \ln  z     \right)  
\quad \mbox{\rm for all}  \quad z \in [1,\infty)  \,.
\end{eqnarray*} 
Proposition is proven.
\hfill $\square$
\smallskip

To complete the proof of the theorem we write
\begin{eqnarray*} 
&  &
\int_{ 0}^{1- e^{-t}}  r^{ 2s-n(\frac{1}{p}-\frac{1}{q}) }
|K_0(r,t)|  \,  dr \\
& \leq & 
e^{ -t[2s-n(\frac{1}{p}-\frac{1}{q})] }\int_{ 0}^{e^{t }-1}  y^{ 2s-n(\frac{1}{p}-\frac{1}{q}) }
\frac{1}{ [(e^{t }-1)^2 -  y^2]\sqrt{(e^{t }+1)^2 - y^2} }\\
&   &
\times   \Bigg|  \big(  e^{t} - e^{2t} - iM(1 -     e^{ 2t} -  y^2) \big) 
F \Big(\frac{1}{2}+iM   ,\frac{1}{2}+iM  ;1; \frac{ ( e^{t }-1)^2 -y^2 }{( e^{t }+1)^2 -y^2 }\Big) \\
&  &
\hspace{1cm}  +   \big(  e^{2t} -1+  y^2 \big)\Big( \frac{1}{2}-iM\Big)
F \Big(-\frac{1}{2}+iM   ,\frac{1}{2}+iM  ;1; \frac{ ( e^{t }-1)^2 -y^2 }{( e^{t }+1)^2 -y^2 }\Big) \Bigg|  \,  d y \\
& \leq & 
Ce^{ -t[\frac{1}{2} + 2s-n(\frac{1}{p}-\frac{1}{q})] } (e^t-1)^{1+ 2s-n(\frac{1}{p}-\frac{1}{q}) }  (1+ t  )\,.
\end{eqnarray*}
Thus,
\begin{eqnarray*} 
&  &
\| (-\bigtriangleup )^{-s} u(x,t) \|_{ { L}^{  q} ({\mathbb R}^n)  }  \\
& \leq & 
C \Big( e ^{\frac{t}{2}} (1- e^{-t}) ^{ 2s-n(\frac{1}{p}-\frac{1}{q}) } 
+ e^{ -t[\frac{1}{2} + 2s-n(\frac{1}{p}-\frac{1}{q})] } (e^t-1)^{1+ 2s-n(\frac{1}{p}-\frac{1}{q}) }  (1+ t  )\Big) \| \varphi_0  (x) \|_{ { L}^{p} ({\mathbb R}^n)  }\\
& \leq & 
C  (1+ t  )e ^{\frac{t}{2}} (1- e^{-t}) ^{ 2s-n(\frac{1}{p}-\frac{1}{q}) } \| \varphi_0  (x) \|_{ { L}^{p} ({\mathbb R}^n)  }.
\end{eqnarray*}
Theorem is proven. \hfill $\square$

\begin{small}

\end{small}
\end{document}